\documentclass{amsart}
\usepackage{amsmath,amsthm,amssymb,latexsym,stmaryrd,bbm}
\usepackage[e]{esvect}
\usepackage{mathrsfs}
\usepackage{enumerate, xcolor}
\usepackage[a4paper, margin=3cm]{geometry}
\usepackage[all,2cell]{xy}
\xyoption{2cell}

\usepackage[active]{srcltx}
\usepackage[parfill]{parskip}
\UseTwocells
\usepackage[
colorlinks=true,
linkcolor=black, 
anchorcolor=black,
citecolor=black, 
urlcolor=black, 
]{hyperref}

\usepackage{todonotes,comment}
\usepackage{tikz-cd}
\tikzcdset{ampersand replacement=\&}
\usepackage{quiver}
\usepackage[normalem]{ulem}
\newtheorem{theorem}{Theorem}[section]
\newtheorem{example}[theorem]{Example}
\newtheorem{prop}[theorem]{Proposition}
\newtheorem{lemma}[theorem]{Lemma}
\newtheorem{corollary}[theorem]{Corollary}
\newtheorem{conjecture}[theorem]{Conjecture}

\newtheorem*{theorem*}{Main Theorem}

\theoremstyle{definition}
\newtheorem{defn}[theorem]{Definition}

\theoremstyle{remark}
\newtheorem{remark}{Remark}[section]

\newcommand{\mdef}{\refstepcounter{theorem} 
	\medskip \noindent ({\bf \thetheorem}) }

\newcommand{\Ab}{\mathsf{Ab}}
\newcommand{\oneA}{\mathsf{A}}
\newcommand{\oneB}{\mathsf{B}}
\newcommand{\oneC}{\mathsf{C}}
\newcommand{\oneD}{\mathsf{D}}

\newcommand{\oneV}{\mathsf{V}}

\newcommand{\Vectk}{{\mathsf{Vect}}_{\mathbbm{k}}}

\newcommand{\Cat}{{\mathsf{Cat}}}

\newcommand{\F}{\mathrm{F}}
\newcommand{\G}{\mathrm{G}}

\newcommand{\K}{\mathrm{K}}
\def\H{\mathrm{H}}

\newcommand{\mon}{\otimes}
\newcommand{\monob}{\mon_0}
\newcommand{\monmor}{\mon_1}
\newcommand{\monnat}{\mon_2}

\newcommand{\desc}{\chi}
\newcommand{\what}{\underline{\otimes}}

\newcommand{\tenrep}{\circledcirc}
\newcommand{\tenrepC}{{\odot}}

\newcommand{\triv}{\bigstar}
\newcommand{\trivob}{\star}

\newcommand{\bfK}{\mathbf{K}}
\newcommand{\bfL}{\mathbf{L}}
\newcommand{\bfM}{\mathbf{M}}
\newcommand{\bfN}{\mathbf{N}}

\newcommand{\bfT}{\mathbf{T}}

\newcommand{\kk}{\mathbbm{k}}
\newcommand{\ob}[1]{\mathtt{#1}}

\newcommand{\A}{\mathfrak{A}_\kk}

\newcommand{\twoAk}{\mathscr{A}_{\mathbbm{k}}}

\newcommand{\Aa}{\mathcal{A}}
\newcommand{\Bb}{\mathcal{B}}
\newcommand{\Cc}{{\mathcal C}}
\newcommand{\Dd}{{\mathcal D}}
\newcommand{\Ee}{{\mathcal E}}
\newcommand{\Ff}{{\mathcal F}}
\newcommand{\Gg}{{\mathcal G}}
\newcommand{\Hh}{{\mathcal H}}

\newcommand{\Vv}{\mathcal{V}}

\newcommand{\hcirc}{\circ_{h}}

\newcommand{\bx}{\boxtimes}
\newcommand{\bxob}{\bx_0}
\newcommand{\bxone}{\bx_{1}}
\newcommand{\bxtwo}{\bx_{2}}

\newcommand{\add}{\operatorname{\mathsf{add}}}
\newcommand{\Add}{\operatorname{\mathsf{Add}}}

\newcommand{\amod}{\text{-}{\mathsf{amod}}}

\newcommand{\ass}{\mathnormal{a}}
\newcommand{\run}{\mathnormal{r}}
\newcommand{\lun}{\ell}

\newcommand{\asstwo}{\boldsymbol{\mathnormal{a}}}
\newcommand{\runtwo}{\boldsymbol{\mathnormal{r}}}
\newcommand{\luntwo}{\boldsymbol{\ell}}

\newcommand{\lambdatwo}{\boldsymbol{\lambda}}
\newcommand{\mutwo}{\boldsymbol{\mu}}

\newcommand{\rhotwo}{\boldsymbol{\rho}}

\newcommand{\assthree}{\underline{\boldsymbol{\mathnormal{a}}}}
\newcommand{\runthree}{\underline{\boldsymbol{\mathnormal{r}}}}
\newcommand{\lunthree}{\underline{\boldsymbol{\ell}}}

\newcommand{\lambdathree}{\underline{\boldsymbol{\lambda}}}
\newcommand{\muthree}{\underline{\boldsymbol{\mu}}}
\newcommand{\pithree}{\underline{\boldsymbol{\pi}}}
\newcommand{\rhothree}{\underline{\boldsymbol{\rho}}}

\newcommand{\One}{\mathbf{1}}

\newcommand{\ppmx}[1]{}

\newcommand{\ftx}[1]{}

\newcommand{\ftsout}{\bgroup\markoverwith{\textcolor{purple}{\rule[0.6ex]{2pt}{0.4pt}}}\ULon} 

\newcommand{\myaddress}%
{\parbox{3.3in}{\footnotesize \begin{center} 
			School of Mathematics, \\ University of East Anglia, 
			Norwich, NR4 7TJ, UK.\end{center}}}

\newcommand{\ignore}[1]{}

\newcommand{\Z}{\mathbb Z}      

\newcommand{\Ob}{\operatorname{Ob}}

\usepackage{mathtools}
\usepackage[foot]{amsaddr}

\title[Tensor products of $2$-representations]{Tensor products of $2$-representations}

\author{Vanessa Miemietz}
\address{V.M.: School of Engineering, Mathematics and Physics, University of East Anglia, Norwich NR4 7TJ, United Kingdom, {\tt  v.miemietz@uea.ac.uk}} 

\author{Fiona Torzewska}
\address{F.T.: School of Mathematics, University of Bristol, Heilbronn Institute of Mathematics}

\begin{document}

\begin{abstract}
In this paper we define comonoidal structures on additive $\Bbbk$-linear $2$-categories which guarantee the existence of a monoidal structure on the corresponding $2$-categories of $2$-representations.
\end{abstract}
\maketitle

\tableofcontents

\section*{Introduction}

In classical representation theory, there are two different types of tensor product being considered. On the one hand, given a ring or an algebra $A$ and a right, respectively left, $A$-module $M$, resp.\ $N$, we can consider their relative tensor product $M\otimes_A N$, which no longer carries an $A$-module structure unless one of $M$ and $N$ carries an $A$-$A$-bimodule structure. To ensure that every $A$-module is also an $A$-$A$-bimodule, it is necessary to require $A$ to be commutative. On the other hand, in certain areas of representation theory, like representation theory of groups or Lie algebras, one can take the external tensor product over the ground field $\Bbbk$ of two representations on $V$ resp.\ $W$ and obtain another representation on $V\otimes_\Bbbk W$. This is obtained via a bialgebra structure on the algebra that acts. Examples of such categories play an important role in mathematical physics, e.g.\ categories of small quantum groups at roots of unity provide examples of $3$-dimensional TQFTs, see e.g.\ \cite{RT}. More generally,
categories of representations with an external tensor product give rise to examples of tensor categories, see e.g.\ \cite[Chapter 5]{EGNO}, and have been intensely studied.

In higher representation theory, the two distinct tensor products have distinct generalisations. 
Categorifying the commutative approach is achieved by considering a Deligne tensor product of module categories over a braided base, and is developed for the semisimple case in \cite{DR}.

Many of the interesting categories with an external tensor product, such as representations of groups or quantum groups, have additive (non-semisimple) categorifications. A lot of work has been done on categorifying tensor products of specific $2$-representations for categorified quantum groups, with the eventual aim of categorifying $3$-dimensional TQFTs to $4$-dimensional ones following ideas by \cite{CF}. These approaches, however, do not use a categorified bialgebra structure, and are instead ad hoc, see e.g.\ \cite{SS, We, LW, Mc}. According to statements of Rouquier, a categorification of the Hopf algebra structure on quantum groups requires a setting of $A_\infty$-structures. Nevertheless, it seems valuable to axiomatise an appropriate categorification of a bialgebra that would allow a tensor product on categories of $2$-representations in general.

In this paper, we develop a rigorous, constructive framework for the purely additive categorification of the bialgebra approach. We establish the general algebraic machinery for comonoidal structures on locally additive $\Bbbk$-linear $2$-categories, providing a generalized method to tensor $2$-representations via a comonoidal pullback mechanism, without requiring the base $2$-category to be braided, semisimple, or restricted to a single object. Our main result is the following.

\begin{theorem*}
Let $\Cc$ be a locally additive $\Bbbk$-linear $2$-category equipped with a comonoidal structure. Then the $2$-category $\Cc\amod$ of strictly unital birepresentations of $\Cc$ has the structure of a monoidal $2$-category.
\end{theorem*}

We prove that, in particular, categorifications of groups do fit into our framework. It would be interesting to find other natural examples of categorifications of bialgebras in our sense.

The paper is organized as follows. In Section \ref{sec:mon2catof1cats}, we recall the foundational definitions of locally additive $\Bbbk$-linear categories and their external tensor product. In particular, we prove the folklore statement that locally additive $\Bbbk$-linear categories form a monoidal $2$-category. In Section  \ref{sec:add2cats}, we provide details for the external tensor product of locally $\Bbbk$-linear, locally additive $2$-categories. In particular, despite the fact that these $2$-categories do not assemble to a $2$-category, the ingredients for a monoidal structure such as associators, unitors, pentagonators and $2$-unitors all exist and satisfy the expected relations.
In Section \ref{sec:comon}, we formally define comonoidal $2$-categories, introducing the pseudofunctor $\Delta$, the counit $\Theta$, and the critical adjoint equivalences that enforce their coherence. Section \ref{sec:mon2catof2reps} contains the primary construction of the paper: defining the pullback tensor product $\bfM \tenrepC \bfN$ for $2$-representations and proving its monoidal coherence. Section \ref{sec:ex} provides concrete examples of this framework.

\section{Background}

\subsection{Notational conventions}
Let $\Bbbk$ be a field. By abuse of notation, we will also use $\Bbbk$ for the category with one object $\bullet_\Bbbk$, whose endomorphisms are given by $\Bbbk\cdot \mathbbm{1}_{\bullet_\Bbbk}$.

For a category $\oneC$ we generally use $\Ob(\oneC)$ to denote the collection of objects, and for $X,Y\in \Ob(\oneC)$ we use $\oneC(X,Y)$ to denote the collection of morphisms from $X$ to $Y$.
Depending on the context we may denote elements of $\oneC(X,Y)$ by only a symbol, $f$ for example, or by a triple $f\colon X \to Y$ indicating the source and target.
Our convention is that composition in $\oneC$ is a map
\begin{align*}
\circ \colon \oneC(Y,Z) \times \oneC(X,Y)&\mapsto \oneC(X,Z)\\
(g,f)&\mapsto g\circ f.
\end{align*}

We denote the category of finite dimensional vector spaces with linear maps by $\Vectk$. 


Given a monoidal category $K$, a category $\oneC$ {\em enriched} in $K$ consists of, a set $\Ob(\oneC)$,  for each pair of objects a $\oneC(X,Y)\in \Ob(K)$, a composition $\circ \colon\oneC(Y,Z)\otimes \oneC(X,Y)\to \oneC(X,Z)$ in $K$ and an identity in each $C(X,X)$, such that composition satisfies an appropriate notion of associativity and unitality.

A category enriched in the category of 
{$\Bbbk$-vector spaces} $\Vectk$ is called a {\it $\Bbbk$-linear category}.

Given two categories $\oneC$ and $\oneD$ enriched in the same monoidal category $\oneV$, a {\em $\oneV$-enriched functor} $\F\colon \oneC \to \oneD$ consists of a functor such that the function $\F_{X,Y} \colon \oneC(X,Y) \to \oneD(\F_0(X),\F_0(Y))$ is a morphism in $\oneV$ for each pair $(X,Y)\in \Ob(\oneC)\times \Ob(\oneC)$.
A functor enriched in 
{$\Bbbk$-vector spaces} $\Vectk$ is called a {\it $\Bbbk$-linear functor}.

A $2$-category is a category enriched in $\Cat$, where $\Cat$ denotes the category of small categories.
Here {\it $2$-category} means {\it strict $2$-category}. The weak notion is called a {\it bicategory}. Note any bicategory is biequivalent to a $2$-category, see for example \cite[2.3]{Leinster}. We will often use the convention that $2$-categories are denoted using calligraphic font as $\Cc$, $\Dd$, \ldots. 

We denote by $\triv$ the trivial $2$-category which has one object $\trivob$, one $1$-morphism $\mathbbm{1}_{\trivob}$ and one $2$-morphism $\mathbbm{1}_{\mathbbm{1}_{\trivob}}$.

\subsection{Background on $2$-categories}

\begin{defn}
Let $\Cc$ and $\Dd$ be $2$-categories. A {\em strictly unital pseudofunctor} $\Ff\colon \Cc\to \Dd$ consists of
\begin{itemize}
\item for any $\ob{i}\in \Ob{\Cc}$, an assignment of an object $\Ff(\ob{i})\in \Ob(\Dd)$;
\item for each pair of objects $\ob{i},\ob{j}\in \Cc$, a {\em local functor}
$\Ff=\Ff_{\ob{j},\ob{i}} \colon \Cc(\ob{i},\ob{j})\to\Dd(\Ff(\ob{i}),\Ff(\ob{j}))$
such that $\Ff_{\ob{i},\ob{i}}(\mathbbm{1}_{\ob{i}}) = \mathbbm{1}_{\Ff(\ob{i})}$
\item for each pair of composable $1$-morphisms $\F\in \Cc(\ob{j},\ob{k}),\G\in \Cc(\ob{i},\ob{j})$ a natural $2$-isomorphism
$$\Ff(\F)\Ff(\G) \xrightarrow{\Ff^2_{\F,\G}} \Ff(\F\G)$$
such that for any $\F\in \Cc(\ob{j},\ob{k}),\G\in \Cc(\ob{i},\ob{j}),\H\in\Cc(\ob{h},\ob{i})$
$$\xymatrix{
\Ff(\F)\Ff(\G)\Ff(\H)\ar^{\mathbbm{1}_{\Ff(\F)}\hcirc \Ff^2_{\G,\H}}[rr]  
\ar_{\Ff^2_{\F,\G}\hcirc\mathbbm{1}_{\Ff(\H)}}[d] 
&&\Ff(\F)\Ff(\G\H)\ar^{\Ff^2_{\F,\G\H}}[d]
\\
\Ff(\F\G)\Ff(\H)\ar_{\Ff^2_{\F\G,\H}}[rr] 
&&\Ff(\F\G\H)
}$$
commutes.
\end{itemize}
The natural isomorphism $\Ff^2$ is called {\em lax functoriality} and commutativity of the above diagram is called {\em lax associativity}.

A {\em $2$-functor} is a strictly unital pseudofunctor such that each $\Ff^2_{\F,\G}$ is an identity.
\end{defn}

\mdef
The composition of strictly unital pseudofunctors is again a strictly unital pseudofunctor, and moreover, this composition is strictly associative and unital by \cite[4.1.30]{JY}.

\mdef\label{strongnatother1}
We remark that the data of $\Ff^2$ can be viewed as a family of invertible natural transformations $\Ff_{\ob{k},\ob{j},\ob{i}} $
\[\begin{tikzcd}[ampersand replacement=\&]
	{\Cc(\ob{j},\ob{k}) \times \Cc(\ob{i},\ob{j})} \&\& {\Cc(\ob{i},\ob{k})} \\
	\\
	{\Dd(\Ff(\ob{j}),\Ff(\ob{k})) \times \Dd(\Ff(\ob{i}),\Ff(\ob{j}))} \&\& {\Dd(\Ff(\ob{i}),\Ff(\ob{k})).}
	\arrow["\circ", from=1-1, to=1-3]
	\arrow["{\Ff_{\ob{k},\ob{j}} \times \Ff_{\ob{j},\ob{i}} }"', from=1-1, to=3-1]
	\arrow["{\Ff_{\ob{k},\ob{i}}}", from=1-3, to=3-3]
	\arrow["{{\Ff^2_{\ob{k},\ob{j},\ob{i}}}}"', between={0.4}{0.6}, Rightarrow, from=3-1, to=1-3]
	\arrow["\circ"', from=3-1, to=3-3]
\end{tikzcd}\]


\begin{defn}\label{strongnat}
Let $\Cc$ and $\Dd$ be $2$-categories and $\Ff,\Gg\colon \Cc \to \Dd$ strictly unital pseudofunctors. A {\em strong natural transformation} $\Phi\colon \Ff\to\Gg$ consists of
\begin{itemize}
\item for each object $\ob{i}\in \Cc$,  a $1$-morphism $\Phi_{\ob{i}}\in \Dd(\Ff(\ob{i}), \Gg(\ob{i}))$
\item for each $\G\in  \Cc(\ob{i},\ob{j})$, a natural $2$-isomorphism $\Phi_{\G}\colon \Gg(\G)\Phi_{\ob{i}} \to \Phi_{\ob{j}}\Ff(\G)$
such that, for $\F\in \Cc(\ob{j},\ob{k})$, 
$$\Phi_{\F\G} (\Gg^2_{\F,\G}\hcirc \mathbbm{1}_{\Phi_{\ob{i}}}) =( \mathbbm{1}_{\Phi_{\ob{k}}}\hcirc\Ff^2_{\F,\G})(\Phi_{\F}\hcirc\mathbbm{1}_{\Ff(\G)})(\mathbbm{1}_{\Gg(\F)}\hcirc\Phi_{\G})$$
as morphisms $\Gg(\F)\Gg(\G)\Phi_{\ob{i}} \to \Phi_{\ob{k}}\Ff(\F\G)$.
\end{itemize}
The latter condition is called lax naturality.
\end{defn}

\mdef\label{strongnatother}
Let $\Cc$ be a $2$-category and $\K\in \Cc(\ob{j}, \ob{k})$. For any $\ob{i}\in \Ob(\Cc)$, there is a functor $\Cc(\ob{i},\ob{j})\to \Cc(\ob{i},\ob{k})$ which on objects acts as $\G\mapsto \K\circ G$. We denote this by $\mathnormal{L}_\K$.
Similarly, denote by $\mathnormal{R}_\K$ the corresponding functor which acts by pre-composition.
The collection of natural $2$-isomorphisms $\Phi_{\G}$ is the same as a collection of natural transformations of functors
$$\Phi_{\ob{j},\ob{i}} \colon \mathnormal{R}_{ \phi_{\ob{i}}}\Gg_{\ob{j}\ob{i}}  \to  \mathnormal{L}_{ \Phi_{\ob{j}}}  \Ff_{\ob{j}\ob{i}}$$
of functors  $\Cc(\ob{i},\ob{j})\to\Dd(\Ff(\ob{i}),\Gg(\ob{j}))$.

The following is Definition 4.2.15 from \cite{JY}, noting that in loc.\ cit.\ this is called horizontal composition, while we call it vertical composition.
\begin{defn}\label{de:compSNTs}
	Let $\Bb$ and $\Cc$ be bicategories, $\Ff,\Ff',\Ff''\colon \Bb\to \Cc$ be pseudofunctors, and $\alpha\colon \Ff\to \Ff'$, $\alpha'\colon \Ff'\to \Ff''$ be strong natural transformations. Then the \emph{vertical composition} of $\alpha$ and $\alpha'$ is the strong natural transformation $\alpha'\alpha\colon \Ff\to \Ff''$ where, for an object $\ob{i}\in \Ob(\Bb)$, $(\alpha'\alpha)_{\ob{i}} = \alpha'_{\ob{i}}\alpha_\ob{i}$ and for a morphism $f\in \Bb(\ob{i},\ob{j})$, $(\alpha'\alpha)_{f}$ is given by
	\[\Ff''(f)\alpha'_{\ob{i}}\alpha_{\ob{i}}\xrightarrow{\alpha'_{f}\star{\alpha_i}} \alpha'_{\ob{j}}\Ff'(f)\alpha_{\ob{i}}\xrightarrow{\alpha'_{\ob{j}}\star \alpha_{f}} \alpha'_{\ob{j}}\alpha_{\ob{j}}\Ff(f).\]	
\end{defn}

\mdef
It follows from \cite[Def.4.4.10]{JY} that in the case of $\Cc$ and $\Dd$ being $2$-categories, vertical composition of strong natural transformations is strictly associative and strictly unital.

\begin{defn}\label{modif}
Let $\Cc, \Dd$ be $2$-categories, $\Ff,\Gg\colon \Cc \to \Dd$ strictly unital pseudofunctors and $\Phi,\Psi\colon \Ff\to \Gg$ strong natural transformations. A {\em modification} $\alpha\colon \Phi\to \Psi$ consists of, for each object $\ob{i}\in \Cc$, a $2$-morphism
$$\alpha_{\ob{i}} \colon \Phi_{\ob{i}}\to\Psi_{\ob{i}}$$
such that, for any $\F\in  \Cc(\ob{i},\ob{j})$,
$$(\alpha_{\ob{j}} \hcirc \mathbbm{1}_{\Ff(\F)})\Phi_\F = \Psi_\F(\mathbbm{1}_{\Gg(\F)} \hcirc\alpha_{\ob{i}} )$$
as morphisms in $\Gg(\F)\Phi_i\to \Psi_j\Ff(\F)$.\\
The latter condition is called the modification axiom.
\end{defn}

The following is Definition 11.1.1 from \cite{JY}.
\begin{defn}\label{de:whiskering}
	Let $\Aa$, $\Bb$, $\Cc$ and $\Dd$ be bicategories,
	$\Ff\colon \Aa\to \Bb$, $\Gg,\Gg'\colon \Bb\to \Cc$ and $\Hh\colon \Cc \to \Dd$ be pseudofunctors, and $\alpha\colon \Gg\to \Gg'$ a strong natural transformation. Then the \emph{pre-whiskering} of $\alpha$ with $\Ff$ is defined as the strong natural transformation $\alpha \star \Ff\colon \Gg\Ff\to \Gg'\Ff$ where, for an object $\ob{i}\in \Ob(\Aa)$, $(\alpha\star\Ff)_{\ob{i}}$ is given by \[ \Gg\Ff(\ob{i})\xrightarrow{\alpha_{\Ff(\ob{i})}} \Gg'\Ff(\ob{i})\] 
	and for a morphism $f\in \Aa(\ob{i},\ob{j})$, $(\alpha\star F)_{f}$ is given by
\[\Gg'\Ff(f)\alpha_{\Ff(\ob{i})}\xrightarrow{\alpha_{F(f)}} \alpha_{\Ff(\ob{j})}\Gg\Ff(f).\]
The \emph{post-whiskering} of $\alpha$ with $\Hh$ is the strong natural transformation $\Hh \star \alpha\colon \Hh\Gg\to \Hh\Gg'$ where, for an object $\ob{i}'\in \Ob(\Bb)$, $(\Hh\star\alpha)_{\ob{i}'}$ is given by  
\[\Hh\Gg(\ob{i}')\xrightarrow{\Hh(\alpha_{\ob{i}'})} \Hh\Gg'(\ob{i}')\]
and $g\in \Bb(\ob{i}',\ob{j}')$, $(\Hh\star\alpha)_{g}$ is given by the composition
 \[\Hh\Gg'(g)\Hh(\alpha_{\ob{i}'})\xrightarrow{\Hh^2_{\Gg'(g),\alpha_{\ob{i}'}}}\Hh(\Gg'(g)\alpha_{\ob{i}'})\xrightarrow{\Hh(\alpha_{g})} \Hh(\alpha_{\ob{j}'}\Gg(g))\xrightarrow{(\Hh^2_{\alpha_{\ob{j}'},\Gg(g)})^{-1}}\Hh(\alpha_{\ob{j}'})\Hh\Gg(g).\] 
\end{defn}


	

The following lemmas will be needed later.

\begin{lemma}\label{le:whiskerassoc}
	Whiskering of strong natural transfromations is associative. Precisely,
let $\Aa$, $\Bb$, $\Cc$ and $\Dd$ be $2$-categories, $\Ff\colon \Aa\to \Bb$, $\Gg\colon \Bb\to \Cc$, $\Gg'\colon \Bb\to \Cc$ and $\Hh\colon \Cc \to \Dd$ be pseudofunctors, and $\alpha\colon \Gg\to \Gg'$ a strong natural transformation.
Then \[(\Hh \star \alpha)\star\Ff = \Hh \star (\alpha \star \Ff).\] 
\end{lemma}
\begin{proof}
Looking at Definition \ref{de:whiskering}, we see that, for both sides of the identity, we have that the component at $\ob{i}\in \Ob(\Aa)$ is given by
\[\Hh\Gg\Ff(\ob{i})\xrightarrow{\Hh(\alpha_{\Ff(\ob{i})})} \Hh\Gg'\Ff(\ob{i}),\]
and the component at $f\in \Aa(\ob{i},\ob{j})$ is given by the composition
\begin{multline*}\Hh\Gg'\Ff(f)\Hh(\alpha_{\Ff(\ob{i})})\xrightarrow{\Hh^2_{\Gg'\Ff(f),\alpha_{\Ff(\ob{i})}}}\Hh(\Gg'\Ff(f)\alpha_{\Ff(\ob{i})})\xrightarrow{\Hh(\alpha_{\Ff(f)})} \Hh(\alpha_{\Ff(\ob{j})}\Gg\Ff(f))\\
	\xrightarrow{(\Hh^2_{\alpha_{\Ff(\ob{j})},\Gg\Ff(f)})^{-1}}\Hh(\alpha_{\Ff(\ob{j})})\Hh\Gg(f).\end{multline*}
\end{proof}

\begin{lemma}\label{le:partialinterchange}
Let $\Aa$, $\Bb$, $\Cc$ and $\Dd$ be $2$-categories, $\Ee\colon \Aa\to \Bb$ and $\Ff,\Ff',\Ff''\colon \Bb\to \Cc$ be pseudofunctors, $\Gg\colon \Cc\to \Dd$ be a $2$-functor, and $\alpha\colon \Ff\to \Ff'$, $\alpha'\colon \Ff'\to \Ff''$ be strong natural transformations. This data is represented by the following diagram.
\[\begin{tikzcd}[ampersand replacement=\&]
	\Aa \& \Bb \&\& \Cc \& \Dd
	\arrow["{\mathcal{E}}", from=1-1, to=1-2]
	\arrow[""{name=0, anchor=center, inner sep=0}, "\Ff", curve={height=-24pt}, from=1-2, to=1-4]
	\arrow[""{name=1, anchor=center, inner sep=0}, "{\Ff''}"', curve={height=24pt}, from=1-2, to=1-4]
	\arrow[""{name=2, anchor=center, inner sep=0}, "{\Ff'}"'{pos=0.3}, from=1-2, to=1-4]
	\arrow["\Gg", from=1-4, to=1-5]
	\arrow["\alpha", between={0.3}{0.7}, Rightarrow, from=0, to=2]
	\arrow["{\alpha'}", between={0.3}{0.7}, Rightarrow, from=2, to=1]
\end{tikzcd}\]
Then
\begin{itemize}
	\item[(i)] $ (\Gg\star\alpha')(\Gg\star\alpha)= \Gg\star(\alpha'\alpha)$, and
	\item[(ii)] $(\alpha'\star\Ee)(\alpha\star\Ee) = (\alpha'\alpha)\star\Ee$.
\end{itemize}
\end{lemma}
\begin{proof}
	(i) Applying definitions Definitions \ref{de:compSNTs} and \ref{de:whiskering} for the simplified setting where all categories are $2$-categories and $\Gg$ is a $2$-functor, we see that, for both sides of the identity, we have that the component at $\ob{i}\in \Ob(\Bb)$ is given by
	\[\Gg\Ff(\ob{i})\xrightarrow{\Gg(\alpha'_{\ob{i}}\alpha_{\ob{i}})} \Gg\Ff''(\ob{i})\]
	and, for a morphism $f\in \Bb(\ob{i},\ob{j})$
	we have that the component at $f$ is given by the composition
	\[
	\Gg\Ff''(f)\Gg(\alpha'_{\ob{i}}\alpha_{\ob{i}})\xrightarrow{\Gg((\mathbbm{1}_{\alpha'_{\ob{j}}}\hcirc \alpha_{f})(\alpha'_{f}\hcirc\mathbbm{1}_{\alpha_i}) )} \Gg(\alpha'_{\ob{j}}\alpha_{\ob{j}})(\Ff(f)),
	\]
	where the given morphism is a morphism between the given source and target because $\Gg$ is a $2$-functor.\\
	(ii) Again applying Definitions \ref{de:compSNTs} and \ref{de:whiskering}, 
	we have that the component at $\ob{i}\in \Ob(\Aa)$ is given by
	\[\Ff\Ee(\ob{i})\xrightarrow{\alpha'_{\Ee(\ob{i})}\alpha_{\Ee(\ob{i})}} \Ff''\Ee(\ob{i})\]
	and, for a morphism $f\in \Aa(\ob{i},\ob{j})$
	we have that the component at $f$ is given by the composition
	\[
	\Ff''\Ee(f)\alpha'_{\Ee(\ob{i})}\alpha_{\Ee(\ob{i})}\xrightarrow{(\mathbbm{1}_{\alpha'_{\Ee(\ob{j})}}\hcirc \alpha_{\Ee(f)})(\alpha'_{\Ee(f)}\hcirc\mathbbm{1}_{\alpha_{\Ee(\ob{i})}})} \alpha'_{\Ee(\ob{j})}\alpha_{\Ee(\ob{j})}\Ff\Ee(f).
	\]
\end{proof}

\begin{defn}\label{de:Bicat}
Given $2$-categories $\Cc,\Dd$, an {\em adjoint equivalence} 
between two strictly unital pseudofunctors $\Ff,\Gg\colon \Cc \to \Dd$ 
is an internal adjoint equivalence in the $2$-category whose objects are $\Ff$, and $\Gg$, whose $1$-cells are strong natural transformations and whose $2$-cells are modifications, see \cite[Def.6.1.1]{JY}. In other words, an adjoint equivalence
consists of a tuple $(\Phi,\Phi^*, \eta^\Phi, \varepsilon^\Phi)$ where $\Ff \underset{\Phi^*}{\overset{\Phi}{\rightleftarrows} }\Gg$ are strong natural transformations, and $ \eta^\Phi\colon \mathbbm{1}_{\Ff}\to \Phi^*\Phi$, $\varepsilon^\Phi\colon \Phi\Phi^*\to  \mathbbm{1}_{\Gg}$ are invertible modifications satisfying the triangle equality for adjunctions, i.e.\
\[\xymatrix{
\Phi\ar^{\mathbbm{1}_{\Phi}\hcirc \eta^\Phi}[rr]\ar@{=}[drr]&&\Phi\Phi^*\Phi\ar^{\varepsilon^\Phi\hcirc \mathbbm{1}_{\Phi}}[d]\\
&&\Phi
}
\xymatrix{
\Phi^*\ar^{ \eta^\Phi\hcirc \mathbbm{1}_{\Phi^*}}[rr]\ar@{=}[drr]&&\Phi^*\Phi\Phi^*\ar^{\mathbbm{1}_{\Phi^*}\hcirc\varepsilon^\Phi}[d]\\
&&\Phi^*
}.
\]
The modifications $\eta^\Phi$ and $\varepsilon^\Phi$ are also called the unit and counit, respectively, of the adjunction.
\end{defn}


\mdef\label{pa:datamon2cat}In Sections  \ref{sec:mon2catof1cats} and \ref{sec:mon2catof2reps} below we will give two explicit constructions of monoidal $2$-categories, following \cite[12.1.3]{JY}. In particular we will give tuples $(\Cc,\otimes, \boldsymbol{1}, a, l, r, \pi, \mu, \lambda, \rho)$ consisting of the following, which in each case are examples of the data required by \cite{JY}.
	The input category $\Cc$ will always be a $2$-category, and $\otimes\colon \Cc\times \Cc \to \Cc$ and $\boldsymbol{1}\colon \triv \to \Cc$ will always be $2$-functors.
	For us, $a$ will be an adjoint equivalence between $\otimes(\otimes \times \mathbbm{1})$ and $\otimes(\mathbbm{1}\times \otimes)$, and  $l$, resp.\ $r$, will be adjoint equivalences between $\otimes (\boldsymbol{1}\times \mathbbm{1})$ and $\mathbbm{1}$, resp.\ $\otimes(\mathbbm{1}\times\boldsymbol{1})$ and $\mathbbm{1}$. 
The $2$-unitors  $\mu$, $\lambda$, $\rho$ will be modifications  {with component $2$-cells 
\[(\mathbbm{1}_C\otimes l_D)\, a_{C,\boldsymbol{1}(\trivob), D } \,(\run_C^*\otimes\mathbbm{1}_D )\xrightarrow{\mu_{C,D}}\mathbbm{1}_{C\otimes D}\]
\[l_C\otimes \mathbbm{1}_D\xrightarrow{\lambda_{C,D}}l_{C\otimes D}\, a_{\boldsymbol{1}(\trivob), C, D} \]
\[\mathbbm{1}_C\otimes r_D^* \xrightarrow{\rho_{C,D}} a_{C,D,\boldsymbol{1}(\trivob)} \,r^*_{C\otimes D}\]
for objects $C, D\in \Cc$}
and the $2$-associator $\pi$ will always be an identity modification with component identity $2$-cells 
\[\pi\colon (\mathbbm{1}_{A}\otimes a_{B, C, D}) a_{A,B\otimes C, D}( a_{A, B, C}\otimes\mathbbm{1}_{D}) = a_{A, B, C\otimes D}  a_{A\otimes B, C, D} \] 
for objects $A,B,C,D\in \Cc$.

\subsection{Additive categories}

We note some relevant properties and notational conventions on additive categories here, for a more complete introductory reference see \cite{Ritter} for example.

\mdef
We denote by $\Ab$ the monoidal category whose objects are abelian groups, whose morphisms are group homomorphisms and whose monoidal product is given by the tensor product over $\Z$ when viewing the objects as $\Z$-modules. Thus, the monoidal unit is $\Z$.

\mdef\label{pa:defaddcat} Recall that a category is  {\em preadditive} if it is enriched in $\Ab$. It is called  {\em additive} if it further admits all finite coproducts. We use {\em (pre)additive functor} to mean an $\Ab$-enriched functor between (pre)additive categories.

\mdef\label{pa:addnote} We note some well-known properties of additive categories that we will make use of later.
Firstly, all finite coproducts exist, where finite coproduct means colimit over any finite set, including, in particular, the empty set. This implies the existence of an initial object.
Secondly, recall that in an additive category, the product of any finite set of objects also exists and is naturally isomorphic to the coproduct. We denote a biproduct of $X$ and $Y$ by $X\oplus Y$. We denote the corresponding inclusion maps by $\iota_X, \iota_Y$ and the projections by $p_X, p_Y$.
In general, we use matrices of morphisms to describe morphisms between direct sums, i.e.\ for a morphism $f\colon X_1\oplus\cdots \oplus X_n \to Y_1\oplus \cdots \oplus Y_m$, we write $f=(f_{ij})$ where $f_{ij} = p_{Y_i}f\iota_{X_j}$.   Using this notation we may observe that the additivity of a functor $\F\colon \oneC\to \oneD$ between additive categories implies that, for any $X,Y\in \oneC$, there is an isomorphism $\F(X\oplus Y)\cong \F(X)\oplus\F(Y)$.

\mdef The following constructions are treated carefully in \cite[Sec.2.2]{Ritter}. We note that it is an abuse of notation to denote the formal tuples below using the same symbol used to denote direct sums above. By construction these formal tuples are direct sums, and thus where direct sums already exist, the two objects are isomorphic.

\mdef\label{pa:addclo} Given a preadditive category $\oneC$, one can construct an additive category $\Add \oneC$, with tuples of objects in $\oneC$ as objects and formal matrices whose entries are morphisms of A as morphisms. Our convention is to write object tuples in $\Add\oneC$ as $X=X_1\oplus \ldots \oplus X_n$ where each $X_i\in \Ob(C)$. To ease notation, we also write $\bigoplus_{i=1}^n X_i$ for the same (ordered) tuple. Composition is then given by usual matrix multiplication, together with componentwise application of the composition in $\oneC$. The category $\Add \oneC$ has the property that it comes with the natural inclusion $\iota_\oneC\colon\oneC \to \Add\oneC$ and that, for any preadditive functor $\mathrm{F}$ from $\oneC$ to an additive category $\oneD$, there is a unique (up to isomorphism) additive functor, which we also denote by $\mathrm{F}$ (rather than $\Add \F$), from  $\Add \oneC$ to $\oneD$ making the diagram
$$\xymatrix{
 \oneC\ar^{\mathrm{F}}[drr]\ar@{^{(}->}^{\iota_\oneC}[rr]&&\Add \oneC\ar@{^{}-->}^{\mathrm{F}}[d]\\
 && \oneD
}$$
commute which acts strictly with respect to the formal biproducts in $\Add \oneC$. 
Observe that if $\F$ is a $\Bbbk$-linear functor between $\Bbbk$-linear categories, the induced functor is again $\Bbbk$-linear.

\mdef\label{pa:AddF} In general, given preadditive categories $\oneC,\oneD$ and a preadditive functor $\F\colon \oneC\to \oneD$, applying the above construction to $\iota_{\oneD}\F \colon\oneC\to\Add\oneD$ induces an additive functor 
$$\Add\F\colon\Add\oneC\to \Add\oneD.$$

\mdef\label{pa:i*} Notice that, if $\oneC$ is additive, the inclusion $ \iota_{\oneC}\colon\oneC\to\Add \oneC$ is an equivalence. We denote by $\iota_{\oneC}^*$ the canonical inverse equivalence (which is bidajoint to $\iota_{\oneC}$). This is given by sending the object $(X_{11}\oplus\cdots\oplus X_{1r_1})\oplus\cdots \oplus (X_{s1}\oplus\cdots\oplus X_{sr_s})$ in $\Add\oneC$ to $X_{11}\oplus\cdots\oplus X_{1r_1}\oplus\cdots \oplus X_{s1}\oplus\cdots\oplus X_{sr_s}$ in $\oneC$.

\mdef\label{pa:indnattrans} Given a preadditive category $\oneC$, an additive category $\oneD$, two preadditive functors $\mathrm{F},\mathrm{G}\colon \oneC\to \oneD$ and a natural transformation $\alpha\colon \mathrm{F}\to\mathrm{G}$, there is an induced natural transformation between the 
induced additive functors $\mathrm{F},\mathrm{G}\colon \Add(\oneC)\to \oneD$, whose component at an object $X_1\oplus\cdots \oplus X_n$ is just given by the diagonal matrix $\operatorname{diag}(\alpha_{X_1}, \dots,\alpha_{X_n} )$.

\mdef Given a full subcategory $\oneD$ of an additive category $\oneC$, we denote by $\add \oneD$ the additive closure of $\oneD$ inside $\oneC$, meaning the closure under direct sums, direct summands and isomorphisms. Note that $\add\oneD$ and $\Add\oneD$ are naturally equivalent, provided $\oneD$ is already closed under direct summands in $\oneC$, in particular, if it is idempotent complete.
\section{Monoidal structure on the $2$-category of additive  $\Bbbk$-linear categories} \label{sec:mon2catof1cats}


The main result of this section is Theorem~\ref{th:mon2catof1cats}
which proves the existence of a monoidal $2$-category with base $2$-category $\A$ of $\Bbbk$-linear additive categories (constructed in \eqref{pa:catA}). 
In Section \ref{sec:tensorprod_kfincats} we give a mapping $\A\times \A\to \A$ and in 
Section \ref{ss:mon2catof1cats}, we confirm that this leads to a $2$-functor, give all other required data as listed in \eqref{pa:datamon2cat}, and prove that this yields a monoidal $2$-category.

\subsection{Data involved in the tensor product of additive  $\Bbbk$-linear categories}\label{sec:tensorprod_kfincats}\hfill\noindent 
\vspace*{-\baselineskip}

\mdef \label{pa:catA}
Let $\A$ denote the $2$-category whose
\begin{itemize}
\item objects are small $\Bbbk$-linear additive categories;
\item morphism categories are categories of $\Bbbk$-linear (additive) functors and natural transformations.
\end{itemize}
Here, horizontal composition is given by composition of functors and identities are the identity functors.

\mdef\label{naiveotimes} Given a pair of $\Bbbk$-linear categories $\oneC$ and $\oneD$, there exists a $\Bbbk$-linear category $\oneC\diamond \oneD$ whose objects are pairs of objects $(X,Y)\in \Ob(\oneC)\times \Ob( \oneD)$ and whose morphism spaces are given by 
$$(\oneC\diamond\oneD)((X,Y), (X',Y')) = \oneC(X,X')\otimes_\Bbbk \oneD(Y,Y')$$
with component-wise composition. We call this the naive tensor product of $\oneC$ and $\oneD$.
Observe that it comes with a canonical functor $$\oneC\times \oneD \to  \oneC\diamond \oneD$$  given by the identity on objects and the canonical map coming from the universal property of tensor products on morphism spaces.

Given $\Bbbk$-linear categories $\oneC,\oneC',\oneD,\oneD'$ and $\Bbbk$-linear functors $\F\colon\oneC\to \oneC', \G\colon\oneD\to \oneD'$, we define a functor 
$$\F\diamond \G \colon \oneC\diamond \oneD\to \oneC'\diamond \oneD'$$ 
given by component-wise application of $\F$ and $\G$, respectively.

Given four $\Bbbk$-linear additive categories $\oneC,\oneC',\oneD,\oneD'$, four $\Bbbk$-linear functors $\F,\F'\colon\oneC\to \oneC', \G,\G'\colon\oneD\to \oneD'$, and natural transformations $\alpha\colon \F\to\F', \beta\colon \G\to \G'$, there is an induced natural transformation
$$\alpha\diamond\beta\colon \F\diamond\G\to \F'\diamond\G'$$
whose component morphism on an object $(X,Y)$ is given by $\alpha_{X}\otimes \beta_{Y}$.

\mdef\label{otimesdef}
If $\oneC$ and $\oneD$ are $\Bbbk$-linear additive categories, we define their tensor product 
$$\oneC\monob \oneD:=\Add\left(\oneC\diamond \oneD\right),$$
noting that this is by construction a $\Bbbk$-linear additive category.


Given $\Bbbk$-linear additive categories $\oneC,\oneC',\oneD,\oneD'$ and $\Bbbk$-linear functors $\F\colon\oneC\to \oneC', \G\colon\oneD\to \oneD'$, we define a functor $\F\monmor\G$ determined by the diagram
$$\xymatrix{
\oneC\diamond \oneD\ar^{\iota_{\oneC\diamond \oneD}}[rr]\ar_{\F\diamond \G}[d]&&\oneC\monob \oneD\ar^{\F\monmor\G}[d] \\
    \oneC'\diamond \oneD'\ar^{\iota_{\oneC'\diamond \oneD'}}[rr]&&\oneC'\monob \oneD'
}$$
applying the universal property of the additive closure (see \eqref{pa:addclo}) to the composition $\iota_{(\oneC'\diamond \oneD')}\circ (\F\diamond \G)$. In other words, $\F\monmor\G = \Add  (\F\diamond \G)$, as in (\ref{pa:AddF}). This is by construction a $\Bbbk$-linear  functor.

Likewise, {extending (\ref{pa:indnattrans})}, given four $\Bbbk$-linear additive categories $\oneC,\oneC',\oneD,\oneD'$, four $\Bbbk$-linear functors $\F,\F'\colon\oneC\to \oneC', \G,\G'\colon\oneD\to \oneD'$, and natural transformations $\alpha\colon \F\to\F', \beta\colon \G\to \G'$, there is an induced natural transformation
$$\alpha\monnat\beta\colon \F\monmor\G\to \F'\monmor\G'$$
whose component morphism on an object $(X_1,Y_1)\oplus\cdots\oplus (X_n,Y_n)$ is given by the diagonal matrix 
$\operatorname{diag}((\alpha\diamond\beta)_{(X_1,Y_1)}, \dots, (\alpha\diamond\beta)_{(X_n,Y_n)})=\operatorname{diag}(\alpha_{X_1}\otimes \beta_{Y_1}, \dots, \alpha_{X_n}\otimes \beta_{Y_n})$.

\mdef\label{addadj} We now consider some natural morphisms that will be important in later constructions.
Given $\Bbbk$-linear categories $\oneC, \oneD$ and considering the diagram
$$\xymatrix{
\oneC\diamond \oneD \ar^{\iota_{\oneC\diamond \oneD}}[rr] 
\ar^{\iota_{\oneC}\diamond \iota_{\oneD}}[dr]
&& \Add (\oneC \diamond \oneD)\ar@{-->}^{\Add (\iota_{\oneC}\diamond \iota_{\oneD})}[dd]\\
& \Add \oneC\diamond \Add\oneD  \ar^{\iota_{\Add \oneC\diamond \Add\oneD }}[dr]&\\
&&\Add (\Add \oneC\diamond \Add\oneD),
}$$
it is easy to see that the induced, as in (\ref{pa:AddF}), vertical functor $\Add (\iota_{\oneC}\diamond \iota_{\oneD})$ is an equivalence. In particular, if $\oneC$ and  $\oneD$ are already additive, a biadjoint inverse equivalence is given by $\Add (\iota_{\oneC}^*\diamond \iota_{\oneD}^*)$, where $\iota_{\oneC}^*$ is as in (\ref{pa:i*}).
Similarly, if  $\oneC$ is additive, 
$$\Add (\mathbbm{1}_{\oneC}\diamond\iota_{\oneD})\colon \Add(\oneC\diamond\oneD) \to \Add(\oneC\diamond\Add\oneD)$$ is an equivalence, and if $\oneD$ is additive, 
$$\Add (\iota_{\oneC}\diamond\mathbbm{1}_{\oneD} ) \colon \Add(\oneC\diamond\oneD)\to \Add(\Add\oneC\diamond\oneD)$$ 
is an equivalence, with biadjoint inverse equivalences $(\Add (\mathbbm{1}_{\oneC}\diamond \iota_{\oneD}))^*$ and $(\Add (\iota_{\oneC}\diamond\mathbbm{1}_{\oneD} ))^*$ induced by
$$ \Add(\oneC\diamond\Add\oneD)\to \Add(\oneC\diamond\oneD) , \; (X, Y_1\oplus \cdots \oplus Y_n)\mapsto (X, Y_1)\oplus \cdots \oplus (X, Y_n)$$
and $$\Add(\Add\oneC\diamond\oneD) \to \Add(\oneC\diamond\oneD), \;  (X_1\oplus\cdots\oplus X_n, Y)\mapsto (X_1,Y)\oplus\cdots\oplus (X_n, Y),$$
respectively. Notice these strictly commute with any functors induced on the additive closure as in \ref{pa:addclo}.

The inverse equivalences we have constructed will be sufficient to prove that $\monob, \monmor, \monnat$ leads to a monoidal structure. In Sec.\ref{sec:add2cats}, we will need the general case, so we also give this here.
More generally, to define the (bi)adjoint functor to  $\Add (\iota_{\oneC}\diamond \iota_{\oneD})$, we choose the one given by the reverse lexicographic ordering, that is, the one induced by 
\begin{align*}
\Add(\Add\oneC\diamond\Add\oneD) & \to \Add(\oneC\diamond\oneD), \\  (X_1\oplus\dotsi\oplus X_n, Y_1\oplus\dotsi \oplus Y_m)&\mapsto (X_1,Y_1)\oplus\dotsi\oplus (X_1,Y_m)\oplus (X_2,Y_1)\oplus\dotsi\oplus (X_n, Y_m).
\end{align*}

\mdef\label{addadjmor}
Observe that if $\oneC$ is additive, the adjunction morphisms 
$$\mathbbm{1}_{ \Add(\oneC\diamond\oneD)} \leftrightarrows (\Add (\mathbbm{1}_{\oneC}\diamond \iota_{\oneD}))^*\Add (\mathbbm{1}_{\oneC}\diamond \iota_{\oneD})$$
are the identities, while the adjunction morphisms
$$\mathbbm{1}_{\Add(\oneC\diamond\Add\oneD)} \leftrightarrows (\Add (\mathbbm{1}_{\oneC}\diamond \iota_{\oneD}))(\Add (\mathbbm{1}_{\oneC}\diamond \iota_{\oneD}))^*$$
are given by the canonical inverse isomorphisms
$$\bigoplus_i(X_i, \bigoplus_j Y_{ij}) \leftrightarrows \bigoplus_i\bigoplus_j(X_i,  Y_{ij}),$$
and similarly in the case where $\oneD$ is additive.

\subsection{The monoidal $2$-category $\A$}
\label{ss:mon2catof1cats}

In this subsection, we confirm that the tensor products of additive categories, functors and natural transformations above endow $\A$ with the structure of a monoidal $2$-category. 

\begin{prop}\label{tensorone}
The functions $\monob, \monmor, \monnat$ on the levels of objects, $1$-morphisms and $2$-morphisms define a $2$-functor $\otimes \colon\A\times\A\to\A$.
\end{prop}

\proof
It is straightforward to check that $\otimes$ restricts, for all quadruples of categories $\oneC,\oneC',\oneD,\oneD'\in \A$, to a functor
\[\otimes \colon \A(\oneC,\oneC')\times \A(\oneD,\oneD')\to \A(\oneC\monob \oneC',\oneD\monob \oneD').
\]
To check that this is a $2$-functor, we need to verify that for $\Bbbk$-linear functors $\F\in \A(\oneB,\oneC), \F'\in \A(\oneB',\oneC'),\G\in \A(\oneA,\oneB)\G'\in \A(\oneA',\oneB')$, we have
$$\F\G\monmor \F'\G' = (\F\monmor\F')  (\G\monmor\G')$$
and that $\mathbbm{1}_{\oneA\monob\oneB} = \mathbbm{1}_{\oneA}\monmor\mathbbm{1}_{\oneB}$. 
The first condition follows from the fact that the two sides are induced by the diagrams
$$ 
\xymatrix{
(\oneA\diamond \oneA') \ar^{\iota_{(\oneA\diamond\oneA')}}[rr]\ar^{(\F\G\diamond \F'\G')}[d]
&&\oneA\monob \oneA'\ar^{\F\G\monmor \F'\G'}[d]\\
(\oneC\diamond \oneC') \ar^{\iota_{(\oneC\diamond \oneC')}}[rr] 
&&\oneC\monob \oneC'
}\;\;
\xymatrix{
(\oneA\diamond \oneA') \ar^{\iota_{(\oneA\diamond \oneA')}}[rr]\ar^{ (\G\diamond \G')}[d]
&&\oneA\monob \oneA' \ar^{ (\G\monmor\G')}[d]\\
(\oneB\diamond \oneB') \ar^{\iota_{(\oneB\diamond\oneB')}}[rr]\ar^{(\F\diamond \F') }[d]
&&\oneA\monob \oneA' \ar^{(\F\monmor\F')}[d]\\
(\oneC\diamond \oneC') \ar^{\iota_{(\oneC\diamond\oneC')}}[rr] 
&&\oneC\monob \oneC'\\
}$$
and, by construction, $(\F\G\diamond\F'\G') = (\F\diamond \F')(\G\diamond \G')$. 
Likewise $\mathbbm{1}_{\oneA}\monmor\mathbbm{1}_{\oneB}$ is induced by the diagram 
$$ 
\xymatrix{
(\oneA\diamond \oneB) \ar^{\iota_{(\oneA\diamond \oneB)}}[rr]\ar_{(\mathbbm{1}_{\oneA}\diamond \mathbbm{1}_{\oneB})}[d]
&&\oneA\monob \oneB\ar^{\mathbbm{1}_{\oneA}\monmor\mathbbm{1}_{\oneB}}[d]\\
(\oneA\diamond \oneB)  \ar^{\iota_{(\oneA\diamond \oneB)}}[rr] 
&&\oneA\monob \oneB
}$$
and is equal to $\mathbbm{1}_{\oneA\monob\oneB}$ since the latter also makes this diagram commute.
\endproof

\begin{prop}\label{unitone}
There is a $2$-functor 
\begin{align*}\One_{\A}\colon\triv \to \A \colon\qquad \trivob &\mapsto \Add \Bbbk \\
 \mathbbm{1}_\trivob &\mapsto \mathbbm{1}_{\Add \Bbbk}\\
\mathbbm{1}_{\mathbbm{1}_\trivob} &\mapsto \mathbbm{1}_{\mathbbm{1}_{\Add \Bbbk}}.
\end{align*}
\end{prop}
\proof
This is trivially true, since the composition of identity functors is the identity functor.
\endproof

In Proposition~\ref{assone} and \ref{unitorone} below, we give the three adjoint equivalences as described in \eqref{pa:datamon2cat}. Looking at Definition \ref{strongnat} notice that for any $2$-category $\Cc$, the data of a strong natural transformation between $2$-functors from $\Cc$ to $\A$ is given by specifying, for each $\ob{i}\in \Cc$, an additive functor. The natural transformations we give will all be strict, so this is the only data.
The data of a modification is then given by specifying a family of natural transformations associated to each $\ob{i}\in \Cc$.
Further, by definition, an adjoint equivalence between $2$-functors restricts, for each object $\ob{i}\in\Cc$ to a pair of adjoint functors in $\A$. Thus, to define each adjoint equivalences below, we give a family of pairs of adjoint functors in $\A$, indexed by an object in the appropriate $\Cc$ ($\A^3$ for the associator, and $\A$ for the unitors), and then check that these assemble to give adjoint equivalences between the relevant $2$-functors.

\mdef \label{pa:associator} Let $\oneA,\oneB,\oneC \in \Ob(\A)$. There is a canonical isomorphism of categories 
\begin{align*}
\tilde \ass_{\oneA,\oneB,\oneC}\colon(\oneA\diamond\oneB) \diamond \oneC&\to \oneA \diamond (\oneB\diamond\oneC) \\
((X,Y),Z)&\mapsto (X,(Y,Z))\\
(f\otimes g)\otimes h&\mapsto f\otimes (g\otimes h).
\end{align*}
Let $\tilde\ass_{\oneA,\oneB,\oneC}^*$ be the inverse isomorphism. 

We define
$$\ass_{\oneA,\oneB,\oneC}\colon   (\oneA\monob\oneB)\monob \oneC \to  \oneA \monob (\oneB\monob\oneC)$$
by the commutative diagram
$$\xymatrix{
(\oneA\monob\oneB)\monob \oneC  \ar@{=}[d]\ar^{\ass_{\oneA,\oneB,\oneC}}[rr]&&  \oneA \monob (\oneB\monob\oneC)\ar@{=}[d]\\
\Add (\Add (\oneA\diamond \oneB) \diamond \oneC)
\ar^{(\Add (\iota_{ (\oneA\diamond \oneB)}\diamond \mathbbm{1}_{\oneC}))^*}[d]
&& \Add ( \oneA\diamond \Add (\oneB\diamond \oneC))\\
\Add ( (\oneA\diamond\oneB) \diamond\oneC) \ar^{\Add \tilde \ass_{\oneA,\oneB,\oneC}}[rr] 
&& \Add ( \oneA\diamond (\oneB\diamond \oneC))\ar^{\Add (\mathbbm{1}_{\oneA}\diamond\iota_{ (\oneB\diamond \oneC)} )}[u].
\\
}$$

The functor $\ass_{\oneA,\oneB,\oneC}$ has a (two-sided) adjoint which we label $\ass_{\oneA,\oneB,\oneC}^*$,
and which is defined in the obvious way using the (two-sided) adjoint pairs of functors from \eqref{addadj}, and $\Add \tilde\ass_{\oneA,\oneB,\oneC}^*$. 
Notice that, on an object $\bigoplus_{i=1}^n( \bigoplus_{j=1}^{n_i} (X_{ij}, Y_{ij}), Z_i) \in  (\oneA\monob\oneB)\monob \oneC$, the first step in applying $\ass_{\oneA,\oneB,\oneC}$ pulls out the direct sum, the second step internally rebrackets and the third step then does nothing, and hence
$$\ass_{\oneA,\oneB,\oneC}(\bigoplus_{i=1}^n( \bigoplus_{j=1}^{n_i} (X_{ij}, Y_{ij}), Z_i))  = \bigoplus_{i=1}^n\bigoplus_{j=1}^{n_i} (X_{ij}, (Y_{ij}, Z_i)).$$

\mdef\label{pa:unitcounit} We denote the unit and counit of the adjunction $(\ass_{\oneA,\oneB,\oneC},\ass_{\oneA,\oneB,\oneC}^*)$ by $\eta^\ass_{\oneA,\oneB,\oneC}$ and $\varepsilon^\ass_{\oneA,\oneB,\oneC}$, respectively.
The natural transformation $$\varepsilon^\ass_{\oneA,\oneB,\oneC}\colon \Add(\mathbbm{1}_{\oneA}\diamond \iota_{\oneB\diamond \oneC})\Add(\mathbbm{1}_{\oneA}\diamond \iota_{\oneB\diamond \oneC})^*= a_{\oneA,\oneB,\oneC}a^*_{\oneA,\oneB,\oneC}\to  \mathbbm{1}_{\Add(\oneA\diamond \Add(\oneB\diamond \oneC))} $$ has components which are precisely the rebracketing given explicitly in \eqref{addadjmor}, letting $\oneC$ there be $\oneA$, and $\oneD$ be $\oneB\diamond \oneC$.
The map
$\eta^\ass_{\oneA,\oneB,\oneC}$ differs only in that the sum is pulled out from the left hand component.
Observing that all morphisms in $\Add(\oneA\diamond \Add(\oneB\diamond \oneC))=\oneA\monob (\oneB\monob \oneC)$, resp.\ $\Add(\Add(\oneA \diamond \oneB) \diamond \oneC) = (\oneA\monob \oneB)\monob \oneC$, commute with these rebracketings, we see that the components of $\varepsilon^\ass_{\oneA,\oneB,\oneC}$, resp.\ $\eta^\ass_{\oneA,\oneB,\oneC}$ indeed assemble to a natural transformation.

\begin{prop}\label{assone}
There is an adjoint equivalence $(\ass,\ass^*,\eta^{\ass},\varepsilon^{\ass})$ 
$$ \otimes(\otimes \times \mathbbm{1}_{\A})\underset{\ass^*}{ \overset{\ass}{\rightleftarrows}} \otimes(\mathbbm{1}_{\A} \times \otimes)$$ between $2$-functors $\A\times\A\times\A\to \A$.\end{prop}

\proof
The component $1$-cells of $\ass$ (resp.\ $\ass^*$) are given by  $\ass_{\oneA,\oneB,\oneC}$ (resp.\  $\ass_{\oneA,\oneB,\oneC}^*$) defined above. Since each functor in their construction strictly commutes with any $\Bbbk$-linear functor induced as in \eqref{pa:addclo} (cf. \eqref{addadj}), these assemble to a strict natural transformation.

The components of $\varepsilon^\ass$ and unit $\eta^\ass$ are given by $\varepsilon^\ass_{\oneA,\oneB,\oneC}$ and $\eta^\ass_{\oneA,\oneB,\oneC}$, respectively. The fact that these assemble to modifications, and satisfy the triangle identities, follows from the explicit descriptions of the adjunctions between each of the functors in the construction in \eqref{pa:unitcounit}.
In particular that they are given by rebracketings and thus commute with any natural transformation of induced $\Bbbk$-linear functors.
\endproof

\mdef\label{pa:oneunitors} Let $\oneC$ be a $\Bbbk$-linear additive category. We have an obvious $\Bbbk$-linear functor
$$\begin{array}{rrclll}
\tilde{\lun}_{\oneC}\colon & \Bbbk \diamond \oneC &\to&\oneC&& \\
&(\bullet_\Bbbk, X)&\mapsto &X,&& X\in \oneC\\
&\beta \otimes f & \mapsto &\beta f,&&f\in \oneC(X,X), \beta\in \Bbbk
\end{array}$$
with (bi)adjoint $\tilde{\lun}_{\oneC}^*$ given by the inverse function on objects and the assignment $f\mapsto 1\otimes f$ on morphisms.
This defines a $\Bbbk$-linear functor $\lun_{\oneC}\colon (\Add\,\Bbbk)\monob \oneC\to\oneC $
 by the diagram
 $$\xymatrix{ 
(\Add\,\Bbbk) \monob \oneC \ar@{=}[rr]\ar^{\lun_{\oneC}}[rrrrd]&&  \Add((\Add\,\Bbbk) \diamond \oneC)  \ar^{(\Add (\iota_{\Bbbk}\diamond\mathbbm{1}_{\oneC} ))^*}[rr]&& \Add(\Bbbk \diamond \oneC)\ar_{\tilde\lun_{\oneC}}[d]\\
&&&& \oneC,
}$$
where the vertical arrow is induced as in \eqref {pa:addclo}.
The adjoint $\lun_{\oneC}^*$ is induced by the adjoint diagram.

Analogously, the $\Bbbk$-linear functor 
$$\begin{array}{rrclll}
\tilde{\run}_{\oneC}\colon &  \oneC\diamond\Bbbk &\to&\oneC&& \\
&(X,\bullet_\Bbbk)&\mapsto &X,&& X\in \oneC\\
&f\otimes \beta & \mapsto &\beta f,&&f\in \oneC(X,X), \beta\in \Bbbk
\end{array}$$
and its (bi)adjoint $\tilde{\run}_{\oneC}^*$ given by the inverse function on objects and the assignment $f\mapsto f\otimes 1$ on morphisms, 
define $\Bbbk$-linear functors $\run_{\oneC}\colon \oneC\monob\Add\,\Bbbk \to\oneC $ and its adjoint $\run_{\oneC}^*$ by the diagram
$$\xymatrix{ 
\oneC\monob (\Add\,\Bbbk)  \ar@{=}[rr]\ar^{\run_{\oneC}}[rrrrd]&&  \Add(\oneC \diamond (\Add\,\Bbbk))  \ar^{(\Add(\mathbbm{1}_{\oneC} \diamond \iota_{\Bbbk}))^*}[rr]&& \Add(\oneC \diamond \Bbbk)\ar_{\tilde\run_{\oneC}}[d]\\
&&&& \oneC
}$$
and its adjoint.


\begin{prop}\label{unitorone}
There are adjoint equivalences $(\lun, \lun^*, \eta^\lun, \varepsilon^\lun)$
$$\otimes( \One_{\A} \times \mathbbm{1}_{\A})\underset{\lun^*}{ \overset{\lun}{\rightleftarrows}}\mathbbm{1}_{\A} $$ 
and $(\run, \run^*, \eta^\run, \varepsilon^\run)$
$$ \otimes(\mathbbm{1}_{\A}\times \One_{\A})\underset{\run^*}{ \overset{\run}{\rightleftarrows}} \mathbbm{1}_{\A} $$
between $2$-functors $\A\to \A$.
\end{prop}

\proof
The component functors of $\lun, \lun^*,\run,\run^*$ are given by $\lun_{\oneC}, \lun^*_{\oneC},\run_{\oneC},\run^*_{\oneC}$, respectively. These assemble to strict natural transformations since again both functors in their definition commute with any $\Bbbk$-linear functor $\oneC \to \oneC'$.

There are local adunction morphisms 
$$\begin{array}{rrcl}
\eta^\lun_{\oneC}\colon &\mathbbm{1}_{\Add\Bbbk\monob\oneC}&\to & \lun_{\oneC}^*\lun_{\oneC}\\
\varepsilon^\lun_{\oneC}\colon &\lun_{\oneC}\lun^*_{\oneC} &\to& \mathbbm{1}_{\oneC}\\
\eta^\run_{\oneC}\colon&\mathbbm{1}_{\oneC\monob\Add\Bbbk}&\to & \run_{\oneC}^*\run_{\oneC}\\
\varepsilon^\run_{\oneC}\colon &\run_{\oneC}\run^*_{\oneC} &\to &\mathbbm{1}_{\oneC}
\end{array}$$
where
$\varepsilon^\lun_{\oneC}$ and $\varepsilon^\run_{\oneC}$ are the identity natural transformations, while $\eta^\lun_{\oneC}$ and $\eta^\run_{\oneC}$ are given by rebracketing (c.f.\ Section \ref{addadjmor}).
Since rebracketing commutes with all morphisms in $\Add(\oneC\diamond \Add\Bbbk)$ and $\Add(\Add\Bbbk\diamond \oneC)$, these define natural transformations.
That they satisfy the triangle identities follows from the explicit descriptions of the adjunctions between each of the functors in the construction in \eqref{pa:unitcounit}, thus we have a local adjunctions $(\lun_\oneC, \lun_\oneC^*, \eta^\lun_\oneC, \varepsilon^\lun_\oneC)$ and $(\run_\oneC, \run_\oneC^*, \eta^\run_\oneC, \varepsilon^\run_\oneC)$.

To see that the local units and counits assemble to modifications, it suffices to observe that, since they are given by rebracketings, they commute with all natural transformation of functors $\oneC\to \oneC'$. That they satisfy the triangle identities as modifications then follows from the fact that the local units and counits satisfy the triangle identities as natural transformations. 
\endproof

\begin{prop}\label{pentagone}
For $\Bbbk$-linear, additive categories $\oneA, \oneB, \oneC, \oneD$, there are equalities of functors
$$
(\mathbbm{1}_{\oneA}\monmor\ass_{\oneB, \oneC, \oneD}) \ass_{\oneA, \oneB\monob \oneC, \oneD}( \ass_{\oneA, \oneB, \oneC}\monmor\mathbbm{1}_{\oneD}) = \ass_{\oneA, \oneB, \oneC\monob \oneD}  \ass_{\oneA\monob\oneB, \oneC, \oneD} 
$$
and thus the pentagonator $\pi$ can be chosen as the identity modification.
\end{prop}

\proof
Applied to $\bigoplus_{i=1}^n ( \bigoplus_{j=1}^{m_i} (\bigoplus_{k=1}^{r_{ij}} (W_{ijk}, X_{ijk}),Y_{ij}), Z_i  )$ in 
$((\oneA\monob \oneB)\monob \oneC)\monob\oneD$, both functors yield
$\bigoplus_{i=1}^n  \bigoplus_{j=1}^{m_i} \bigoplus_{k=1}^{r_{ij}} (W_{ijk},( X_{ijk},(Y_{ij}, Z_i  )))$ in  $\oneA\monob( \oneB\monob (\oneC\monob\oneD))$, and they likewise agree on morphisms. Thus each natural transformation $\pi_{\oneA, \oneB, \oneC, \oneD}$ can be chosen as the identity, and these assemble to the pentagonator being the identity modification as claimed.
\endproof

\begin{prop}\label{2unitorone}
For any $\Bbbk$-linear additive categories $\oneC,\oneD$, we have isomorphisms of functors
$$(\mathbbm{1}_{\oneC}\monmor \lun_{\oneD})\, \ass_{\oneC,\Add \Bbbk, \oneD} \,(\run^*_{\oneC}\monmor\mathbbm{1}_{\oneD} )\xrightarrow{\mu_{\oneC,\oneD}}\mathbbm{1}_{\oneC\monob\oneD}$$
$$ \lun_{\oneC}\monmor \mathbbm{1}_{\oneD}
 \xrightarrow{\lambda_{\oneC,\oneD}}\lun_{\oneC\monob\oneD}\, \ass_{\Add \Bbbk,\oneC,\oneD } $$
$$\mathbbm{1}_{\oneC}\monmor \run^*_{\oneD}
 \xrightarrow{\rho_{\oneC,\oneD}} \ass_{\oneC,\oneD,\Add \Bbbk } \,\run^*_{\oneC\monob\oneD} $$
and which assemble to modifications $\mu, \lambda, \rho$, the middle, left and right $2$-unitors.
\end{prop}

\proof
Explicitly, the natural isomorphisms are given by 
$$(\mu_{\oneC,\oneD})_{\bigoplus_{i=1}^n(X_i,Y_i)} = \mathbbm{1}_{\bigoplus_{i=1}^n(X_i,Y_i)} $$
and 
$$(\rho_{\oneC,\oneD})_{\bigoplus_{i=1}^n(X_i,Y_i)}=\mathbbm{1}_{ {\bigoplus_{i=1}^n(X_i,(Y_i, \bullet_\Bbbk))}}$$
for an object $\bigoplus_{i=1}^n(X_i,Y_i)\in \oneC\otimes_0\oneD$,
 and the canonical isomorphism
$$ (\lambda_{\oneC,\oneD})_{Z} \colon \bigoplus_{j=1}^m (\bigoplus_{i=1}^{n_j}  X_{i,j}^{\oplus k_{ij}}, Y_j) \to  \bigoplus_{j=1}^m \bigoplus_{i=1}^{n_j} (X_{i,j}, Y_j)^{\oplus k_{ij}}$$
for an object $Z = \bigoplus_{j=1}^m(\bigoplus_{i=1}^{n_j} (\bullet_{\Bbbk}^{\oplus k_{ij}}, X_{i,j}), Y_j)\in (\Add\Bbbk \monob \oneC)\monob \oneD$. Using strict commutativity of the involved functors with any $\Bbbk$-linear functors $\oneC\to\oneC', \oneD\to\oneD'$, it is immediate that these isomorphisms assemble to modifications, as claimed. 
\endproof

\begin{theorem}\label{th:mon2catof1cats}
$\A$ is a monoidal $2$-category with monoidal product $\otimes$ (Proposition \ref{tensorone}), unit $\One_{\A}$  (Proposition \ref{unitone}), associator $\ass$ (Proposition \ref{assone}), left and right unitors
$\lun$ and $\run$ (Proposition \ref{unitorone}), pentagonator $\pi$  (Proposition \ref{pentagone}) and $2$-unitors $\mu, \lambda$ and $\rho$ (Proposition \ref{2unitorone}).
\end{theorem}

\proof
One checks that the diagrams (11.2.14), (11.2.16) and (11.2.17) in \cite[Def.11.2.4]{JY} commute trivially, since, on each object from the domain, each of the $2$-cells in the pasting diagram is given by an identity.
\endproof

\section{Locally additive $\Bbbk$-linear $2$-categories}\label{sec:add2cats}

Denote by $\twoAk$ the collection whose 
\begin{itemize}
		\item objects are locally $\Bbbk$-linear additive $2$-categories (that is, $2$-categories whose morphism categories are in $\A$);
		\item $1$-morphisms are strictly unital $\Bbbk$-linear pseudofunctors; 
		\item $2$-morphisms are strong natural transformations;
		\item $3$-morphisms are modifications.
	\end{itemize}

Note that since the objects of $\A$ are small $\Bbbk$-linear additive categories, we have well-defined functor categories between two objects in $\A$, which are again $\Bbbk$-linear and additive, so $\A$ is indeed an object of $\twoAk$.
	
We remark that this does not form a tricategory, since $1$-, $2$- and $3$-morphisms from an object $\Cc$ to an object $\Dd$ only form a bicategory if $\Cc$ is small (i.e. has a set of objects). Since we do want to consider $\A$ as an object of $\twoAk$, we thus refrain from giving a name to this structure. If we restrict objects to {\em small} locally $\Bbbk$-linear additive $2$-categories, we obtain a tricategory.

We would like to just consider objects, $1$-morphisms and $2$-morphisms, and endow this with a structure of monoidal $2$-category. This is, of course, impossible, since it is not a $2$-category. However, we can define tensor products of $2$-categories in $\twoAk$, strictly unital pseudofunctors, and strong transformations, and find unitors and associators, as well as $2$-unitors and pentagonators, which is enough for our application.


\subsection{Tensor products of $2$-categories, pseudofunctors and strong transformations}
\mdef\label{CboxD}
Given two $2$-categories $\Cc,\Dd$ in $\twoAk$, we construct a new $2$-category $\Cc\bxob \Dd$ which has 
\begin{itemize}
\item objects $(\ob{i},\ob{i}')\in \Ob(\Cc)\times \Ob(\Dd)$;
\item morphism categories $(\Cc\bxob \Dd)((\ob{i},\ob{i}'),(\ob{j},\ob{j}'))= \Cc(\ob{i},\ob{j}) \monob\Dd(\ob{i}',\ob{j}')$ for objects $(\ob{i},\ob{i}'),(\ob{j},\ob{j}')\in  \Ob(\Cc)\times \Ob(\Dd)$.
\end{itemize}

To define horizontal composition, consider the functor
uniquely defined by 
\begin{align*}
 (\Cc(\ob{j},\ob{k})\diamond \Dd(\ob{j}',\ob{k}')) \times  (\Cc(\ob{i},\ob{j})\diamond \Dd(\ob{i}',\ob{j}')) &\to \Cc(\ob{i},\ob{k}) \diamond\Dd(\ob{i}',\ob{k}')\\
( (\G,\G'), (\F, \F')) &\mapsto (\G\F, \G'\F')\\
((\alpha \otimes_\Bbbk \alpha'),(\beta \otimes_\Bbbk \beta'))&\mapsto(\alpha \hcirc \beta)\otimes_\Bbbk (\alpha' \hcirc \beta')
\end{align*}
with the appropriate bilinear extension to general morphisms,
and notice that, using the canonical functor from \eqref{naiveotimes}, this factors over
$(\Cc(\ob{j},\ob{k})\diamond \Dd(\ob{j}',\ob{k}')) \diamond (\Cc(\ob{i},\ob{j})\diamond \Dd(\ob{i}',\ob{j}'))$ giving rise to a functor
\begin{align*}\tilde{\circ}\colon (\Cc(\ob{j},\ob{k})\diamond\Dd(\ob{j}',\ob{k}'))\diamond (\Cc(\ob{i},\ob{j})\diamond \Dd(\ob{i}',\ob{j}')) &  \to \Cc(\ob{i},\ob{k}) \diamond\Dd(\ob{i}',\ob{k}')\\
( (\G,\G'), (\F, \F')) &\mapsto (\G\F, \G'\F')\\
((\alpha \otimes_\Bbbk \alpha')\otimes_\Bbbk(\beta \otimes_\Bbbk \beta'))&\mapsto(\alpha \hcirc \beta)\otimes_\Bbbk (\alpha' \hcirc \beta').
\end{align*}


Horizontal composition $\circ$ is then given by the composition of functors

$$\xymatrix{(\Cc(\ob{j},\ob{k})\monob \Dd(\ob{j}',\ob{k}'))\times (\Cc(\ob{i},\ob{j})\monob \Dd(\ob{i}',\ob{j}'))\ar[d]\\
(\Cc(\ob{j},\ob{k})\monob \Dd(\ob{j}',\ob{k}')) \diamond(\Cc(\ob{i},\ob{j})\monob \Dd(\ob{i}',\ob{j}'))
\ar^{\iota_{(\Cc(\ob{j},\ob{k})\monob \Dd(\ob{j}',\ob{k}')) \diamond (\Cc(\ob{i},\ob{j})\monob \Dd(\ob{i}',\ob{j}'))}}[d]\\
\Add \left((\Cc(\ob{j},\ob{k})\monob \Dd(\ob{j}',\ob{k}'))\diamond (\Cc(\ob{i},\ob{j})\monob \Dd(\ob{i}',\ob{j}'))\right)
\ar^{\left( \Add\left( \iota_{(\Cc(\ob{j},\ob{k})\diamond \Dd(\ob{j}',\ob{k}'))}\diamond \iota_{(\Cc(\ob{i},\ob{j})\diamond \Dd(\ob{i}',\ob{j}'))}  \right)\right)^* }[d]
\\
 \Add\left((\Cc(\ob{j},\ob{k})\diamond \Dd(\ob{j}',\ob{k}')) \diamond  (\Cc(\ob{i},\ob{j})\diamond \Dd(\ob{i}',\ob{j}'))\right)
 \ar^{\Add(\tilde\circ)}[d]\\
  \Cc(\ob{i},\ob{k}) \monob\Dd(\ob{i}',\ob{k}').
} $$

Explicitly, this means that
$$\left( \bigoplus_j(\F_j,\F_j') \right)\left(\bigoplus_i(\G_i,\G_i')\right) = \bigoplus_j\bigoplus_i (\F_j\G_i,\F_j'\G_i') $$
for $\F_j \in \Cc(\ob{j},\ob{k}),\F_j' \in \Cc(\ob{j}',\ob{k}'), \G_i\in \Cc(\ob{i},\ob{j}), \G_i'\in \Cc(\ob{i}',\ob{j}')$.

%
%

\mdef By construction, $ \Cc\bxob \Dd$ is locally $\Bbbk$-linear and additive. It follows directly from the definition of horizontal composition that $\mathbbm{1}_{(\ob{i}, \ob{j})}\in (\Cc\bxob \Dd)((\ob{i}, \ob{j}),(\ob{i}, \ob{j}))$ is given by
$(\mathbbm{1}_{\ob{i}},\mathbbm{1}_{ \ob{j}}).$
Associativity of horizontal composition follows from its above description as induced by component-wise composition in $\Cc$ and $\Dd$.
Likewise, the interchange law reduces to the componentwise interchange laws. 
Thus, the above data defines a function
\begin{align*}
\bxob \colon \Ob(\twoAk) \times \Ob(\twoAk) &\to \Ob(\twoAk)\\
(\Cc,\Dd) &\to \Cc\bxob \Dd.
\end{align*}


\mdef\label{FboxG}
Now assume $\Cc,\Dd, \Cc',\Dd'\in\twoAk$, and let $\Ff\colon \Cc\to \Cc'$ and $\Gg\colon \Dd\to \Dd'$ be strictly unital pseudofunctors with compositionality isomorphisms $\Ff^2, \Gg^2$, respectively.
We wish to construct a strictly unital pseudofunctor $\Ff\bxone\Gg\colon \Cc\bxob \Dd\to \Cc'\bxob \Dd'$. On objects, this is simply defined as 
$$(\Ff\bxone\Gg)(\ob{i},\ob{k}) =  (\Ff(\ob{i}),\Gg(\ob{k})). $$
Denoting, for $\ob{i},\ob{j}\in \Cc, \ob{k},\ob{l}\in \Dd$, by
$\Ff_{\ob{j}\ob{i}}\colon  \Cc(\ob{i},\ob{j})\to \Cc'(\Ff(\ob{i}),\Ff(\ob{j}))$ and $\Gg_{\ob{l}\ob{k}}\colon  \Dd(\ob{k},\ob{l})\to \Dd'(\Gg(\ob{k}),\Gg(\ob{l}))$ the corresponding component functors of $\Ff$ and $\Gg$, respectively, we, moreover, set
$$(\Ff\bxone\Gg)_{(\ob{j},\ob{l})(\ob{i},\ob{k})}= \Ff_{\ob{j}\ob{i}}\monmor\Gg_{\ob{l}\ob{k}}.$$
Note this preserves identities.
The compositionality isomorphism $(\Ff\bxone\Gg)^2$ is given by  
$$(\Ff\bxone\Gg)^2_{(\ob{j},\ob{m}),(\ob{i},\ob{l}),(\ob{h},\ob{k})}=\Ff^2_{\ob{j}\ob{i}\ob{h}}\monnat\Gg^2_{\ob{m}\ob{l}\ob{k}}.$$

\mdef 
By construction, the component functors are $\Bbbk$-linear functors. Strict unitality follows immediately from strict unitality of $\Ff$ and $\Gg$. Moreover, $(\Ff\bxone\Gg)^2$ is a natural transformation since it is the natural transformation constructed in \eqref{otimesdef}, an isomorphism by construction, and it satisfies lax associativity since $\Ff^2$ and $\Gg^2$ do. Thus, $\Ff\bxone \Gg$ defines a strictly unital pseudofunctor from $\Cc\bxob\Dd$ to $\Cc'\bxob \Dd'$.



\mdef
We observe that for $\Bbbk$-linear strictly unital pseudofunctors
$$\Ff \colon\Bb \to \Cc, \Ff'\colon\Bb' \to \Cc', \Gg\colon \Aa\to \Bb, \Gg'\colon \Aa'\to \Bb',$$ we have
$$\Ff\Gg\bxone \Ff'\Gg' = (\Ff\bxone\Ff')  (\Gg\bxone\Gg')$$
and that  $\mathbbm{1}_{\Aa\bxob\Bb} = \mathbbm{1}_{\Aa}\bxone\mathbbm{1}_{\Bb}$.


\mdef
Now let $\Cc,\Dd,\Cc',\Dd'\in\twoAk$, let $\Ff,\Ff'\colon \Cc\to \Cc'$ and $\Gg,\Gg'\colon \Dd\to \Dd'$ be strictly unital pseudofunctors
and $\phi\colon \Ff\to \Ff'$ and $\psi\colon \Gg\to \Gg'$ be strong natural transformations given, for each $\ob{i}\in\Cc,\ob{j}\in\Dd$, by $1$-morphisms $\phi_{\ob{i}}\in \Cc'(\Ff(\ob{i}),\Ff'(\ob{i})), \psi_{\ob{j}}\in \Dd'(\Gg(\ob{j}),\Gg'(\ob{j}))$ and, 
for each $\ob{i},\ob{k}\in\Cc,\ob{j},\ob{l}\in\Dd$, a natural transformation
$$\phi_{\ob{k},\ob{i}} \colon  \mathnormal{R}_{ \phi_{\ob{i}}}\Ff'_{\ob{k}\ob{i}} \to \mathnormal{L}_{ \phi_{\ob{k}}}  \Ff_{\ob{k}\ob{i}} \qquad \psi_{\ob{l},\ob{j}} \colon  \mathnormal{R}_{ \phi_{\ob{j}}}\Gg'_{\ob{l}\ob{j}}  \to  \mathnormal{L}_{ \psi_{\ob{l}}} \Gg_{\ob{l}\ob{j}}$$
of functors  $\Cc(\ob{i},\ob{k})\to\Cc'(\Ff(\ob{i}),\Ff'(\ob{k}))$ and $\Dd(\ob{j},\ob{l})\to \Dd'(\Gg(\ob{j}),\Gg'(\ob{l}))$.

\mdef \label{phiboxpsi}
Define $\phi\bxtwo\psi$ as having, for each object $(\ob{i},\ob{j})\in \Ob(\Cc\bxob \Dd)$ the $1$-morphism $$(\phi_{\ob{i}}, \psi_{\ob{j}}) \in \Cc'(\Ff(\ob{i}),\Ff'(\ob{i}))\monob\Dd'(\Gg(\ob{j}),\Gg'(\ob{j})) = \left(\Cc'\bxob\Dd' \right)((\Ff(\ob{i}),\Gg(\ob{j})),(\Ff'(\ob{i}),\Gg'(\ob{j}))).$$
To define the natural transformations
$$(\phi\bxtwo\psi)_{(\ob{k},\ob{l}),(\ob{i},\ob{j})}\colon  
 \mathnormal{R}_{ (\phi_{\ob{i}}, \psi_{\ob{j}})} \left( \Ff'_{\ob{k}\ob{i}}\monmor \Gg'_{\ob{l}\ob{j}} \right) \to
 \mathnormal{L}_{ (\phi_{\ob{k}}, \psi_{\ob{l}})  }\left( \Ff_{\ob{k}\ob{i}}\monmor \Gg_{\ob{l}\ob{j}} \right)
 $$
notice that 
$$ \mathnormal{L}_{ (\phi_{\ob{k}}, \psi_{\ob{l}}) } =  \mathnormal{L}_{ \phi_{\ob{k}}}\monmor  \mathnormal{L}_{ \psi_{\ob{l}} }  \qquad \text{and}\qquad \mathnormal{R}_{(\phi_{\ob{i}}, \psi_{\ob{j}}) } =  \mathnormal{R}_{ \phi_{\ob{i}}}\monmor  \mathnormal{R}_{ \psi_{\ob{j}} } $$
so 
\begin{equation*}\begin{split}
\mathnormal{L}_{  (\phi_{\ob{k}}, \psi_{\ob{l}})  }\left( \Ff_{\ob{k}\ob{i}}\monmor  \Gg_{\ob{l}\ob{j}} \right) &= \left( \mathnormal{L}_{ \phi_{\ob{k}}}\monmor   \mathnormal{L}_{ \psi_{\ob{l}} }\right)\left( \Ff_{\ob{k}\ob{i}}\monmor  \Gg_{\ob{l}\ob{j}} \right)\\
&= \left( \mathnormal{L}_{ \phi_{\ob{k}}}  \Ff_{\ob{k}\ob{i}}\right)\monmor \left(  \mathnormal{L}_{ \psi_{\ob{l}} } \Gg_{\ob{l}\ob{j}} \right)\end{split}\end{equation*}
and likewise
$$\mathnormal{R}_{ (\phi_{\ob{i}}, \psi_{\ob{j}})} \left( \Ff'_{\ob{k}\ob{i}}\bxone \Gg'_{\ob{l}\ob{j}} \right) =  \left( \mathnormal{R}_{ \phi_{\ob{i}}}  \Ff'_{\ob{k}\ob{i}}\right)\monmor \left(  \mathnormal{R}_{ \psi_{\ob{j}} } \Gg'_{\ob{l}\ob{j}} \right)$$
so we can set 
$$(\phi\bxtwo\psi)_{(\ob{k},\ob{l}),(\ob{i},\ob{j})} = \phi_{\ob{k},\ob{i}}\monnat \psi_{\ob{l},\ob{j}}.$$

\mdef Lax naturality of $ \phi\bxtwo\psi$ follows from lax naturality of $\phi$ and $\psi$, hence $ \phi\bxtwo\psi$ indeed defines a strong natural transformation from $\Ff\bxone \Gg$ to $\Ff'\bxone\Gg'$.
%

%
%
%
%
%
%
%

\subsection{Unitors for the tensor product of $2$-categories}\label{boxunit}

Denote by $ \Vv$ the $2$-category, with one object $\bullet$ whose endomorphism category is given by $\Add \Bbbk$ with horizontal composition given by tensor product, explicitly, $(\bullet_{\Bbbk}^{\oplus n})\circ(\bullet_{\Bbbk}^{\oplus m}) = (\bullet_{\Bbbk}^{\oplus nm})$. Note that $\mathbbm{1}_\bullet = \bullet_{\Bbbk}$ in $\Vv$.

We now turn to the left and right unitors. Recall $\lun_{\Cc}, \lun^*_{\Cc}, \run_{\Cc}$ and $\run^*_{\Cc}$ from \eqref{pa:oneunitors}.

There is a strict $2$-adjunction $(\luntwo_{\Cc}, \luntwo^*_{\Cc}, \eta^{\luntwo_{\Cc}}, \varepsilon^{\luntwo_{\Cc}})$ between $\Vv\bxob\Cc$ and $\Cc$, in the sense of \cite[Sec. I.7.1]{Gray}. In more detail, we have the following.

\begin{prop}\label{unitortwopartsleft}
\begin{enumerate}[(a)]
\item\label{unitortwopartsleft.1}For any $\Cc\in \Ob(\twoAk)$, there is a $\Bbbk$-linear $2$-functor $\luntwo_{\Cc}\colon \Vv\bxob\Cc \to\Cc$ given by $(\bullet, \ob{i})\mapsto \ob{i} $ on objects and by $(\luntwo_{\Cc})_{ (\ob{i}, \ob{j})} = \lun_{\Cc(\ob{i}, \ob{j})}$ on morphism categories.
\item\label{unitortwopartsleft.2} For any $\Cc\in \Ob(\twoAk)$, there is a $\Bbbk$-linear $2$-functor $\luntwo^*_{\Cc}\colon \Cc \to\Vv\bxob\Cc$ given by $\ob{i}\mapsto  (\bullet, \ob{i})$ on objects and by $(\luntwo^*_{\Cc})_{ (\ob{i}, \ob{j})} = \lun^*_{\Cc(\ob{i}, \ob{j})}$ on morphism categories.
\item\label{unitortwopartsleft.3} 
There is a strong natural transformation  $\eta^{\luntwo_{\Cc}}\colon\mathbbm{1}_{\Vv\bxob \Cc}\to  \luntwo_{\Cc}^*\luntwo_{\Cc}$
given, for objects $(\bullet, \ob{i})\in\Vv\bxob\Cc$, by
$$ \eta^{\luntwo_{\Cc}}_{(\bullet, \ob{i})} = \mathbbm{1}_{(\bullet, \ob{i})}$$
together with, for $1$-morphisms $\F_r\in \Cc(\ob{i},\ob{j})$, the isomorphism 
\[
(\eta^\lun_{\Cc(\ob{i}, \ob{j})})_{\bigoplus_r(\mathbbm{1}_\bullet^{\oplus n_r}, \F_r)}\colon \bigoplus_r(\mathbbm{1}_\bullet^{\oplus n_r}, \F_r)\to \luntwo_{\Cc}^*\luntwo_{\Cc}\left(\bigoplus_r (\mathbbm{1}_\bullet^{\oplus n_r}, \F_r)\right).
\]
\item\label{unitortwopartsleft.4} There is a strong (even strict) natural transformation $\varepsilon^{\luntwo_{\Cc}}\colon \luntwo_{\Cc}\luntwo_{\Cc}^*\to \mathbbm{1}_{\Cc}$ given, for objects $\ob{i}\in\Cc$, by
$$\varepsilon^{\luntwo_{\Cc}}_{ \ob{i}} =  \mathbbm{1}_{ \ob{i}}$$ 
together with, for $\F\in \Cc(\ob{i},\ob{j})$, the isomorphism between
 $\luntwo_{\Cc}\luntwo_{\Cc}^*(\F)$ and $\F$ 
 given by $(\varepsilon^\lun_{\Cc(\ob{i}, \ob{j})})_{\F} = \mathbbm{1}_{\F}$.
\item\label{unitortwopartsleft.5} 
The strong natural transformations $\eta^{\luntwo_{\Cc}}$ and $\varepsilon^{\luntwo_{\Cc}}$ satisfy the triangle identities strictly, i.e. the diagrams
\[\xymatrix{
\luntwo_\Cc\ar^{\luntwo_\Cc\star \eta^{\luntwo_\Cc}}[rr]\ar@{=}[drr]&&\luntwo_\Cc\luntwo_\Cc^*\luntwo_\Cc\ar^{\varepsilon^{\luntwo_\Cc}\star\luntwo_\Cc}[d]\\
&&\luntwo_\Cc
}
\xymatrix{
\luntwo_\Cc^*\ar^{ \eta^{\luntwo_\Cc}\star \luntwo_\Cc^*}[rr]\ar@{=}[drr]&&\luntwo_\Cc^*\luntwo_\Cc\luntwo_\Cc^*\ar^{\luntwo_\Cc^*\star\varepsilon^{\luntwo_\Cc}}[d]\\
&&\luntwo_\Cc^*
}
\]
commute strictly in the sense that their appropriate compositions form the identity strong natural transformations on $\luntwo_{\Cc}$ and $\luntwo^*_{\Cc}$, respectively.
\end{enumerate}
\end{prop}

\proof
\eqref{unitortwopartsleft.1} We need to check compositionality and unitality of $\luntwo_{\Cc}$, the latter being obvious. For the former, given $1$-morphisms $(\mathbbm{1}_\bullet^{\oplus n}, \F)$ and $(\mathbbm{1}_\bullet^{\oplus m}, \G)$ in $(\Vv\bxob\Cc)( (\bullet, \ob{i}), (\bullet, \ob{j}))$ and $(\Vv\bxob\Cc)(  (\bullet, \ob{j}),  (\bullet, \ob{k}))$, respectively, one explicitly computes 
\begin{equation*}\begin{split}
\lun_{\Cc(\ob{j}, \ob{k})}(\mathbbm{1}_\bullet^{\oplus m}, \G)\lun_{\Cc(\ob{i}, \ob{j})}(\mathbbm{1}_\bullet^{\oplus n}, \F)
& = \G^{\oplus m}\F^{\oplus n} \\
&= (\G\F)^{\oplus mn}\\
&=\lun_{\Cc(\ob{i}, \ob{k})} (\mathbbm{1}_\bullet^{\oplus mn}, \G\F)\\
&=\lun_{\Cc(\ob{i}, \ob{k})} ((\mathbbm{1}_\bullet^{\oplus m}, \G)(\mathbbm{1}_\bullet^{\oplus n}, \F)).
\end{split}\end{equation*}

\eqref{unitortwopartsleft.2} This is similar to part \eqref{unitortwopartsleft.1}, but easier, since no direct sums appear.

\eqref{unitortwopartsleft.3} Using the explicit description of 
$$(\eta^\lun_{\Cc(\ob{i}, \ob{j})})_{\bigoplus_r(\mathbbm{1}_\bullet^{\oplus n_r}, \F_r)}\colon \bigoplus_r(\mathbbm{1}_\bullet^{\oplus n_r}, \F_r)\to \luntwo_{\Cc}^*\luntwo_{\Cc}\left(\bigoplus_r (\mathbbm{1}_\bullet^{\oplus n_r}, \F_r)\right) =(\mathbbm{1}_\bullet, \bigoplus_r\F_r^{\oplus n_r})$$ we check that 
for $\F_r\in \Cc(\ob{i},\ob{j})$, $\G_s\in \Cc(\ob{j},\ob{k})$
the diagram
$$ \xymatrix{
\bigoplus_s(\mathbbm{1}_\bullet^{\oplus m_s}, \G_s)\bigoplus_r(\mathbbm{1}_\bullet^{\oplus n_r}, \F_s) \ar@{=}[rr]
\ar^{(\eta^\lun_{\Cc(\ob{j}, \ob{k})})_{\bigoplus_s(\mathbbm{1}_\bullet^{\oplus m_s}, \G_s)}\hcirc(\eta^\lun_{\Cc(\ob{i}, \ob{j})})_{\bigoplus_r(\mathbbm{1}_\bullet^{\oplus n_r}, \F_r)}}[d]
&&\bigoplus_s\bigoplus_r(\mathbbm{1}_\bullet^{\oplus m_sn_r}
, \G_s\F_r)
\ar^{(\eta^\lun_{\Cc(\ob{i}, \ob{k})})_{\bigoplus_s\bigoplus_r(\mathbbm{1}_\bullet^{\oplus m_sn_r}, \G_s\F_r)}}[d]\\
(\mathbbm{1}_\bullet, \bigoplus_s\G_s^{\oplus m_s})(\mathbbm{1}_\bullet, \bigoplus_r\F_r^{\oplus n_r})\ar@{=}[rr]&& (\mathbbm{1}_\bullet, \bigoplus_s\bigoplus_r\G_s\F_r^{\oplus m_sn_r})
}$$
commutes and hence $$(\eta^\lun_{\Cc(\ob{j}, \ob{k})})_{\bigoplus_s(\mathbbm{1}_\bullet^{\oplus m_s}, \G_s)}\hcirc(\eta^\lun_{\Cc(\ob{i}, \ob{j})})_{\bigoplus_r(\mathbbm{1}_\bullet^{\oplus n_r}, \F_r)} = (\eta^\lun_{\Cc(\ob{i}, \ob{k})})_{\bigoplus_s\bigoplus_r(\mathbbm{1}_\bullet^{\oplus m_sn_r}, \G_s\F_r)},$$ as required for $\eta^{\luntwo_{\Cc}}$ to be a strong natural transformation.

\eqref{unitortwopartsleft.4} This is similar to part \eqref{unitortwopartsleft.3}, but easier, since all components $(\varepsilon^\lun_{\Cc(\ob{i}, \ob{j})})_{\F} $ are identities.

\eqref{unitortwopartsleft.5}This follows from the fact that the same holds for the components $\eta^\lun_{\Cc(\ob{i}, \ob{j})}$, $\varepsilon^\lun_{\Cc(\ob{i}, \ob{j})}$.
\endproof

Similarly, there is a strict $2$-adjunction $(\runtwo_{\Cc}, \runtwo^*_{\Cc}, \eta^{\runtwo_{\Cc}}, \varepsilon^{\runtwo_{\Cc}})$ between $\Cc\bxob\Vv$ and $\Cc$, i.e.\ we have the following.

\begin{prop}\label{unitortwopartsright}
\begin{enumerate}[(a)]
\item\label{unitortwopartsright.1} For any $\Cc\in \Ob(\twoAk)$, there is a $\Bbbk$-linear  $2$-functor $\runtwo_{\Cc}\colon \Cc\bxob\Vv \to\Cc$ given by $(\ob{i}, \bullet)\mapsto \ob{i} $ on objects and by $(\runtwo_{\Cc})_{ (\ob{i}, \ob{j})} = \run_{\Cc(\ob{i}, \ob{j})}$ on morphism categories.
\item\label{unitortwopartsright.2} For any $\Cc\in \Ob(\twoAk)$, there is a $\Bbbk$-linear $2$-functor $\runtwo^*_{\Cc}\colon \Cc \to\Cc\bxob\Vv$ given by $\ob{i}\mapsto  (\ob{i},\bullet)$ on objects and by $(\runtwo^*_{\Cc})_{ (\ob{i}, \ob{j})} = \run^*_{\Cc(\ob{i}, \ob{j})}$ on morphism categories.
\item\label{unitortwopartsright.3} There is a strong natural transformation  $\eta^{\runtwo_{\Cc}}\colon\mathbbm{1}_{\Cc\bxob \Vv}\to  \runtwo_{\Cc}^*\runtwo_{\Cc}$
given, for objects $( \ob{i},\bullet)\in\Vv\bxob\Cc$, by
$$ \eta^{\runtwo_{\Cc}}_{( \ob{i},\bullet)} = \mathbbm{1}_{( \ob{i},\bullet)}$$
together with, for $1$-morphisms $\F_r\in \Cc(\ob{i},\ob{j})$, the isomorphism from $\runtwo_{\Cc}^*\runtwo_{\Cc}\left(\bigoplus_r (\F_r,\mathbbm{1}_\bullet^{\oplus n_r})\right)$ to $\bigoplus_r( \F_r,\mathbbm{1}_\bullet^{\oplus n_r})$ and given by 
$(\eta^\run_{\Cc(\ob{i}, \ob{j})})_{\bigoplus_r(\F_r,\mathbbm{1}_\bullet^{\oplus n_r})}$.
\item\label{unitortwopartsright.4} There is a strong (even strict) natural transformation $\varepsilon^{\runtwo_{\Cc}}\colon \runtwo_{\Cc}\runtwo_{\Cc}^*\to \mathbbm{1}_{\Cc}$ given, for objects $\ob{i}\in\Cc$, by
$$\varepsilon^{\runtwo_{\Cc}}_{ \ob{i}} =  \mathbbm{1}_{ \ob{i}}$$ 
together with, for $\F\in \Cc(\ob{i},\ob{j})$, the isomorphism from $\F$ to
 $\runtwo_{\Cc}\runtwo_{\Cc}^*(\F)$
 given by $(\varepsilon^\run_{\Cc(\ob{i}, \ob{j})})_{\F} = \mathbbm{1}_{\F}$.
\item\label{unitortwopartsright.5} 
The strong natural transformations $\eta^{\runtwo_{\Cc}}$ and $\varepsilon^{\runtwo_{\Cc}}$ satisfy the triangle inequality strictly, i.e. the diagrams
\[\xymatrix{
\runtwo_\Cc\ar^{\runtwo_\Cc\star \eta^{\runtwo_\Cc}}[rr]\ar@{=}[drr]&&\runtwo_\Cc\runtwo_\Cc^*\runtwo_\Cc\ar^{\varepsilon^{\runtwo_\Cc}\star\runtwo_\Cc}[d]\\
&&\runtwo_\Cc
}
\xymatrix{
\runtwo_\Cc^*\ar^{ \eta^{\runtwo_\Cc}\star \runtwo_\Cc^*}[rr]\ar@{=}[drr]&&\runtwo_\Cc^*\runtwo_\Cc\runtwo_\Cc^*\ar^{\runtwo_\Cc^*\star\varepsilon^{\runtwo_\Cc}}[d]\\
&&\runtwo_\Cc^*
}
\]
commute strictly.
\end{enumerate}
\end{prop}

\proof 
Mutatis mutandis the proof of Proposition \ref{unitortwopartsleft}.
\endproof

\begin{prop}\label{unitortwocommute}
For any strictly unital pseudofunctor $\Ff\colon \Cc\to\Dd$ in $\twoAk$, there exist invertible strong transformations
\[\luntwo_\Ff\colon \Ff \luntwo_\Cc\to \luntwo_\Dd(\mathbbm{1}_{\Vv}\bxone \Ff) \]
and
\[\runtwo_\Ff\colon \Ff\runtwo_\Cc \to \runtwo_\Dd(\Ff\bxone \mathbbm{1}_{\Vv})\]
\end{prop}

\proof
To define $\luntwo_\Ff$, note first that both the domain and codomain pseudofunctors act in the same way on objects, sending $(\bullet,\ob{i})$ to $\Ff(\ob{i})$, so the $1$-cell components of $\luntwo_\Ff$ may be chosen to be identities.
For a $1$-morphism $\bigoplus_s(\mathbbm{1}_\bullet^{\oplus a_s},\F_s)$ in $\Vv\bxob\Cc$, we have
\[
\Ff \luntwo_\Cc \bigl(\bigoplus_s(\mathbbm{1}_\bullet^{\oplus a_s},\F_s)\bigr) = 
\Ff(\bigoplus_s\F_s^{\oplus a_s}),
\]
while
\[
\luntwo_{\Dd}(\mathbbm{1}_{\Vv}\bxone \Ff)
\bigl(\bigoplus_s(\mathbbm{1}_\bullet^{\oplus a_s},\F_s)\bigr)
=\bigoplus_s(\Ff(\F_s)^{\oplus a_s})
\]
thus we may choose the $2$-cell components of $\luntwo_\Ff$ to be given by the canonical isomorphisms between these given by the additivity of $\Ff$. To check lax naturality, note that, for given $1$-morphisms $\F=\bigoplus_s(\mathbbm{1}_\bullet^{\oplus a_s},\F_s)$ and $\G=\bigoplus_t(\mathbbm{1}_\bullet^{\oplus b_t},\G_t)$ in $\Vv\bxob\Cc$, the diagram
\[\xymatrix{
\Ff(\bigoplus_s\F_s^{\oplus a_s})\Ff(\bigoplus_t\G_t^{\oplus b_t}) \ar^{\Ff^2_{\luntwo_\Cc(\F),\luntwo_\Cc(\G)}}[rr]\ar_{(\luntwo_\Ff)_{\F}\hcirc(\luntwo_\Ff)_{\G}}[d]&&\Ff(\bigoplus_s\bigoplus_t(\F_s\G_t)^{\oplus a_sb_t})\ar^{(\luntwo_\Ff)_{\F\G}}[d]\\
\bigoplus_s(\Ff(\F_s)^{\oplus a_s})\bigoplus_t(\Ff(\G_t)^{\oplus b_t}) \ar_{\Ff^2}[rr]&& \bigoplus_s\bigoplus_t(\Ff(\F_s\G_t)^{\oplus a_sb_t})\\
}\]
where the bottom horizontal arrow is the matrix with entries $\Ff^2_{\F_s,\G_t}$, commutes by additivity of $\Ff$.

Similarly, to define $\runtwo_\Ff\colon \Ff\runtwo_\Cc \to \run_{\Dd}(\Ff\bxone \mathbbm{1}_{\Vv})$, note that both pseudofunctors again act in the same way on objects, namely by sending $(\ob{i},\bullet)$ to $\Ff(\ob{i})$, so the $1$-cell components of $\runtwo_\Ff$ are the identities.
For a $1$-morphism $\bigoplus_s(\F_s, \mathbbm{1}_\bullet^{\oplus a_s})$ in $\Cc\bxob \Vv$, we have 
$$ (\Ff\runtwo_\Cc)(\bigoplus_s(\F_s, \mathbbm{1}_\bullet^{\oplus a_s})) = \Ff(\bigoplus_s\F_s^{\oplus a_s})$$
and 
$$\runtwo_{\Dd}(\Ff\bxone \mathbbm{1}_{\Vv})(\bigoplus_s(\F_s, \mathbbm{1}_\bullet^{\oplus a_s})) = \runtwo_{\Dd}(\bigoplus_s(\Ff(\F_s), \mathbbm{1}_\bullet^{\oplus a_s})) = \bigoplus_s \Ff(\F_s)^{\oplus a_s},$$
thus the $2$-cell components of $\runtwo_\Ff$ are chosen as the isomorphism between these given by the additivity of $\Ff$. Lax naturality follows analogously to that of $\luntwo_\Ff$.
\endproof

\subsection{Associator for the tensor product of $2$-categories}\label{boxass}

In this subsection we show that there is a strict $2$-adjunction $(\asstwo_{\Cc}, \asstwo^*_{\Cc}, \eta^{\asstwo_{\Cc}}, \varepsilon^{\asstwo_{\Cc}})$ between $(\Aa\bxob\Bb)\bxob\Cc$ and $\Aa\bxob(\Bb\bxob\Cc)$, i.e.\ we have the following.

\begin{prop}\label{ass2fun} Let $\Aa,\Bb,\Cc \in \Ob(\twoAk)$.
\begin{enumerate}[(a)]
\item\label{ass2fun.1} 
The assignment
$$((\ob{i},\ob{j}),\ob{k})\mapsto (\ob{i},(\ob{j},\ob{k}))$$
for objects 
and
$$(\asstwo_{\Aa,\Bb,\Cc})_{((\ob{i}',\ob{j}'),\ob{k}'),((\ob{i},\ob{j}),\ob{k})} = \ass_{\Aa(\ob{i},\ob{i}'),\Bb(\ob{j},\ob{j}'),\Cc(\ob{k},\ob{k}')}$$
for morphism categories (using notation from \eqref{pa:associator}) defines a $\Bbbk$-linear strictly unital pseudofunctor 
$$\asstwo_{\Aa,\Bb,\Cc}\colon (\Aa\bxob\Bb)\bxob\Cc \to \Aa\bxob(\Bb\bxob\Cc).$$
\item\label{ass2fun.2} 
The assignment
$$(\ob{i},(\ob{j},\ob{k}))\mapsto ((\ob{i},\ob{j}),\ob{k})$$
for objects 
and
$$(\asstwo^*_{\Aa,\Bb,\Cc})_{(\ob{i}',(\ob{j}',\ob{k}')),(\ob{i},(\ob{j},\ob{k}))} = \ass^*_{\Aa(\ob{i},\ob{i}'),\Bb(\ob{j},\ob{j}'),\Cc(\ob{k},\ob{k}')}$$
for morphism categories defines a $\Bbbk$-linear strictly unital pseudofunctor 
$$\asstwo^*_{\Aa,\Bb,\Cc}\colon \Aa\bxob(\Bb\bxob\Cc) \to  (\Aa\bxob\Bb)\bxob\Cc.$$
\item \label{ass2fun.3}
There is a strong natural transformation $\eta^{\asstwo_{\Aa,\Bb,\Cc}}\colon \mathbbm{1}_{ \Aa\bxob(\Bb\bxob\Cc)} \to \asstwo_{\Aa,\Bb,\Cc}^*\asstwo_{\Aa,\Bb,\Cc}$ given, for objects $((\ob{i},\ob{j}),\ob{k})\in \Aa\bxob(\Bb\bxob\Cc)$, by
$$\eta^{\asstwo_{\Aa,\Bb,\Cc}}_{((\ob{i},\ob{j}),\ob{k})} = \mathbbm{1}_{((\ob{i},\ob{j}),\ob{k})}$$
and, for each $1$-morphism 
{\tiny
$$\bigoplus_r \left(\bigoplus_s(\F_{1rs},\F_{2rs})  ,\F_{3r}\right)\in ((\Aa\bxob\Bb)\bxob\Cc)(((\ob{i},\ob{j}),\ob{k}),((\ob{i}',\ob{j}'),\ob{k}')) = (\Aa(\ob{i},\ob{i}')\monob \Bb(\ob{j},\ob{j}'))\monob \Cc(\ob{k},\ob{k}'),$$
}
the isomorphism 
$$\left(\eta^{\ass}_{\Aa(\ob{i},\ob{i}'), \Bb(\ob{j},\ob{j}'), \Cc(\ob{k},\ob{k}')}\right)_{\bigoplus_r \left(\bigoplus_s(\F_{1rs},\F_{2rs})  ,\F_{3r}\right)}.$$
\item \label{ass2fun.4} There is a strong natural transformation $\varepsilon^{\asstwo_{\Aa,\Bb,\Cc}}\colon  \asstwo_{\Aa,\Bb,\Cc}\asstwo_{\Aa,\Bb,\Cc}^*  \to\mathbbm{1}_{\Aa\bxob(\Bb\bxob\Cc)}$ given, for objects $(\ob{i},(\ob{j},\ob{k}))\in \Aa\bxob(\Bb\bxob\Cc)$, by
$$\varepsilon^{\asstwo_{\Aa,\Bb,\Cc}}_{(\ob{i},(\ob{j},\ob{k}))} = \mathbbm{1}_{(\ob{i},(\ob{j},\ob{k}))}$$
and, for each 
{\tiny
$$\bigoplus_r \left( \F_{1r}, \bigoplus_s (\F_{2rs} ,\F_{3rs})\right)\in \Aa\bxob(\Bb\bxob\Cc)((\ob{i},(\ob{j},\ob{k})),(\ob{i}',(\ob{j}',\ob{k}'))) = \Aa(\ob{i},\ob{i}')\monob (\Bb(\ob{j},\ob{j}')\monob \Cc(\ob{k},\ob{k}')),$$
}
the isomorphism 
$$\left(\varepsilon^{\ass}_{\Aa(\ob{i},\ob{i}'), \Bb(\ob{j},\ob{j}')), \Cc(\ob{k},\ob{k}')}\right)_{\bigoplus_r \left( \F_{1r}, \bigoplus_s (\F_{2rs} ,\F_{3rs})\right)}$$
\item\label{ass2fun.5} The strong natural transformations satisfy the triangle identities strictly, 
 i.e. the diagrams
\[\xymatrix{
\asstwo_\Cc\ar^{\asstwo_\Cc\star \eta^{\asstwo_\Cc}}[rr]\ar@{=}[drr]&&\asstwo_\Cc\asstwo_\Cc^*\asstwo_\Cc\ar^{\varepsilon^{\asstwo_\Cc}\star\asstwo_\Cc}[d]\\
&&\asstwo_\Cc
}
\xymatrix{
\asstwo_\Cc^*\ar^{ \eta^{\asstwo_\Cc}\star \asstwo_\Cc^*}[rr]\ar@{=}[drr]&&\asstwo_\Cc^*\asstwo_\Cc\asstwo_\Cc^*\ar^{\asstwo_\Cc^*\star\varepsilon^{\asstwo_\Cc}}[d]\\
&&\asstwo_\Cc^*
}
\]
commute strictly.

 in the sense that their appropriate compositions are the identity strong natural transformation on $\asstwo_{\Aa,\Bb,\Cc}$ and $\asstwo^*_{\Aa,\Bb,\Cc}$, respectively.
\end{enumerate}
\end{prop}

\proof
\eqref{ass2fun.1} Strict unitality is obvious from the definition, so we need to give an isomorphism
{\tiny
\begin{equation*}\begin{split}
(\asstwo_{\Aa,\Bb,\Cc}&)_{(\ob{i}'',(\ob{j}'',\ob{k}'')),(\ob{i}',(\ob{j}',\ob{k}'))} \left(\bigoplus_p \left(\bigoplus_q(\G_{1pq},\G_{2pq})  ,\G_{3p}\right) \right)(\asstwo_{\Aa,\Bb,\Cc})_{(\ob{i}',(\ob{j}',\ob{k}')),(\ob{i},(\ob{j},\ob{k}))}  \left(\bigoplus_r \left(\bigoplus_s(\F_{1rs},\F_{2rs})  ,\F_{3r}\right) \right)\\
&\to(\asstwo_{\Aa,\Bb,\Cc})_{(\ob{i}'',(\ob{j}'',\ob{k}'')),(\ob{i},(\ob{j},\ob{k}))} \left( \left(\bigoplus_p \left(\bigoplus_q(\G_{1pq},\G_{2pq})  ,\G_{3p}\right) \right) \left(\bigoplus_r \left(\bigoplus_s(\F_{1rs},\F_{2rs})  ,\F_{3r}\right) \right) \right).
\end{split}\end{equation*}
}
To this end, using the definition of $\asstwo_{\Aa,\Bb,\Cc}$ on morphism categories, we explicitly compute
{\tiny
\begin{equation*}\begin{split}
\asstwo_{\Aa(\ob{i}',\ob{i}''),\Bb(\ob{j}',\ob{j}''),\Cc(\ob{k}',\ob{k}'')}&\left(\bigoplus_p \left(\bigoplus_q(\G_{1pq},\G_{2pq})  ,\G_{3p}\right) \right)\asstwo_{\Aa(\ob{i},\ob{i}'),\Bb(\ob{j},\ob{j}'),\Cc(\ob{k},\ob{k}')} \left(\bigoplus_r \left(\bigoplus_s(\F_{1rs},\F_{2rs})  ,\F_{3r}\right) \right)\\
&=\left( \bigoplus_p \bigoplus_q\left(\G_{1pq},(\G_{2pq}  ,\G_{3p})\right)  \right)\left( \bigoplus_r \bigoplus_s\left(\F_{1rs},(\F_{2rs} ,\F_{3r})\right)  \right)\\
&=\bigoplus_p \bigoplus_q \bigoplus_r \bigoplus_s\left( \G_{1pq} \F_{1rs}, (\G_{2pq} \F_{2rs}  ,\G_{3p} \F_{3r} ) \right)\\
\end{split}\end{equation*}
}
and
{\tiny
\begin{equation*}\begin{split}
\asstwo_{\Aa(\ob{i},\ob{i}''),\Bb(\ob{j},\ob{j}''),\Cc(\ob{k},\ob{k}'')} &\left( \left(\bigoplus_p \left(\bigoplus_q(\G_{1pq},\G_{2pq})  ,\G_{3p}\right) \right) \left(\bigoplus_r \left(\bigoplus_s(\F_{1rs},\F_{2rs})  ,\F_{3r}\right) \right) \right)\\
&=\asstwo_{\Aa(\ob{i},\ob{i}''),\Bb(\ob{j},\ob{j}''),\Cc(\ob{k},\ob{k}'')} 
\left( \bigoplus_p\bigoplus_r \left(\bigoplus_q(\G_{1pq},\G_{2pq}) \bigoplus_s(\F_{1rs},\F_{2rs}) ,\G_{3p}\F_{3r}\right) \right)\\
&=\asstwo_{\Aa(\ob{i},\ob{i}''),\Bb(\ob{j},\ob{j}''),\Cc(\ob{k},\ob{k}'')} 
\left( \bigoplus_p\bigoplus_r \left(\bigoplus_q\bigoplus_s(\G_{1pq}\F_{1rs},\G_{2pq}\F_{2rs}),\G_{3p}\F_{3r}\right) \right)\\
&=\bigoplus_p\bigoplus_r \bigoplus_q\bigoplus_s\left(\G_{1pq}\F_{1rs},(\G_{2pq}\F_{2rs},\G_{3p}\F_{3r})\right).
\end{split}\end{equation*}
}
There is an obvious isomorphism, given by resummation, which is clearly natural and satisfies lax associativity, and defines the associativity constraint $\asstwo_{\Aa,\Bb,\Cc}^2$.

\eqref{ass2fun.2} Mutatis mutandis the proof of \eqref{ass2fun.1}.

\eqref{ass2fun.3} We need to check lax naturality. 
Using 
\begin{equation*}\begin{split}
\asstwo_{\Aa,\Bb,\Cc}\asstwo_{\Aa,\Bb,\Cc}^*\left(\bigoplus_r \left( \F_{1r}, \bigoplus_s (\F_{2rs} ,\F_{3rs})\right)\right) =\bigoplus_r \bigoplus_s\left( \F_{1r},  (\F_{2rs} ,\F_{3rs})\right)
\end{split}\end{equation*}
and 
\begin{equation*}\begin{split}
\left(\bigoplus_p \left( \G_{1p}, \bigoplus_q (\G_{2pq} ,\G_{3pq})\right)\right)\left(\bigoplus_r \left( \F_{1r}, \bigoplus_s (\F_{2rs} ,\F_{3rs})\right)\right)\\
\qquad= \bigoplus_p\bigoplus_r \left(  \G_{1p} \F_{1r}, \left(   \bigoplus_q  \bigoplus_r (\G_{2pq}\F_{2rs},\G_{3pq} \F_{3rs}) \right)  \right)
\end{split}\end{equation*}
is is easy to see that
{\tiny

$$\xymatrix{
\left(\bigoplus_p \left( \G_{1p}, \bigoplus_q (\G_{2pq} ,\G_{3pq})\right)\right)\left(\bigoplus_r \left( \F_{1r}, \bigoplus_s (\F_{2rs} ,\F_{3rs})\right)\right)\ar@{=}[r]\ar^{\eta\hcirc\eta}[d]& \bigoplus_p\bigoplus_r \left(  \G_{1p} \F_{1r}, \left(   \bigoplus_q  \bigoplus_s(\G_{2pq}\F_{2rs},\G_{3pq} \F_{3rs}) \right)  \right)\ar^\eta[dd]\\
\left(\bigoplus_p\bigoplus_q  \left( \G_{1p}, (\G_{2pq} ,\G_{3pq})\right)\right)\left(\bigoplus_r\bigoplus_s \left( \F_{1r},  (\F_{2rs} ,\F_{3rs})\right)\right) \ar@{=}[d]& \\
\bigoplus_p \bigoplus_q \bigoplus_r  \bigoplus_s \left(  \G_{1p} \F_{1r},  (\G_{2pq}\F_{2rs},\G_{3pq} \F_{3rs})  \right) \ar^{(\asstwo_{\Aa,\Bb,\Cc}\asstwo_{\Aa,\Bb,\Cc}^*)^2}[r]&\bigoplus_p\bigoplus_r  \bigoplus_q  \bigoplus_s\left(  \G_{1p} \F_{1r},  (\G_{2pq}\F_{2rs},\G_{3pq} \F_{3rs}) \right) 
}$$
}

commutes, as required.

\eqref{ass2fun.4} Mutatis mutandis the proof of \eqref{ass2fun.3}.

\eqref{ass2fun.5} This follows from the fact that the same holds for the components $\eta^{\ass}_{\Aa(\ob{i},\ob{i}'), \Bb(\ob{j},\ob{j}')), \Cc(\ob{k},\ob{k}')}$, $\varepsilon^{\ass}_{\Aa(\ob{i},\ob{i}'), \Bb(\ob{j},\ob{j}')), \Cc(\ob{k},\ob{k}')}$.
\endproof

\begin{prop}
Let $\Aa,\Bb,\Cc, \Aa',\Bb',\Cc'$ be $2$-categories in $\twoAk$ and consider strictly unital $\Bbbk$-linear pseudofunctors $\Ff\colon \Aa\to\Aa', \Gg\colon \Bb\to\Bb', \Hh\colon \Cc\to\Cc'$.
Then \[(\Ff\bxone(\Gg\bxone\Hh)) \asstwo_{\Aa,\Bb,\Cc} =  \asstwo_{\Aa',\Bb',\Cc'} ((\Ff\bxone\Gg)\bxone\Hh).\]
\end{prop}

\proof
For a $1$-morphism $\bigoplus_r \left(\bigoplus_s(\F_{1rs},\F_{2rs})  ,\F_{3r}\right)$ in $(\Aa\bxob\Bb)\bxob\Cc $, we compute
\begin{equation*}\begin{split}
(\Ff\bxone(\Gg\bxone\Hh))& \asstwo_{\Aa,\Bb,\Cc}\left(\bigoplus_r \left(\bigoplus_s(\F_{1rs},\F_{2rs})  ,\F_{3r}\right)\right)\\
&= (\Ff\bxone(\Gg\bxone\Hh))\bigl(\bigoplus_r \bigoplus_s\left(\F_{1rs},(\F_{2rs}  ,\F_{3r})\right)\bigr)\\
&=\bigoplus_r \bigoplus_s\left(\Ff(\F_{1rs}),(\Gg(\F_{2rs})  ,\Hh(\F_{3r}))\right)
\end{split}\end{equation*}
and 
\begin{equation*}\begin{split}
 \asstwo_{\Aa',\Bb',\Cc'}& ((\Ff\bxone\Gg)\bxone\Hh)\left(\bigoplus_r \left(\bigoplus_s(\F_{1rs},\F_{2rs})  ,\F_{3r}\right)\right)\\
&= \asstwo_{\Aa',\Bb',\Cc'}\left(\bigoplus_r \left(\bigoplus_s(\Ff(\F_{1rs}),\Gg(\F_{2rs}))  ,\Hh(\F_{3r})\right)\right)\\
&=\bigoplus_r \bigoplus_s\left(\Ff(\F_{1rs}),(\Gg(\F_{2rs})  ,\Hh(\F_{3r}))\right).
\end{split}\end{equation*}
A similar computation shows that they agree on $2$-morphisms.
\endproof

\subsection{Higher coherences}\label{boxhigher}

\begin{prop}\label{pentagtwo}
For each quadruple $\Aa, \Bb, \Cc, \Dd$ of objects in $\twoAk$, there are equalities of pseudofunctors
$$
(\mathbbm{1}_{\Aa}\bxone\asstwo_{\Bb, \Cc, \Dd}) \asstwo_{\Aa, \Bb\bxob \Cc, \Dd}( \asstwo_{\Aa, \Bb, \Cc}\bxone\mathbbm{1}_{\Dd}) = \asstwo_{\Aa, \Bb, \Cc\bxob \Dd}  \asstwo_{\Aa\bxob\Bb, \Cc, \Dd}.
$$
\end{prop}

\proof
The identity of the pseudofunctors applied to objects in $\Aa, \Bb, \Cc, \Dd$ is immediate, while the identity of the component functors follows from Proposition \ref{pentagone}. The associativity constraints of the pseudofunctors are the canonical isomorphisms given by permutation of direct summands and hence agree as well.
\endproof

\begin{prop}\label{2unitortwo}
For any objects $\Cc,\Dd$, we have invertible strong natural transformations of pseudofunctors
$$(\mathbbm{1}_{\Cc}\bxone \luntwo_{\Dd})\, \asstwo_{\Cc,\Vv, \Dd} \,(\runtwo^*_{\Cc}\bxone\mathbbm{1}_{\Dd} )\xrightarrow{\mutwo_{\Cc,\Dd}}\mathbbm{1}_{\Cc\bxob\Dd}$$
$$ \luntwo_{\Cc}\bxone \mathbbm{1}_{\Dd}
 \xrightarrow{\lambdatwo_{\Cc,\Dd}}\luntwo_{\Cc\otimes_0\Dd}\, \asstwo_{\Vv,\Cc,\Dd } $$
$$\mathbbm{1}_{\Cc}\bxone \runtwo^*_{\Dd}
 \xrightarrow{\rhotwo_{\Cc,\Dd}} \ass_{\Cc,\Dd,\Vv } \,\runtwo^*_{\Cc\bxob\Dd} .$$
\end{prop}

\proof
On objects $\ob{i}\in \Cc,\ob{j}\in\Dd$, the component $1$-cells are the respective identity $1$-morphisms, that is, $\mathbbm{1}_{(\ob{i}, \ob{j})}$ for $\mutwo_{\Cc,\Dd}, \rhotwo_{\Cc,\Dd}$ and $\mathbbm{1}_{(\bullet, \ob{i}, \ob{j})}$ for $\lambdatwo_{\Cc,\Dd}$. For $1$-morphisms $\bigoplus_{i}(\F_i,\G_i)$ in $(\Cc\bxob\Dd)((\ob{i}, \ob{j}), (\ob{i}', \ob{j}'))$, the required compatibility isomorphisms for $\mutwo$ and $\rhotwo$ are given by 
\[(\mutwo_{\Cc,\Dd})_{\bigoplus_{i}(\F_i,\G_i)} = \mathbbm{1}_{\bigoplus_{i}(\F_i,\G_i)} \qquad (\rhotwo_{\Cc,\Dd})_{\bigoplus_{i}(\F_i,\G_i)} = \mathbbm{1}_{\bigoplus_{i}(\F_i,\G_i)}\]
which clearly define a strong (even strict) natural transformation. 

For $1$-morphisms $ \bigoplus_{j}(\bigoplus_{i}(\mathbbm{1}_\bullet^{\oplus k_{ij}}, \F_{i,j}), \G_j)    $ in $((\Vv\bxob\Cc)\bxob\Dd)(((\bullet,\ob{i}), \ob{j}), ((\bullet, \ob{i}'), \ob{j}'))$, the compatibility morphisms are given by 
$$(\lambdatwo_{\Cc,\Dd})_{\bigoplus_{j}(\bigoplus_{i}(\mathbbm{1}_\bullet^{\oplus k_{ij}}, \F_{i,j}), \G_j)  }
= (\lambda_{\Cc(\ob{i}, \ob{i}'),\Dd(\ob{j}, \ob{j}')})_{\bigoplus_{j}(\bigoplus_{i} (\mathbbm{1}_\bullet^{\oplus k_{ij}}, \F_{i,j}), \G_j)  }.$$
To check that this defines a strong natural transformation, we use the fact that $\mathbbm{1}_{\Cc}\bxone \runtwo^*_{\Dd}$ is a $2$-functor as well as the computation
\begin{equation*}\begin{split}
\bigoplus_{j}\left(\bigoplus_{i} (\mathbbm{1}_\bullet^{\oplus k_{ij}}, \F_{i,j}), \G_j\right) 
\bigoplus_{b}\left(\bigoplus_{a} (\mathbbm{1}_\bullet^{\oplus c_{ab}}, \F'_{a,b}), \G'_b\right)    \\
\qquad = \bigoplus_{j}\bigoplus_b\left(\bigoplus_{i} \bigoplus_{a} (\mathbbm{1}_\bullet,  \F_{i,j}\F'_{a,b})^{\oplus k_{ij}c_{ab}},\G_j\G'_b    \right)
\end{split}\end{equation*}
to check that the horizontal composition of the relevant components of $\lambdatwo_{\Cc,\Dd}$ maps the leftt hand-side to
$$\bigoplus_{j}\bigoplus_{i} (\F_{i,j}, \G_j)^{\oplus k_{ij}}
\bigoplus_{b}\bigoplus_{a} ( \F'_{a,b}, \G'_b)^{\oplus c_{ab}} 
= \bigoplus_{j}\bigoplus_{i}\bigoplus_{b}\bigoplus_{a}(\F_{i,j} \F'_{a,b}, \G_j \G'_b)^{\oplus k_{ij}c_{ab}}$$
while the right hand-side, under the relevant component of $\lambdatwo_{\Cc,\Dd}$ gets mapped to 
$$\bigoplus_{j}\bigoplus_{b}\bigoplus_{i}\bigoplus_{a}(\F_{i,j} \F'_{a,b}, \G_j \G'_b)^{\oplus k_{ij}c_{ab}}$$
and applying the relevant component of the isomorphism $(\luntwo_{\Cc\otimes_0\Dd}\, \asstwo_{\Vv,\Cc,\Dd })^2$, which is induced by the component  $\asstwo_{\Vv,\Cc,\Dd }^2$ since $\luntwo_{\Cc\otimes_0\Dd}$ is a $2$-functor, to the former, we see that we precisely obtain the latter, as required for a strong natural transformation.

\endproof

%

\section{Comonoidal structure on $2$-categories}\label{sec:comon}

\subsection{Some auxiliary strong transformations}

Recall the $2$-functors $\luntwo_\Cc\colon \Vv \bxob \Cc\to \Cc$ and  $\runtwo_\Cc\colon \Cc\bxob \Vv\to \Cc$ from Propositions \ref{unitortwopartsleft} and \ref{unitortwopartsright} respectively, and of $\asstwo_{\Aa,\Bb,\Cc}$ from Proposition~\ref{ass2fun}. Recall further the strong transformations $\luntwo_\Ff\colon \luntwo_{\Dd}(\mathbbm{1}_{\Vv}\bxone \Ff) \to \Ff \luntwo_\Cc$, $\runtwo_\Ff\colon \Ff\runtwo_\Cc \to \runtwo_\Dd(\Ff\bxone \mathbbm{1}_{\Vv})$ from Proposition \ref{unitortwocommute} and 
$\mutwo_{\Cc,\Dd}, \lambdatwo_{\Cc,\Dd}$ and $\rhotwo_{\Cc,\Dd}$ from Proposition \ref{2unitortwo}.

In the following, we will need slight variations of $\mutwo_{\Cc,\Dd}$ and $\rhotwo_{\Cc,\Dd}$, which we record here.


\subsubsection{}\label{phimdef}
There is an invertible strong natural transformation 
\[\tilde\mutwo_{\Cc,\Cc}\colon \runtwo_\Cc\bxone \mathbbm{1}_\Cc \to (\mathbbm{1}_\Cc \bxone\luntwo_{\Cc})\asstwo_{\Cc,\Vv,\Cc}\]
of pseudofunctors from 
$(\Cc\bxob\Vv)\bxob \Cc$ to $\Cc\bxob\Cc$, given by
\[\tilde\mutwo_{\Cc,\Cc} = ( ( (\mathbbm{1}_\Cc \bxone\luntwo_{\Cc})\asstwo_{\Cc,\Vv,\Cc})\star \varepsilon^{\runtwo_{\Cc}})(\mutwo_{\Cc,\Cc}^{-1}\star (\runtwo_\Cc\bxone \mathbbm{1}_\Cc)).\]

Explicitly, its $1$-cell component is the identity (since both pseudofunctors map an object $((\ob{i},\bullet ),\ob{j})$ to $(\ob{i},\ob{j})$) and, for a $1$-morphism $\bigoplus_s(\bigoplus_t(\F_{st}, \mathbbm{1}_\bullet^{\oplus a_{st}}), \G_{s} )$ in $(\Cc\bxob\Vv)\bxob \Cc$,  its $2$-cell component is given by the canonical isomorphism between 
\[ (\runtwo_\Cc\bxone \mathbbm{1}_\Cc)(\bigoplus_s(\bigoplus_t(\F_{st}, \mathbbm{1}_\bullet^{\oplus a_{st}}), \G_{s} )) 
= \bigoplus_s(\bigoplus_t \F_{st}^{\oplus a_{st}} ,\G_{s} ) \]
and 
\begin{equation*}\begin{split}  
(\mathbbm{1}_\Cc\bxone \luntwo_{\Cc})\asstwo_{\Cc,\Vv,\Cc}    (\bigoplus_s(\bigoplus_t(\F_{st}, \mathbbm{1}_\bullet^{\oplus a_{st}}), \G_{s} ))
& =  (\mathbbm{1}_\Cc\bxone \luntwo_{\Cc})  (\bigoplus_s\bigoplus_t(\F_{st} ,(\mathbbm{1}_\bullet^{\oplus a_{st}},\G_{s} ))) \\&= \bigoplus_s\bigoplus_t(\F_{st} ,\G_{s}^{\oplus a_{st}}). 
\end{split}\end{equation*}

\subsubsection{}\label{phirdef}
There is an invertible strong natural transformation 
\[\tilde\rhotwo_{\Cc,\Cc}\colon \runtwo_{\Cc\bxob\Cc} \to (\mathbbm{1}_\Cc\bxone \runtwo_\Cc)\asstwo_{\Cc,\Cc,\Vv}\]
between pseudofunctors from $(\Cc\bxob\Cc)\bxob \Vv$ to $\Cc\bxob\Cc$ given by
\[  \tilde\rhotwo_{\Cc,\Cc} = (  ( (\mathbbm{1}_\Cc\bxone \runtwo_\Cc)\asstwo_{\Cc,\Cc,\Vv})\star \varepsilon^{\runtwo_{\Cc\bxob\Cc}}   )  ( (\mathbbm{1}_\Cc\bxone \runtwo_\Cc)\star \rhotwo_{\Cc,\Cc} \star  \runtwo_{\Cc\bxone\Cc}  )(   (\mathbbm{1}_{\mathbbm{1}_\Cc}\bxtwo \eta^{\runtwo_\Cc})\star \runtwo_{\Cc\bxone\Cc}) .  \]
Explicitly, its $1$-cell component is the identity (since both pseudofunctors map an object $((\ob{i},\ob{j}), \bullet)$ to $(\ob{i},\ob{j})$) and, for a $1$-morphism $\bigoplus_s(\bigoplus_t(\F_{st}, \G_{st}),  \mathbbm{1}_\bullet^{\oplus a_s})$ in $(\Cc\bxob\Cc)\bxob \Vv$, its $2$-cell component is given by the canonical isomorphism between 
\[ \runtwo_{\Cc\bxob\Cc}(\bigoplus_s(\bigoplus_t(\F_{st}, \G_{st}),  \mathbbm{1}_\bullet^{\oplus a_s})) = \bigoplus_s(\bigoplus_t(\F_{st}, \G_{st})^{\oplus a_s}) \]
and 
\begin{equation*}\begin{split}  
(\mathbbm{1}_\Cc\bxone \runtwo_\Cc)\asstwo_{\Cc,\Cc,\Vv}   (\bigoplus_s(\bigoplus_t(\F_{st}, \G_{st}),  \mathbbm{1}_\bullet^{\oplus a_s}))& =  (\mathbbm{1}_\Cc\bxone \runtwo_\Cc) (\bigoplus_s\bigoplus_t(\F_{st}, (\G_{st},  \mathbbm{1}_\bullet^{\oplus a_s}))) \\&= \bigoplus_s\bigoplus_t(\F_{st}, \G_{st}^{\oplus a_s}). 
\end{split}\end{equation*}

\subsection{Definition of a comonoidal structure}

Let $\Cc\in \twoAk$ be small. We will say $\Cc$ has a {\em comonoidal structure} if there are strictly unital $\Bbbk$-linear pseudofunctors $\Delta\colon \Cc\to\Cc\bxob\Cc$ and $\Theta\colon \Cc\to \Vv$ together with adjoint equivalences 
$$\delta\colon\asstwo_{\Cc,\Cc,\Cc}(\Delta\bxone  \mathbbm{1}_{\Cc})\Delta \rightleftarrows ( \mathbbm{1}_{\Cc}\bxone \Delta)\Delta \colon \delta^*,$$
$$\theta^\ell\colon \luntwo_\Cc(\Theta\bxone  \mathbbm{1}_{\Cc})\Delta \rightleftarrows \mathbbm{1}_{\Cc} \colon \theta^{\ell*}$$
and $$\theta^r\colon \runtwo_\Cc(\mathbbm{1}_{\Cc}\bxone  \Theta)\Delta \rightleftarrows \mathbbm{1}_{\Cc} \colon \theta^{r*},$$ respectively, as depicted in the diagrams

$$\xymatrix{&&\Cc\ar^{\Delta}[dll]\ar^{\Delta}[drr] &&\\
\Cc\bxob\Cc \ar^{\Delta\bxone  \mathbbm{1}_{\Cc}}[dd] &&&\ar@2{->}@<5pt>^{\delta^*}[dll]&\Cc\bxob\Cc \ar_{ \mathbbm{1}_{\Cc}\bxone \Delta}[dd]\\
&\ar@2{->}@<5pt>^{\delta}[urr]&&&\\
(\Cc\bxob\Cc)\bxob\Cc\ar^{\asstwo_{\Cc,\Cc,\Cc}}[rrrr]&&&&\Cc\bxob(\Cc\bxob\Cc)}
$$


$$
\xymatrix{
\Cc\ar^{\Delta}[rrr]\ar_{\mathbbm{1}_\Cc}[rrrdd] &&&\Cc\bxob\Cc\ar^{\Theta\bxone  \mathbbm{1}_{\Cc}}[d] \ar@2{->}@<5pt>^{\theta^\ell}[dl]\\
&&\ar@2{->}@<5pt>^{\theta^{\ell*}}[ur]&\Vv\bxob\Cc\ar^{\luntwo_{\Cc}}[d]\\
&&&\Cc
}
 \quad\text{and}\quad
\xymatrix{
\Cc\ar^{\Delta}[rrr] \ar_{\mathbbm{1}_\Cc}[rrrdd] &&&\Cc\bxob\Cc \ar^{\mathbbm{1}_{\Cc}\bxone \Theta }[d]\ar@2{->}@<5pt>^{\theta^r}[dl]\\
&&\ar@2{->}@<5pt>^{\theta^{r*}}[ur]&\Cc\bxob\Vv \ar^{\runtwo_{\Cc}}[d]\\
&&&\Cc.
}
 $$
 We denote the units and counits  these adjoint equivalences by $\eta^{\delta}, \varepsilon^{\delta},\eta^{\theta^\ell}, \varepsilon^{\theta^\ell},\eta^{\theta^r}$ and $ \varepsilon^{\theta^r}$, respectively.
 
 Moreover, we require that
\begin{enumerate} 
\item\label{middlecond}
 the strong transformation defined by the pasting diagram
 \begin{align*}
\scalebox{0.8}{
\begin{tikzcd}[ampersand replacement=\&]
	\&\& \Cc \&\& \\
	{\Cc\bxone \Cc} \&\&\&\& {\Cc\bxone \Cc} \\
	\& {(\Cc\bxob \Cc )\bxob \Cc} \&\& {\Cc\bxob(\Cc\bxob \Cc )} \\
	\\
	\& {(\Cc\bxob \Vv )\bxob \Cc} \&\& {\Cc\bxob(\Vv\bxob \Cc )} \\
	\\
	{\Cc\bxob \Cc} \&\&\&\& {\mathcal{C}\boxtimes_0 \mathcal{C}}
	\arrow["\Delta"', from=1-3, to=2-1]
	\arrow["\Delta", from=1-3, to=2-5]
	\arrow["\delta", between={0.4}{0.6}, Rightarrow, from=2-1, to=2-5]
	\arrow["{{\Delta\bxone\mathbbm{1}_\Cc}}", from=2-1, to=3-2]
	\arrow[""{name=0, anchor=center, inner sep=0}, equals, from=2-1, to=7-1]
	\arrow["{{\mathbbm{1}_\Cc\bxone\Delta}}"', from=2-5, to=3-4]
	\arrow[""{name=1, anchor=center, inner sep=0}, equals, from=2-5, to=7-5]
	\arrow["{{\asstwo_{\Cc,\Cc,\Cc}}}", from=3-2, to=3-4]
	\arrow[""{name=2, anchor=center, inner sep=0}, "{{(\mathbbm{1}_\Cc\bxone\Theta)\bxone \mathbbm{1}_\Cc}}", from=3-2, to=5-2]
	\arrow[""{name=3, anchor=center, inner sep=0}, "{{\mathbbm{1}_\Cc\bxone(\Theta\bxone \mathbbm{1}_\Cc)}}"', from=3-4, to=5-4]
	\arrow[between={0.4}{0.6}, equals, from=5-2, to=3-4]
	\arrow["{{\asstwo_{\Cc,\Vv,\Cc}}}", from=5-2, to=5-4]
	\arrow["{{\runtwo_\Cc\bxone \mathbbm{1}_\Cc}}", from=5-2, to=7-1]
	\arrow["{{\mathbbm{1}_\Cc\bxone\luntwo_\Cc}}"', from=5-4, to=7-5]
	\arrow["{\tilde\mutwo_{\Cc,\Cc}}"{pos=0.7}, between={0.6}{0.8}, Rightarrow, from=7-1, to=5-4]
	\arrow[equals, from=7-1, to=7-5]
	\arrow["{{\theta^{r*}\bxtwo \mathbbm{1}_{\mathbbm{1}_\Cc}}}"{pos=0.6}, between={0.4}{0.6}, Rightarrow, from=0, to=2]
	\arrow["{{\mathbbm{1}_{\mathbbm{1}_\Cc}\bxtwo \theta^\ell}}"{pos=0.4}, between={0.4}{0.6}, Rightarrow, from=3, to=1]
\end{tikzcd}}
\end{align*}
is equal to the identity strong transformation on $\Delta$,
\item\label{leftcond} the two strong transformations from $\luntwo_{\Cc\bxob\Cc}\asstwo_{\Vv,\Cc,\Cc}  ((\Theta\bxone \mathbbm{1}_\Cc)\bxone \mathbbm{1}_\Cc)(\Delta\bxone\mathbbm{1}_\Cc)\Delta$ to $\Delta$ defined by the parts to the left, respectively right, of the highlighted blue pseudofunctor in the diagram 
 \begin{align*}
\scalebox{0.8}{
\begin{tikzcd}[ampersand replacement=\&]
	\&\& \Cc \&\&\&\& \\
	\& {\Cc\bxone \Cc} \&\& {\Cc\bxob \Cc} \\
	\\
	\& {(\Cc\bxob \Cc )\bxob \Cc} \&\& {\Cc\bxob(\Cc\bxob \Cc )} \& {\Vv\bxob \Cc} \&\& \Cc \\
	\\
	{\Cc\bxob \Cc} \& {(\Vv\bxob \Cc )\bxob \Cc} \&\& {\Vv\bxob(\Cc\bxob \Cc )} \\
	\& {\Vv\bxob(\Cc\bxob \Cc )} \\
	\\
	\& {\Cc\bxob \Cc} \&\& {\Cc\bxob \Cc}
	\arrow["\Delta"', color={rgb,255:red,92;green,92;blue,214}, from=1-3, to=2-2]
	\arrow["\Delta", from=1-3, to=2-4]
	\arrow[""{name=0, anchor=center, inner sep=0}, curve={height=-30pt}, equals, from=1-3, to=4-7]
	\arrow["\delta", between={0.4}{0.6}, Rightarrow, from=2-2, to=2-4]
	\arrow["{{\Delta\bxone\mathbbm{1}_{\Cc}}}", color={rgb,255:red,92;green,92;blue,214}, from=2-2, to=4-2]
	\arrow[""{name=1, anchor=center, inner sep=0}, curve={height=12pt}, equals, from=2-2, to=6-1]
	\arrow["{{\mathbbm{1}_\Cc\bxone\Delta}}"', from=2-4, to=4-4]
	\arrow[""{name=2, anchor=center, inner sep=0}, "{{\Theta \bxone \mathbbm{1}_{\Cc}}}", from=2-4, to=4-5]
	\arrow["{{\asstwo_{\Cc,\Cc,\Cc}}}", from=4-2, to=4-4]
	\arrow[""{name=3, anchor=center, inner sep=0}, "{{(\Theta\bxone\mathbbm{1}_\Cc)\bxone \mathbbm{1}_\Cc}}", color={rgb,255:red,92;green,92;blue,214}, from=4-2, to=6-2]
	\arrow[between={0.4}{0.6}, equals, from=4-4, to=4-5]
	\arrow["{{\Theta\bxone (\mathbbm{1}_\Cc\bxone\mathbbm{1}_\Cc)}}"', from=4-4, to=6-4]
	\arrow["{{\luntwo_\Cc}}"', from=4-5, to=4-7]
	\arrow["{{\mathbbm{1}_\Vv \bxone \Delta}}", from=4-5, to=6-4]
	\arrow["\Delta", curve={height=-30pt}, from=4-7, to=9-4]
	\arrow[curve={height=12pt}, equals, from=6-1, to=9-2]
	\arrow[between={0.4}{0.6}, equals, from=6-2, to=4-4]
	\arrow["{{\luntwo_\Cc\bxone \mathbbm{1}_\Cc}}"', from=6-2, to=6-1]
	\arrow["{{\asstwo_{\Vv,\Cc,\Cc}}}", from=6-2, to=6-4]
	\arrow["{{\asstwo_{\Vv,\Cc,\Cc}}}", color={rgb,255:red,92;green,92;blue,214}, from=6-2, to=7-2]
	\arrow[""{name=4, anchor=center, inner sep=0}, "{{\luntwo_{\Cc \bxob \Cc}}}"', from=6-4, to=9-4]
	\arrow["{\lambdatwo_{\Cc,\Cc}^{-1}}"{description}, between={0.1}{0.9}, Rightarrow, from=7-2, to=6-1]
	\arrow["{{\luntwo_{\Cc \bxob \Cc}}}", color={rgb,255:red,92;green,92;blue,214}, from=7-2, to=9-2]
	\arrow[between={0.4}{0.6}, equals, from=9-2, to=6-4]
	\arrow[equals, from=9-2, to=9-4]
	\arrow["{{\theta^{\luntwo}\bxtwo \mathbbm{1}_{\mathbbm{1}_{\Cc}}}}"'{pos=0.7}, between={0.3}{0.7}, equals, 2tail reversed, from=1, to=3]
	\arrow["{{\theta^{\luntwo}}}", between={0.4}{0.6}, Rightarrow, from=2, to=0]
	\arrow["{\luntwo_\Delta^{-1}}", between={0.4}{0.6}, Rightarrow, from=4, to=4-7]
\end{tikzcd}
}
\end{align*}
are equal,
\item\label{rightcond} the  two strong transformations from $\Delta$  to $(\mathbbm{1}_\Cc\bxone \runtwo_\Cc) ( \mathbbm{1}_\Cc\bxone (\mathbbm{1}_\Cc\bxone\Theta))(\mathbbm{1}_\Cc\bxone\Delta)\Delta $, defined by the parts to the left, respectively right, of the highlighted blue pseudofunctor in the diagram  \begin{align*}
\scalebox{0.8}{
\begin{tikzcd}[ampersand replacement=\&]
	\&\&\& \Cc \&\& \\
	\&\& {\Cc\bxob\Cc} \&\& {\Cc\bxob\Cc} \\
	\\
	\Cc \& {\Cc\bxob\Vv} \& {(\Cc\bxob\Cc)\bxob\Cc} \&\& {\Cc\bxob(\Cc\bxob\Cc)} \\
	\\
	{\Cc\bxob\Cc} \&\& {(\Cc\bxob\Cc)\bxob\Vv} \&\& {\Cc\bxob(\Cc\bxob\Vv)} \\
	\&\& {(\Cc\bxob\Cc)\bxob\Vv} \&\&\& {\Cc\bxob\Cc} \\
	\&\& {\Cc\bxob(\Cc\bxob\Vv)} \\
	\&\& {\Cc\bxob\Cc}
	\arrow["\Delta"', from=1-4, to=2-3]
	\arrow["\Delta", color={rgb,255:red,92;green,92;blue,214}, from=1-4, to=2-5]
	\arrow[""{name=0, anchor=center, inner sep=0}, curve={height=30pt}, equals, from=1-4, to=4-1]
	\arrow[""{name=1, anchor=center, inner sep=0}, "{{\mathbbm{1}_\Cc\bxone \Theta}}", from=2-3, to=4-2]
	\arrow[""{name=2, anchor=center, inner sep=0}, "{{\Delta\bxone \mathbbm{1}_{\Cc}}}", from=2-3, to=4-3]
	\arrow["{{\mathbbm{1}_{\Cc}\bxone \Delta}}"', color={rgb,255:red,92;green,92;blue,214}, from=2-5, to=4-5]
	\arrow[""{name=3, anchor=center, inner sep=0}, curve={height=-24pt}, equals, from=2-5, to=7-6]
	\arrow["\Delta"', from=4-1, to=6-1]
	\arrow["{\runtwo_\Delta}", between={0.4}{0.6}, Rightarrow, from=4-1, to=6-3]
	\arrow["{{\runtwo_\Cc}}", from=4-2, to=4-1]
	\arrow[between={0.4}{0.6}, equals, from=4-2, to=4-3]
	\arrow["{{\Delta\bxone \mathbbm{1}_{\Vv}}}", from=4-2, to=6-3]
	\arrow["{{\asstwo_{\Cc,\Cc,\Cc}}}"', from=4-3, to=4-5]
	\arrow["{{\mathbbm{1}_{\Cc\bxob\Cc}\bxone \Theta}}", from=4-3, to=6-3]
	\arrow["{{\mathbbm{1}_\Cc\bxone(\mathbbm{1}_\Cc\bxone \Theta)}}"', color={rgb,255:red,92;green,92;blue,214}, from=4-5, to=6-5]
	\arrow[curve={height=18pt}, equals, from=6-1, to=9-3]
	\arrow[between={0.4}{0.6}, equals, from=6-3, to=4-5]
	\arrow["{{\runtwo_{\Cc\bxob\Cc}}}"', from=6-3, to=6-1]
	\arrow["{{\asstwo_{\Cc,\Cc,\Vv}}}"', from=6-3, to=6-5]
	\arrow[equals, from=6-3, to=7-3]
	\arrow["{{\mathbbm{1}_\Cc\bxone \runtwo_\Cc}}"', color={rgb,255:red,92;green,92;blue,214}, from=6-5, to=7-6]
	\arrow[""{name=4, anchor=center, inner sep=0}, "{{\asstwo_{\Cc,\Cc,\Vv}}}", from=7-3, to=8-3]
	\arrow[""{name=5, anchor=center, inner sep=0}, "{{\mathbbm{1}_\Cc\bxone \runtwo_\Cc}}", from=8-3, to=9-3]
	\arrow[curve={height=12pt}, equals, from=9-3, to=7-6]
	\arrow["{{(\theta^r)^*}}"', between={0.4}{0.6}, Rightarrow, from=0, to=1]
	\arrow["\delta", between={0.4}{0.6}, Rightarrow, from=2, to=2-5]
	\arrow["{{\mathbbm{1}_{\Cc}\bxone (\theta^r)*}}"{pos=0.6}, between={0.4}{0.6}, Rightarrow, from=3, to=6-5]
	\arrow["{\tilde\rhotwo_{\Cc,\Cc}}", between={0.4}{0.6}, Rightarrow, from=6-1, to=4]
	\arrow[between={0.4}{0.6}, equals, from=5, to=6-5]
\end{tikzcd}
}
\end{align*}
are equal, and
\item\label{pentacond} the two strong transformations
 \begin{align*}
\scalebox{0.6}{
\begin{tikzcd}[ampersand replacement=\&]
	\&\&\& \Cc \&\&\& \\
	\\
	\& {\Cc\bxob\Cc} \&\&\&\& {\Cc\bxob\Cc} \\
	\\
	{(\Cc\bxob\Cc)\bxob \Cc} \&\& {(\Cc\bxob\Cc)\bxob \Cc} \&\& {\Cc\bxob(\Cc\bxob\Cc)} \&\& {\Cc\bxob(\Cc\bxob\Cc)} \\
	\\
	{((\Cc\bxob\Cc)\bxob\Cc)\bxob \Cc} \&\& {(\Cc\bxob(\Cc\bxob\Cc))\bxob\Cc} \&\& {\Cc\bxob ((\Cc\bxob\Cc)\bxob\Cc)} \&\& {\Cc\bxob (\Cc\bxob(\Cc\bxob\Cc))} \\
	\\
	\\
	\\
	\&\&\&\&\&\& {}
	\arrow["\Delta"', from=1-4, to=3-2]
	\arrow["\Delta", from=1-4, to=3-6]
	\arrow["{{\Delta\bxone\mathbbm{1}_\Cc}}"', from=3-2, to=5-1]
	\arrow[""{name=0, anchor=center, inner sep=0}, "{{\Delta\bxone\mathbbm{1}_\Cc}}", from=3-2, to=5-3]
	\arrow["{{\mathbbm{1}_\Cc\bxone \Delta}}"', from=3-6, to=5-5]
	\arrow["{{\mathbbm{1}_\Cc\bxone \Delta}}", from=3-6, to=5-7]
	\arrow[""{name=1, anchor=center, inner sep=0}, "{{(\Delta\bxone\mathbbm{1}_\Cc)\bxone\mathbbm{1}_\Cc}}"', from=5-1, to=7-1]
	\arrow["{{\asstwo_{\Cc,\Cc,\Cc}}}", from=5-3, to=5-5]
	\arrow["{{(\mathbbm{1}_\Cc\bxone \Delta)\bxone\mathbbm{1}_\Cc}}"', from=5-3, to=7-3]
	\arrow[""{name=2, anchor=center, inner sep=0}, "{{\mathbbm{1}_\Cc\bxone (\Delta\bxone\mathbbm{1}_\Cc)}}", from=5-5, to=7-5]
	\arrow["{{\mathbbm{1}_\Cc\bxone(\mathbbm{1}_\Cc\bxone \Delta)}}", from=5-7, to=7-7]
	\arrow["{{\asstwo_{\Cc,\Cc,\Cc}\bxone\mathbbm{1}_\Cc}}", from=7-1, to=7-3]
	\arrow[between={0.4}{0.6}, equals, from=7-3, to=5-5]
	\arrow["{{\asstwo_{\Cc,\Cc\bxob\Cc,\Cc}}}", from=7-3, to=7-5]
	\arrow["{{\mathbbm{1}_{\Cc}\bxone\asstwo_{\Cc,\Cc,\Cc}}}", from=7-5, to=7-7]
	\arrow["\delta", between={0.4}{0.6}, Rightarrow, from=0, to=3-6]
	\arrow["{{\delta\bxtwo\mathbbm{1}_{\mathbbm{1}_\Cc}}}", between={0.4}{0.6}, Rightarrow, from=1, to=5-3]
	\arrow["{{\mathbbm{1}_{\mathbbm{1}_\Cc}\bxtwo \delta}}", between={0.4}{0.6}, Rightarrow, from=2, to=5-7]
\end{tikzcd}}\end{align*}
and 
\begin{align*}
\scalebox{0.6}{
\begin{tikzcd}[ampersand replacement=\&]
	\&\&\& \Cc \&\&\& \\
	\\
	\& {\Cc\bxob\Cc} \&\& {\Cc\bxob\Cc} \&\& {\Cc\bxob\Cc} \\
	\\
	{(\Cc\bxob\Cc)\bxob \Cc} \&\& {\Cc\bxob(\Cc\bxob\Cc)} \&\& {(\Cc\bxob\Cc)\bxob \Cc} \&\& {\Cc\bxob(\Cc\bxob \Cc)} \\
	\\
	{((\Cc\bxob\Cc)\bxob\Cc)\bxob\Cc} \&\&\& {(\Cc\bxob \Cc)\bxob(\Cc\bxob\Cc)} \&\&\& {\Cc \bxob(\Cc\bxob(\Cc\bxob\Cc))}
	\arrow["\Delta"', from=1-4, to=3-2]
	\arrow["\Delta", from=1-4, to=3-4]
	\arrow["\Delta"', from=1-4, to=3-6]
	\arrow[""{name=0, anchor=center, inner sep=0}, "{{\Delta\bxone\mathbbm{1}_\Cc}}", from=3-2, to=5-1]
	\arrow["{{\mathbbm{1}_\Cc\bxone \Delta}}"', from=3-4, to=5-3]
	\arrow[""{name=1, anchor=center, inner sep=0}, "{{\Delta\bxone\mathbbm{1}_\Cc}}", from=3-4, to=5-5]
	\arrow["{{\mathbbm{1}_\Cc\bxone \Delta}}"', from=3-6, to=5-7]
	\arrow["{{\asstwo_{\Cc,\Cc,\Cc}}}", from=5-1, to=5-3]
	\arrow["{{(\Delta\bxone\mathbbm{1}_\Cc)\bxone\mathbbm{1}_\Cc}}"', from=5-1, to=7-1]
	\arrow[between={0.4}{0.6}, equals, from=5-3, to=5-5]
	\arrow["{{\Delta\bxone (\mathbbm{1}_\Cc\bxone\mathbbm{1}_\Cc)}}", from=5-3, to=7-4]
	\arrow["{{\asstwo_{\Cc,\Cc,\Cc}}}", from=5-5, to=5-7]
	\arrow["{{(\mathbbm{1}_\Cc\bxone\mathbbm{1}_\Cc)\bxone \Delta}}", from=5-5, to=7-4]
	\arrow["{{\mathbbm{1}_\Cc\bxone(\mathbbm{1}_\Cc\bxone \Delta)}}", from=5-7, to=7-7]
	\arrow[between={0.4}{0.6}, equals, from=7-1, to=5-3]
	\arrow["{{\asstwo_{\Cc\bxob\Cc,\Cc,\Cc}}}", from=7-1, to=7-4]
	\arrow[between={0.4}{0.6}, equals, from=7-4, to=5-7]
	\arrow["{{\asstwo_{\Cc,\Cc,\Cc\bxob\Cc}}}", from=7-4, to=7-7]
	\arrow["\delta", between={0.4}{0.6}, Rightarrow, from=0, to=3-4]
	\arrow["\delta", between={0.4}{0.6}, Rightarrow, from=1, to=3-6]
\end{tikzcd}}
\end{align*}
coincide, where the sources agree by Proposition \ref{pentagtwo}.
\end{enumerate}

Note that all pasting diagrams in the definition indeed uniquely define strong natural transformations by Lemma \ref{le:partialinterchange}.

\subsection{Immediate properties.}

\subsubsection{}
 Notice that, in particular, using the notation $\Delta(\ob{i}) = (\ob{i}_\ell,\ob{i}_r)$, part of the data of $\theta^\ell, \theta^{\ell*},\theta^r,\theta^{r*}$ are, for each $\ob{i}\in \Cc$, $1$-morphisms $$\theta^\ell_{\ob{i}}\in \Cc(\ob{i}_r,\ob{i}),\quad\theta^{\ell*}_{\ob{i}}\in \Cc(\ob{i},\ob{i}_r),\quad \theta^r_{\ob{i}}\in \Cc(\ob{i}_\ell,\ob{i}),\quad\theta^{r*}_{\ob{i}}\in \Cc(\ob{i},\ob{i}_\ell),$$ respectively. 
 
 Likewise, since $\asstwo_{\Cc,\Cc,\Cc}(\Delta\bxone \mathbbm{1}_\Cc)\Delta (\ob{i}) = (\ob{i}_{\ell\ell},(\ob{i}_{\ell r},\ob{i}_{r}))$ and
 $(\mathbbm{1}_\Cc\bxone \Delta)\Delta (\ob{i})  = (\ob{i}_{\ell},(\ob{i}_{r\ell},\ob{i}_{rr}))$, part of the data of $\delta$, consists of, for each $\ob{i}\in \Cc$, a $1$-morphism 
 $$
 \delta_{\ob{i}}\in (\Cc\bxob(\Cc\bxob\Cc))((\ob{i}_{\ell\ell},(\ob{i}_{\ell r},\ob{i}_{r}))  ,  (\ob{i}_{\ell},(\ob{i}_{r\ell},\ob{i}_{rr})) ) = 
 \Cc(\ob{i}_{\ell\ell}, \ob{i}_{\ell})\monob (\Cc(\ob{i}_{\ell r}, \ob{i}_{r \ell})\monob \Cc(\ob{i}_{r}, \ob{i}_{r r}))
 $$
 and part of the data of $\delta^*$ consistes of, for each $\ob{i}\in \Cc$, a $1$-morphism 
 $$
 \delta^*_{\ob{i}}\in (\Cc\bxob(\Cc\bxob\Cc))( (\ob{i}_{\ell},(\ob{i}_{r\ell},\ob{i}_{rr})) , (\ob{i}_{\ell\ell},(\ob{i}_{\ell r},\ob{i}_{r}))  ) = 
 \Cc(\ob{i}_{\ell},\ob{i}_{\ell\ell} )\monob (\Cc(\ob{i}_{r \ell}, \ob{i}_{\ell r})\monob \Cc(\ob{i}_{r r}, \ob{i}_{r})).
 $$
 Moreover, the $2$-cell data of $\Theta$ gives  (coherent) isomorphisms
\begin{equation}\label{cohisos}\begin{split}
\theta^\ell_{\ob{j}}\circ \left(\luntwo_\Cc(\Theta\bxone  \mathbbm{1}_{\Cc})\Delta (\F) \right) & \xleftarrow{\theta^\ell_\F}  \F \circ  \theta^\ell_{\ob{i}} \qquad \text{ in } \Cc(\ob{i}_r,\ob{j})\\
 \theta^r_{\ob{j}}\circ \left(\runtwo_\Cc(  \mathbbm{1}_{\Cc}\bxone\Theta)\Delta (\F) \right) & \xleftarrow{\theta^r_\F}   \F \circ  \theta^r_{\ob{i}} \qquad \text{ in } \Cc(\ob{i}_\ell,\ob{j})\\
 \left(\luntwo_\Cc(\Theta\bxone  \mathbbm{1}_{\Cc})\Delta (\F) \right) \circ  \theta^{\ell *}_{\ob{i}} & \xrightarrow{\theta^{\ell *}_\F}  \theta^{\ell *}_{\ob{j}} \circ \F\qquad \text{ in } \Cc(\ob{i},\ob{j}_r)\\
 \left(\runtwo_\Cc(  \mathbbm{1}_{\Cc}\bxone\Theta)\Delta (\F)\right) \circ  \theta^{r *}_{\ob{i}} & \xrightarrow{\theta^{r *}_\F}  \theta^{r *}_{\ob{j}} \circ \F\qquad \text{ in } \Cc(\ob{i},\ob{j}_\ell)\\
  \end{split}\end{equation}
 for any $\F\in \Cc(\ob{i},\ob{j})$.
Similarly, for any $\F\in \Cc(\ob{i},\ob{j})$, we obtain coherent isomorphisms
\begin{equation}\label{cohisodelta}
\delta_{\ob{j}}\circ \left(\asstwo_{\Cc,\Cc,\Cc}(\Delta\bxone  \mathbbm{1}_{\Cc})\Delta(\F)\right)  \xleftarrow{\delta_\F}  \left(( \mathbbm{1}_{\Cc}\bxone \Delta)\Delta(\F)\right) \circ \delta_{\ob{i}}.
\end{equation}

By  \cite[Example 6.2.7]{JY}
\begin{itemize}
\item for any $\ob{i}$ in $\Cc$, $(\theta^\ell_{\ob{i}}, \theta^{\ell*}_{\ob{i}}, \eta^{\theta^\ell}_{\ob{i}}, \varepsilon^{\theta^\ell}_{\ob{i}})$, 
$(\theta^r_{\ob{i}}, \theta^{r*}_{\ob{i}}, \eta^{\theta^r}_{\ob{i}}, \varepsilon^{\theta^r}_{\ob{i}})$and
$(\delta_{\ob{i}}, \delta^*_{\ob{i}}, \eta^\delta_{\ob{i}}, \varepsilon^\delta_{\ob{i}})$ are internal adjoint equivalences in $\Cc$, respectively $\Cc\bxob(\Cc\bxob\Cc)$ (in the sense of \cite[Definition 6.1.1]{JY}) and
\item for any $\ob{i}, \ob{j}$ in $\Cc$, there are equivalences of categories 
\begin{align*}
\Cc(\ob{i}_r, \ob{j}_r) &\to \Cc(\ob{i}, \ob{j}), & \F &\mapsto\theta^\ell_{\ob{j}}\F\theta^{\ell*}_{\ob{i}}\\
\Cc(\ob{i}_\ell, \ob{j}_\ell) &\to \Cc(\ob{i}, \ob{j}), & \F &\mapsto\theta^r_{\ob{j}}\F\theta^{r*}_{\ob{i}}\\
\end{align*}
as well as an equivalence of categories
\begin{align*}
(\Cc\bxob(\Cc\bxob\Cc))((\ob{i}_{\ell\ell},(\ob{i}_{\ell r},\ob{i}_{r})),(\ob{j}_{\ell\ell},(\ob{j}_{\ell r},\ob{j}_{r}))) &\to (\Cc\bxob(\Cc\bxob\Cc))(  (\ob{i}_{\ell},(\ob{i}_{r\ell},\ob{i}_{rr})) , (\ob{j}_{\ell},(\ob{j}_{r\ell},\ob{j}_{rr}))  )\\
\F &\mapsto \delta_{\ob{j}}\F\delta^*_{\ob{i}}.
\end{align*}
\end{itemize}

 In particular, the units and counits of the adjunctions provide isomorphisms
 \begin{equation}\label{thetaeqs}
 \theta^\ell_{\ob{i}}\theta^{\ell*}_{\ob{i}} \xrightarrow{\varepsilon^{\theta^\ell}_{\ob{i}}} \mathbbm{1}_{\ob{i}} \qquad \theta^{\ell*}_{\ob{i}}\theta^\ell_{\ob{i}} \xleftarrow{\eta^{\theta^\ell}_{\ob{i}}} \mathbbm{1}_{\ob{i}_r}\qquad \theta^r_{\ob{i}}\theta^{r*}_{\ob{i}}\xrightarrow{\varepsilon^{\theta^r}_{\ob{i}}} \mathbbm{1}_{\ob{i}} \qquad \theta^{r*}_{\ob{i}}\theta^r_{\ob{i}}\xleftarrow{\eta^{\theta^r}_{\ob{i}}} \mathbbm{1}_{\ob{i}_\ell}
 \end{equation} 
 \begin{equation}\label{deltaeqs}
 \delta_{\ob{i}}\delta^*_{\ob{i}}  \xrightarrow{\varepsilon^{\delta}_{\ob{i}}} \mathbbm{1}_{(\ob{i}_{\ell},(\ob{i}_{r\ell},\ob{i}_{rr}))}   \qquad \delta^*_{\ob{i}}\delta_{\ob{i}}  \xleftarrow{\eta^{\delta}_{\ob{i}}} \mathbbm{1}_{ (\ob{i}_{\ell\ell},(\ob{i}_{\ell r},\ob{i}_{r}))} , 
 \end{equation} 
 
 \subsubsection{}\label{generaladj}
 Now consider the situation where 
 \begin{itemize}
 \item $\Dd$ is another small locally additive $\Bbbk$-linear $2$-category,
\item $\Gamma, \Gamma'$ are strictly unital, $\Bbbk$-linear pseudofunctors from $\Cc$ to $\Dd$,
\item $\beta\colon \Gamma \to \Gamma'$ and $\beta^*\colon \Gamma'\to \Gamma$ are strong natural transformations,
\item $\eta^\beta\colon \mathbbm{1}_\Gamma \to \beta^*\beta$ and $\varepsilon^\beta\colon \beta \beta^*\to  \mathbbm{1}_\Gamma'$ are modifications
\end{itemize}
such that $(\beta,\beta^*,\eta^\beta, \varepsilon^\beta)$ form an adjoint equivalence between $\Gamma$ and $\Gamma'$.

First suppose for simplicity that $ \mathbbm{1}_{\Gamma(\ob{i})} $ is indecomposable. We claim that this implies that  $\beta_{\ob{i}}$ and $\beta^*_{\ob{i}}$ are indecomposable. To see this, suppose $\beta_{\ob{i}} = \bigoplus_s \beta_{\ob{i},s}$ and $\beta^*_{\ob{i}} = \bigoplus_t \beta^*_{\ob{i},t}$ are decompositions into indecomposables.
Then the isomorphisms $\beta_{\ob{i}}\beta^*_{\ob{i}} \xrightarrow{\varepsilon^{\beta}_{\ob{i}}} \mathbbm{1}_{\Gamma'(\ob{i})}$ and  $ \beta^*_{\ob{i}}\beta_{\ob{i}}  \xleftarrow{\eta^{\beta}_{\ob{i}}} \mathbbm{1}_{\Gamma(\ob{i})} $ imply that there are unique pairs $(s,t)$ and $(t', s')$ such that 

$$\beta_{\ob{i},p}\beta^*_{\ob{i},q} \cong \begin{cases} 0 & \text{ if }(p,q)\neq (s,t) \\ \mathbbm{1}_{\Gamma'(\ob{i})}& \text{ if }(p,q) = (s,t) \end{cases} \quad \text{and}\quad \beta^*_{\ob{i},q}\beta_{\ob{i},p} \cong \begin{cases} 0 & \text{ if }(q,p)\neq (t',s') \\ \mathbbm{1}_{\Gamma(\ob{i})}& \text{ if }(q,p)= (t',s') \end{cases}.$$

Moreover, considering

$$\beta_{\ob{i}} \xrightarrow{\mathbbm{1}_{\beta_{\ob{i}}}\hcirc \eta^{\beta}_{\ob{i}}}  \beta_{\ob{i}}\beta^*_{\ob{i}}\beta_{\ob{i}} \xrightarrow{\varepsilon^{\beta}_{\ob{i}}\hcirc\mathbbm{1}_{\beta_{\ob{i}}} }\beta_{\ob{i}}, $$

the triangle identity implies that the composition is the identity on each $\beta_{\ob{i},p}$. However, the first morphism maps the direct summand $\beta_{\ob{i},p}$ of $\beta_{\ob{i}}$  to $\beta_{\ob{i},p}\beta^*_{\ob{i},t'}\beta_{\ob{i},s'}=0$ for any $p\neq s$. This implies that $\beta_{\ob{i}}$ and (by considering the other triangle identity) $\beta^*_{\ob{i}}$ are indecomposable, as claimed.

Now assume that $\mathbbm{1}_{\Gamma(\ob{i})} = \bigoplus_s  \mathbbm{1}_{\Gamma(\ob{i}),s}$ is a decomposition of $\mathbbm{1}_{\Gamma(\ob{i})} =\Gamma(\mathbbm{1}_{\ob{i}})$ into indecomposable $1$-morphisms. Given the equivalence of categories $\Dd(\Gamma(\ob{i}), \Gamma(\ob{j})) \to \Dd(\Gamma'(\ob{i}), \Gamma'(\ob{j}))$ from \cite[Example 6.2.7]{JY}, this implies that $ \mathbbm{1}_{\Gamma'(\ob{i})}$ has a decomposition into indecomposables given by $\mathbbm{1}_{\Gamma'(\ob{i})} = \bigoplus_s  \mathbbm{1}_{\Gamma'(\ob{i}),s}$ with the same indexing set, such that the restriction of $\varepsilon^{\beta}_{\ob{i}}$ yields an isomorphism $\beta_{\ob{i}}\mathbbm{1}_{\Gamma(\ob{i}),s}\beta^*_{\ob{i}}\xrightarrow{\sim} \mathbbm{1}_{\Gamma'(\ob{i}),s}$.

Setting $\beta_{\ob{i},s}=\beta_{\ob{i}}\mathbbm{1}_{\Gamma(\ob{i}),s}$ as a direct summand of $\beta_{\ob{i}}$, restriction of the triangle identity to this direct summand yields that the composition
$$\beta_{\ob{i},s}  \xrightarrow{\mathbbm{1}_{\beta_{\ob{i}, s}}\hcirc \eta^{\beta}_{\ob{i}}}   \beta_{\ob{i},s}\beta^*_{\ob{i}}\beta_{\ob{i}}  \xrightarrow{\varepsilon^{\beta}_{\ob{i}}\hcirc\mathbbm{1}_{\beta_{\ob{i}}} }\mathbbm{1}_{\Gamma'(\ob{i}),s}\beta_{\ob{i}} $$
is the identity, so $\beta_{\ob{i}}\mathbbm{1}_{\Gamma(\ob{i}),s}=\beta_{\ob{i},s} = \mathbbm{1}_{\Gamma'(\ob{i}),s}\beta_{\ob{i}}  =\mathbbm{1}_{\Gamma'(\ob{i}),s}\beta_{\ob{i}} \mathbbm{1}_{\Gamma(\ob{i}),s}$.

Similarly, setting $\beta_{\ob{i},s}^* =  \beta^*_{\ob{i}}\mathbbm{1}_{\Gamma'(\ob{i}),s}$, we obtain $\beta^*_{\ob{i}}\mathbbm{1}_{\Gamma'(\ob{i}),s}=\beta^*_{\ob{i},s} = \mathbbm{1}_{\Gamma(\ob{i}),s}\beta^*_{\ob{i}}  =\mathbbm{1}_{\Gamma(\ob{i}),s}\beta^*_{\ob{i}} \mathbbm{1}_{\Gamma'(\ob{i}),s}$. 
Moreover, the adjunction $(\beta_{\ob{i}} , \beta^*_{\ob{i}} )$ restricts to adjunctions $(\beta_{\ob{i},s} , \beta^*_{\ob{i},s} )$ and, given that by assumption $\mathbbm{1}_{\Gamma(\ob{i}),s}$ and $\mathbbm{1}_{\Gamma'(\ob{i}),s}$ are indecomposble, the above claim shows that $\beta_{\ob{i}} = \bigoplus_s\beta_{\ob{i},s}$ and $\beta^*_{\ob{i}} = \bigoplus_s \beta^*_{\ob{i},s}$ are decompositions into indecomposables.

\subsubsection{}
We now apply this to our adjoint equivalences $(\theta^\ell,\theta^{\ell *},\eta^{\theta^\ell}, \varepsilon^{\theta^\ell})$, $(\theta^r,\theta^{r *},\eta^{\theta^r}, \varepsilon^{\theta^r})$ and $(\delta,\delta^*,\eta^\delta, \varepsilon^\delta)$ in turn.

Suppose $\mathbbm{1}_{\ob{i}} = \bigoplus_{s=1}^n  \mathbbm{1}_{\ob{i},s}$ is a decomposition into indecomposable $1$-morphisms. Considering $\theta^\ell$ and $\theta^r$, we see that this implies decompositions $$\mathbbm{1}_{\ob{i}_r} = \bigoplus_{s=1}^n  \mathbbm{1}_{\ob{i}_r,s} \quad \mathbbm{1}_{\ob{i}_\ell} = \bigoplus_{s=1}^n  \mathbbm{1}_{\ob{i}_\ell,s}$$ $$\theta^\ell_{\ob{i}} = \bigoplus_{s=1}^n  \theta^\ell_{\ob{i},s} \quad \theta^{\ell *}_{\ob{i}} = \bigoplus_{s=1}^n  \theta^{\ell *}_{\ob{i},s} \quad \theta^r_{\ob{i}} = \bigoplus_{s=1}^n  \theta^r_{\ob{i},s} \quad \theta^{r *}_{\ob{i}} = \bigoplus_{s=1}^n  \theta^{r *}_{\ob{i},s} $$
where each summand is indecomposable and the adjunction restricts to isomorphisms $ \theta^\ell_{\ob{i},s} \theta^{\ell *}_{\ob{i},s} \to  \mathbbm{1}_{\ob{i},s}$, etc.

Now consider the adjoint equivalence $(\delta,\delta^*,\eta^\delta, \varepsilon^\delta)$ and set $\Gamma = \asstwo_{\Cc,\Cc,\Cc}(\Delta\bxone  \mathbbm{1}_{\Cc})\Delta$ and $\Gamma' =  ( \mathbbm{1}_{\Cc}\bxone \Delta)\Delta$. For $\Gamma(\ob{i}) = (\ob{i}_{\ell\ell},(\ob{i}_{\ell r},\ob{i}_{r}))$ and $\Gamma'(\ob{i}) =  (\ob{i}_{\ell},(\ob{i}_{r\ell},\ob{i}_{rr}))$, the previous paragraph implies that 
$$\mathbbm{1}_{\Gamma(\ob{i})} =  \bigoplus_{s_1=1}^n \bigoplus_{s_2=1}^n \bigoplus_{s_3=1}^n (\mathbbm{1}_{\ob{i}_{\ell\ell}, s_1},(\mathbbm{1}_{\ob{i}_{\ell r}, s_2},\mathbbm{1}_{\ob{i}_{r}, s_3})) \quad \mathbbm{1}_{\Gamma'(\ob{i})} =  \bigoplus_{s_1=1}^n \bigoplus_{s_2=1}^n \bigoplus_{s_3=1}^n (\mathbbm{1}_{\ob{i}_{\ell}, s_1},(\mathbbm{1}_{\ob{i}_{r \ell }, s_2},\mathbbm{1}_{\ob{i}_{rr}, s_3}))$$
such that for each $j\in \{1,2, 3\}$ the relevant combination of $\theta$ matches the corresponding summands indexed by $s_j\in \{1,\ldots, n\}$ to those indexed by the same number, e.g. $\mathbbm{1}_{\ob{i}_{\ell r}, s_2} $ is the image of
 $ \theta^{\ell}_{\ob{i}_\ell, s_2}  \theta^{r}_{\ob{i},s_2}  \mathbbm{1}_{\ob{i}, s_2}\theta^{r *}_{\ob{i}, s_2}  \theta^{\ell *}_{\ob{i}_\ell, s_2} $ under $\varepsilon^{\theta^\ell}_{\ob{i}_\ell}( \mathbbm{1}_{ \theta^{\ell}_{\ob{i}_\ell, s_2} }\hcirc\varepsilon^{\theta^r}_{\ob{i}}\hcirc \mathbbm{1}_{ \theta^{\ell *}_{\ob{i}_\ell, s_2} } )$, etc..
 In particular, each summand in the above compositions is indecomposable.
 
 Now consider $\delta_{\ob{i}}$ and  $\delta^*_{\ob{i}}$. Section \ref{generaladj} implies that there exist decompositions into indecomposables given by 
 
 $$\delta_{\ob{i}} =  \bigoplus_{s_1=1}^n \bigoplus_{s_2=1}^n \bigoplus_{s_3=1}^n (\delta_{\ob{i},1, s_1},(\delta_{\ob{i},2, s_2},\delta_{\ob{i},3, s_3})) \quad 
 \delta^*_{\ob{i}} =  \bigoplus_{s_1=1}^n \bigoplus_{s_2=1}^n \bigoplus_{s_3=1}^n (\delta^*_{\ob{i},1, s_1},(\delta^*_{\ob{i}, 2,s_2},\delta^*_{\ob{i},3, s_3}))$$
 such that the adjunction morphisms $ \eta^{\delta}_{\ob{i}}, \varepsilon^{\delta}_{\ob{i}}$ restrict to isomorphisms
$$(\mathbbm{1}_{\ob{i}_{\ell\ell}, s_1},(\mathbbm{1}_{\ob{i}_{\ell r}, s_2},\mathbbm{1}_{\ob{i}_{r}, s_3})) \xrightarrow{\sim}(\delta^*_{\ob{i},1, s_1}\delta_{\ob{i},1, s_1},(\delta^*_{\ob{i},2, s_2}\delta_{\ob{i}, 2,s_2},\delta^*_{\ob{i} ,3,s_3}\delta_{\ob{i}, 3, s_3}))$$
and
$$ (\delta_{\ob{i}, 1,s_1}\delta^*_{\ob{i},1, s_1},(\delta_{\ob{i}, 2,s_2}\delta^*_{\ob{i},2, s_2},\delta_{\ob{i},3, s_3}\delta^*_{\ob{i},3, s_3}))\xrightarrow{\sim} (\mathbbm{1}_{\ob{i}_{\ell}, s_1},(\mathbbm{1}_{\ob{i}_{r \ell }, s_2},\mathbbm{1}_{\ob{i}_{rr}, s_3})),$$
respectively.
In particular, we obtain component isomorphisms 
$$\mathbbm{1}_{\ob{i}_{\ell\ell}, s_1} \to\delta^*_{\ob{i}, 1,s_1}\delta_{\ob{i}, 1,s_1}\quad  \mathbbm{1}_{\ob{i}_{\ell r}, s_2} \to\delta^*_{\ob{i},2, s_2}\delta_{\ob{i}, 2,s_2}  \quad  \mathbbm{1}_{\ob{i}_{r}, s_3} \to\delta^*_{\ob{i}, 3, s_3}\delta_{\ob{i}, 3, s_3},$$
$$\delta_{\ob{i}, 1,s_1}\delta^*_{\ob{i}, 1,s_1}\to \mathbbm{1}_{\ob{i}_{\ell}, s_1}\quad \delta_{\ob{i}, 2,s_2}\delta^*_{\ob{i}, 2, s_2}\to \mathbbm{1}_{\ob{i}_{r\ell}, s_2}   \quad  \delta_{\ob{i}, 3, s_3}\delta^*_{\ob{i}, 3, s_3}\to \mathbbm{1}_{\ob{i}_{rr}, s_3}. $$

Note further that we can write $\delta_{\ob{i}}= (\delta_{\ob{i},1},(\delta_{\ob{i},2}, \delta_{\ob{i},3}))$ where we set $\delta_{\ob{i},j} = \bigoplus_{s_j=1}^n\delta_{\ob{i},j, s_j}$ for $j\in  \{1,2, 3\}$.

Since  $\Theta(\mathbbm{1}_{\ob{i}}) =\mathbbm{1}_{\bullet}$ for all $\ob{i}\in \Cc$, the isomorphisms  in \eqref{thetaeqs} also imply that $\Theta(\F) = \mathbbm{1}_{\bullet}$ for any $\F\in \{\theta^\ell_{\ob{i}}, \theta^{\ell*}_{\ob{i}}, \theta^r_{\ob{i}}, \theta^{r*}_{\ob{i}}\}$. Moreover, since we can write $\delta_{\ob{i}}= (\delta_{\ob{i},1},(\delta_{\ob{i},2}, \delta_{\ob{i},3}))$ and similarly for $\delta_{\ob{i}}^*$, we see that $\delta_{\ob{i},j}\delta_{\ob{i},j}^*$ is isomorphic to the identity on the corresponding object for $j\in  \{1,2, 3\}$, hence we record that 

\begin{equation}\label{Thetatoone}
\Theta(\F) = \mathbbm{1}_{\bullet} \quad \text{ for }\quad \F\in \{\theta^\ell_{\ob{i}}, \theta^{\ell*}_{\ob{i}}, \theta^r_{\ob{i}}, \theta^{r*}_{\ob{i}},\delta_{\ob{i},j},\delta_{\ob{i},j}^*\,\vert\, j=1,2,3\}.
\end{equation}

\section{Monoidal structure on the $2$-category of birepresentations}\label{sec:mon2catof2reps}

{ 
We denote by $\Cc\amod$ the $2$-category of strictly unital birepresentations of a $2$-category $\Cc\in \twoAk$. Explicitly, objects are strictly unital birepresentations, i.e.\ strictly unital pseudofunctors, from $\Cc$ to $\A$, $1$-cells are morphisms of strictly unital birepresentations, i.e.\ strong natural transformations, and $2$-morphisms are modifications.
}

\subsection{External tensor products of birepresentations}

%
%
%
%
%
%

Let  $\Cc$ and $\Dd$  be $2$-categories in $\twoAk$ and let  $\bfM\colon \Cc\to \A$ and $ \bfN\colon \Dd\to \A$ be strictly unital birepresentations. Given these are $\Bbbk$-linear strictly unital pseudofunctors and $\A$ is an object of $\twoAk$, they are, in particular, $1$-cells in $\twoAk$.
Thus we can consider their tensor product $\bfM\bxone\bfN \colon \Cc\bxob\Dd \to \A\bxob \A$. 
We would like to obtain a strictly unital birepresentation of $ \Cc\bxob\Dd$, i.e. a strictly unital pseudofunctor from $ \Cc\bxob\Dd$ to $\A$. To see that this is possible, 
note that we have a canonical morphism $\desc\colon \A\times\A \to \A\bxob \A$. Indeed, recall that
\begin{itemize}
\item objects in both $2$-categories are pairs $(\oneC, \oneD)$ of objects $\oneC,\oneD$ in $\A$;
\item morphism categories in $\A\times \A$ are $\A(\oneC, \oneC')\times \A(\oneD, \oneD')$ while those in $\A\bxob \A$ are given by  $\A(\oneC, \oneC')\monob \A(\oneD, \oneD')$
\end{itemize}
so we can define $\desc$ as the identity on objects and the canonical functor given by the tensor product of $1$-categories on morphism categories. 
 We then need the following lemma.

\begin{lemma}
The $2$-functor $\A\times \A\xrightarrow{\otimes} \A$ (Prop.~\ref{tensorone}) factors over $\A\bxob \A$, i.e. there is a $2$-functor $\what\colon \A\bxob \A \to \A$ making the diagram
$$\xymatrix{ 
\A\times \A \ar_{\desc}[dr] \ar^{\otimes}[rr]&&\A\\
& \A\bxob \A \ar^{\what}[ur]&
}$$
commute.
\end{lemma}

\proof
Recalling the description of objects and morphism categories in $\A\bxob \A$ above, we define $\what (\oneC, \oneD)= \oneC\monob \oneD$ on the level of objects. We now define the functor $$\what_{(\oneD,\oneD'),(\oneC,\oneC')}\colon \A(\oneC, \oneC')\monob \A(\oneD, \oneD')\to \A(\oneC\monob \oneC',\oneD\monob\oneD')$$ as follows. Observe that there is a functor 
\begin{align*}
	\what'_{(\oneD,\oneD'),(\oneC,\oneC')}\colon \A(\oneC, \oneC')\diamond \A(\oneD, \oneD')&\to \A(\oneC\monob \oneC',\oneD\monob\oneD'),\\
(F,G)&\mapsto F\monmor G,\\
\alpha \otimes \beta &\mapsto \alpha \monnat \beta.\end{align*}
Now define $\what_{(\oneD,\oneD'),(\oneC,\oneC')}=\Add{(\what'_{(\oneD,\oneD'),(\oneC,\oneC')})}$.

It is straightforward to check that $\what$ is a $2$-functor, which makes the diagram in the lemma commute.


\endproof

We define $\bfM\tenrep_0\bfN = \what\circ (\bfM\bxone\bfN)\colon \Cc\bxob\Dd \to \A $. In particular, 
\begin{itemize}
\item $(\bfM\tenrep_0\bfN)(\ob{i}, \ob{k})= \bfM(\ob{i})\monob\bfN(\ob{k})$ for objects $\ob{i}\in\Cc, \ob{k}\in\Dd$, 
\item $(\bfM\tenrep_0\bfN)(\F, \G) = \bfM(\F)\monmor\bfN(\G)$  for $1$-morphisms $\F$ and $\G$ in $\Cc$ and $\Dd$, respectively, and
\item $(\bfM\tenrep_0\bfN)(\alpha\otimes \beta) = \bfM(\alpha)\monnat \bfN(\beta)$ for $2$-morphisms $\alpha$ and $\beta$ in $\Cc$ and $\Dd$, respectively.
\end{itemize}

For $1$-cells $\Phi\colon \bfM\to\bfM'$ and $\Psi\colon \bfN\to \bfN'$ in $\Cc\amod$  and $\Dd\amod$, respectively, we define their external tensor product $\Phi\tenrep_1\Psi\colon \bfM\tenrep_0\bfN\to \bfM'\tenrep_0\bfN'$ by $\Phi\tenrep_1\Psi = \what\star(\Phi\bxtwo\Psi)$, noting that $\Phi$ and $\Psi$ are $2$-morphisms in $\twoAk$. Explicitly,
\begin{itemize}
\item $(\Phi\tenrep_1\Psi)_{(\ob{i}, \ob{k})}= \Phi_{\ob{i}}\monmor\Psi_{\ob{k}}$ for objects $\ob{i}\in\Cc, \ob{k}\in\Dd$, and
\item $(\Phi\tenrep_1\Psi)_{(\ob{j}, \ob{l})(\ob{i}, \ob{k})}$ has components  
\begin{equation*}\begin{split}
((\Phi\tenrep_1\Psi)_{(\ob{j}, \ob{l})(\ob{i}, \ob{k})})_{\F,\G}&= (\what \star (\Phi\bxtwo\Psi)_{(\ob{j}, \ob{l})(\ob{i}, \ob{k})})_{\F,\G}\\
&=\what((\Phi_{\ob{j}\ob{i}})_\F\otimes(\Psi_{\ob{l}\ob{k}})_\G)\\
&=(\Phi_{\ob{j}\ob{i}})_\F\monnat(\Psi_{\ob{l}\ob{k}})_\G
\end{split}\end{equation*}
for objects $\ob{i},\ob{j}\in\Cc, \ob{k},\ob{l}\in\Dd$ and $1$-morphisms $\F\in\Cc(\ob{i},\ob{j}), \G\in \Dd(\ob{k},\ob{l})$, while on general direct sums of $1$-morphisms, we obtain diagonal matrices as in \eqref{otimesdef}.
\end{itemize}


Now let $\Phi, \Phi\colon \bfM\to\bfM'$ and $\Psi,\Psi'\colon \bfN\to \bfN'$ be $1$-cells in $\Cc\amod$  and $\Dd\amod$, respectively, and let $\theta\colon \Phi\to\Phi'$ and $\tau\colon \Psi\to\Psi'$ be modifications. These are given by the data of
a natural transformation $\theta_{\ob{i}}\colon  \Phi_{\ob{i}}\to \Phi'_{\ob{i}}$ and $\tau_{\ob{k}}\colon  \Psi_{\ob{k}}\to \Psi'_{\ob{k}}$, respectively, so we define $(\theta\tenrep_2\tau)_{(\ob{i}, \ob{k})} = \theta_{\ob{i}}\monnat\tau_{\ob{k}}$.

It is easy to check that with the notation from above, $\bfM\tenrep_0\bfN$ is a strictly unital birepresentation, $\Phi\tenrep_1\Psi$ is a morphism of such birepresentations, and $\theta\tenrep_2\tau$ is a modification of such morphisms.

\subsection{Tensor product of strictly unital birepresentations}

Assume that $\Cc$ is equipped with a comonoidal structure. 

Given $\bfM,\bfN\in\Cc\amod$, we can consider $\bfM \tenrepC_0\bfN$ in $\Cc\amod$ defined by  $\bfM \tenrepC_0\bfN =  (\bfM \tenrep\bfN)\Delta$. Similarly, given $1$-cells $\Phi\colon \bfM\to\bfM'$ and $\Psi\colon \bfN\to \bfN'$ in $\Cc\amod$, we obtain a $1$-cell 
$\Phi\tenrepC_1\Psi\colon \bfM\tenrepC_0\bfN\to \bfM'\tenrepC_0\bfN'$ by setting
$\Phi\tenrepC_1\Psi = (\Phi\tenrep_1\Psi)\star\Delta$. Explicitly,  
\begin{align*}
(\Phi\tenrepC_1\Psi)_{\ob{i}} & = \Phi_{\ob{i}_\ell}\monmor\Psi_{\ob{i}_r}\\
(\Phi\tenrepC_1\Psi)_{\ob{j}\ob{i}} &=\Phi_{\ob{j}_\ell\ob{i}_\ell}\monnat\Psi_{\ob{j}_r\ob{i}_r}
\end{align*}
For modifications $\theta\colon \Phi\to\Phi'$ and $\tau\colon \Psi\to\Psi'$ of $1$-cells $\Phi,\Phi', \Psi,\Psi'$ in $\Cc\amod$, we define a modification $\theta\tenrepC_2\tau = (\theta\tenrep_2\tau)\mathbbm{1}_{\mathbbm{1}_{\Delta}}$. Explicitly,  $(\theta\tenrepC_2\tau)_{\ob{i}} = \theta_{\ob{i}_{\ell}} \monnat \tau_{\ob{i}_r} $.

\begin{prop}\label{tensorthree}
The functions $\tenrepC_0,\tenrepC_1,\tenrepC_2$ on the levels of objects, $1$-cells and $2$-cells define a $2$-functor $\tenrepC \colon \Cc\amod\times\Cc\amod\to\Cc\amod$.
\end{prop}

\proof
Functoriality on morphism categories as well as strict unitality follow directly from the definition and the same properties for $\otimes\colon \A \times \A \to \A$.

To see that for morphisms $\Phi\colon \bfM\to\bfM',\Phi'\colon \bfM'\to\bfM'',\Psi\colon \bfN\to\bfN',\Psi'\colon \bfN'\to\bfN''$, we have 
$$ \Phi'\Phi \tenrepC_1 \Psi'\Psi =  (\Phi' \tenrepC_1 \Psi') (\Phi \tenrepC_1 \Psi)$$
consider the definition and notice that
\begin{equation*}\begin{split}
((\Phi' \tenrep_1 \Psi')&\star\Delta) ((\Phi \tenrep_1 \Psi)\star\Delta)\\
&=((\Phi' \tenrep_1 \Psi') (\Phi \tenrep_1 \Psi))\star\Delta
\end{split}\end{equation*}
by Lemma~\ref{le:partialinterchange}, so it suffices to check that $\Phi'\Phi \tenrep_1 \Psi'\Psi =  (\Phi' \tenrep_1 \Psi') (\Phi \tenrep_1 \Psi)$. Given the definition of $\tenrep_1$, this immediately follows from $2$-functoriality of $\boxtimes$.
%
\endproof

\subsection{Unitors for tensor products of strictly unital birepresentations}

Notice that the $\Bbbk$-linear $2$-category $\Vv$ has a simple strictly unital birepresentation 
\begin{align*}
\bfT_\Vv \colon \Vv\to \A^f \qquad \bullet &\mapsto\Add\Bbbk\\
\mathbbm{1}_\bullet^{\oplus n} = \bullet_\Bbbk^{\oplus n} &\mapsto \mathbbm{1}_{\Add \Bbbk}^{\oplus n}\\
f&\mapsto f,
\end{align*}
where, by slight abuse of notation, we again denote by $f$ the natural transformation canonically induced by $f$.

We define the strictly unital birepresentation $\bfT_\Cc$ of $\Cc$ by $\bfT_\Cc=\bfT_\Vv\Theta$. Then, as in Proposition \ref{unitone}, the following proposition is trivially true.

\begin{prop}\label{unitthree}
There is a $2$-functor 
\begin{align*}\One_{\Cc\amod}\colon\triv \to \Cc\amod \colon\qquad \trivob &\mapsto \bfT_\Cc\\
 \mathbbm{1}_\trivob &\mapsto \mathbbm{1}_{\bfT_\Cc}\\
\mathbbm{1}_{\mathbbm{1}_\trivob} &\mapsto \mathbbm{1}_{\mathbbm{1}_{\bfT_\Cc}}.
\end{align*}
\end{prop}

To define left and right unitors for $\tenrepC$, we first need some intermediate constructions. The data involved in defining the left and right unitors in Definition \ref{unitorthree} will locally consist of certain morphisms of strictly unital birepresentations, to be defined via pasting diagrams \eqref{upasting1} and \eqref{upasting2} in Proposition \ref{unitorthreeM}. These pasting diagrams involve strong natural transformations 
$$\tau_\bfM^\ell \colon \bfT_\Vv\tenrep_0\bfM \rightleftarrows \bfM\luntwo_\Cc \colon \tau_\bfM^{\ell*} $$
$$\tau_\bfM^r \colon \bfM\tenrep_0\bfT_\Vv \rightleftarrows \bfM \runtwo_\Cc \colon \tau_\bfM^{r*} $$
 for any $\bfM \in \Cc\amod$, which we will now construct.

Recall from Proposition~\ref{unitorone} that $\lun_{\Cc(\ob{i},\ob{j})}\colon \Add \Bbbk \monob\Cc(\ob{i},\ob{j})\underset{}{ \overset{}{\rightleftarrows}} \Cc(\ob{i},\ob{j})\hspace{0.3em}\colon\hspace{-0.3em} \lun_{\Cc(\ob{i},\ob{j})}^*$ and $\run_{\Cc(\ob{i},\ob{j})}\colon  \Cc(\ob{i},\ob{j})\monob \Add \Bbbk\underset{}{ \overset{}{\rightleftarrows}} \Cc(\ob{i},\ob{j})\hspace{0.3em}\colon\hspace{-0.3em} \run_{\Cc(\ob{i},\ob{j})}^*$ form adjoint equivalences.

Let $\bfM\in \Cc\amod$. 
For each object $(\bullet,\ob{i})\in \Vv\bxob\Cc$, notice that 
$$(\bfT_\Vv\tenrep_0\bfM)(\bullet,\ob{i})=\what (\bfT_\Vv\bxone \bfM)(\bullet,\ob{i}) = \Add\Bbbk\monob \bfM(\ob{i})$$ and $\bfM \luntwo_\Cc (\bullet,\ob{i}) = \bfM(\ob{i})$, thus we can define the functors 
\begin{equation}\label{tauell1}
(\tau_\bfM^\ell )_{(\bullet,\ob{i})} = \lun_{\bfM(\ob{i})}\colon (\bfT_\Vv\tenrep_0\bfM)(\bullet,\ob{i})\to \bfM \luntwo_\Cc (\bullet,\ob{i})
\end{equation}
and 
\begin{equation}\label{tauell*1}(\tau_\bfM^{\ell*} )_{(\bullet,\ob{i})}= \lun_{\bfM(\ob{i})}^*\colon \bfM \luntwo_\Cc (\bullet,\ob{i})\to (\bfT_\Vv\tenrep_0\bfM)(\bullet,\ob{i}).
\end{equation}

Similarly, observing that $$(\bfM\tenrep_0\bfT_\Vv)(\ob{i},\bullet) =  \what (\bfM\bxone \bfT_\Vv)(\ob{i},\bullet) = \bfM(\ob{i})\monob \Add\Bbbk$$ and $\bfM \runtwo_\Cc (\bullet,\ob{i}) = \bfM(\ob{i})$, we define the functors
\begin{equation}\label{taur1}(\tau_\bfM^r )_{(\ob{i},\bullet)} = \run_{\bfM(\ob{i})} \colon(\bfM\tenrep_0\bfT_\Vv)(\ob{i},\bullet)\to \bfM \runtwo_\Cc (\ob{i},\bullet)
\end{equation}
and 
\begin{equation}\label{taur*1}(\tau_\bfM^{r*} )_{(\ob{i},\bullet)} = \run^*_{\bfM(\ob{i})} \colon\bfM \runtwo_\Cc (\ob{i},\bullet)\to (\bfM\tenrep_0\bfT_\Vv)(\ob{i},\bullet).
\end{equation}

These functors, which are $1$-cells in $\A$, will provide the $1$-cell components of our desired strong natural transformations $\tau_\bfM^\ell, \tau_\bfM^{\ell*}, \tau_\bfM^r$ and $\tau_\bfM^{r*}$, respectively. 

To define the $2$-cell components, consider a $1$-cell $\bigoplus_p(\mathbbm{1}_\bullet^{\oplus s_p}, \F_p)$ in $(\Vv\bxob \Cc)( (\bullet, \ob{i}), (\bullet, \ob{j}))$.

Notice that both $\lun_{\bfM(\ob{j})}((\bfT_\Vv\tenrep_0\bfM)(\bigoplus_p(\mathbbm{1}_\bullet^{\oplus s_p},\F_p)))$ and $\bfM( \luntwo_\Cc (\bigoplus_p(\mathbbm{1}_\bullet^{\oplus s_p},\F_p)))
\lun_{\bfM(\ob{i})}$ are functors $\Add \Bbbk \monob \bfM(\ob{i}) \to \bfM(\ob{j})$.
On an object $(\bullet_\Bbbk^{\oplus n}, X)$ in $\Add\Bbbk\monob \bfM(\ob{i})$,  
\begin{equation*}\begin{split}\left(\lun_{\bfM(\ob{j})}((\bfT_\Vv\tenrep_0\bfM)(\bigoplus_p(\mathbbm{1}_\bullet^{\oplus s_p},\F_p)))\right)(\bullet_\Bbbk^{\oplus n}, X) &= \lun_{\bfM(\ob{j})}(\bigoplus_p((\bullet_\Bbbk^{\oplus n})^{\oplus s_p}, \bfM(\F_p)X)) \\ 
&= \bigoplus_p((\bfM(\F_p)X)^{\oplus n})^{\oplus s_p}
\end{split}\end{equation*}
while
$$\left( \bfM(\luntwo_\Cc (\bigoplus_p(\mathbbm{1}_\bullet^{\oplus s_p},\F_p)))
\lun_{\bfM(\ob{i})}\right) (\bullet_\Bbbk^{\oplus n}, X) = \bfM (\bigoplus_p\F_p^{\oplus s_p}) (X^{\oplus n}).
$$

Thus we can define a natural transformation
\begin{equation}\label{tauell2}(\tau_\bfM^\ell )_{ (\bigoplus_p(\mathbbm{1}_\bullet^{\oplus s_p},\F_p))
} \colon 
	\bfM( \luntwo_\Cc (\bigoplus_p(\mathbbm{1}_\bullet^{\oplus s_p},\F_p)))\lun_{\bfM(\ob{i})}  \to   \lun_{\bfM(\ob{j})}((\bfT_\Vv\tenrep_0\bfM)(\bigoplus_p(\mathbbm{1}_\bullet^{\oplus s_p},\F_p))
)
\end{equation}
as having components given by the natural isomorphisms between these, induced by additivity of $\bfM$ and additivity of each $\bfM(\F_p)$.

Morover, for an object $X\in \bfM(\ob{i})$, 
$$((\bfT_\Vv\tenrep \bfM)(\bigoplus_p(\bullet^{\oplus s_p}_\Bbbk, \F_p)))\lun_{\bfM(\ob{i})}^*(X) = 
((\bfT_\Vv\tenrep \bfM)(\bigoplus_p(\bullet^{\oplus s_p}_\Bbbk, \F_p))) (\bullet_\Bbbk, X)
 = \bigoplus_p(\bullet_\Bbbk^{\oplus s_p}, \bfM(\F_p)X)$$
 while
$$\left(\lun_{\bfM(\ob{j})}^*   (\bfM \luntwo_\Cc (\bigoplus_p(\bullet^{\oplus s_p}_\Bbbk, \F_p)))\right)(X) = \lun_{\bfM(\ob{j})}^*(\bfM (\bigoplus_p\F_p^{\oplus s_p})(X)) = (\bullet_\Bbbk, \bigoplus_p(\bfM (\F_p^{\oplus s_p})(X)))$$
so defining 
\begin{equation}\label{tauell*2}(\tau_\bfM^{\ell*} )_{\bigoplus_p (\bullet^{\oplus s_p}_\Bbbk,\F_p)} \colon 
((\bfT_\Vv\tenrep \bfM)(\mathbbm{1}_\bullet^{\oplus s},\F)) \lun_{\bfM(\ob{i})}^*
\to
\lun_{\bfM(\ob{j})}^*   \bfM \luntwo_\Cc (\mathbbm{1}_\bullet^{\oplus s},\F)
\end{equation}
to have the obvious natural isomorphisms between these (given by additivity of $\bfM$ and resummation) as its components provides a natural transformation, hence a $2$-cell in $\A$.

%

Similarly, for a $1$-morphism $\bigoplus_p(\F_p,\mathbbm{1}_\bullet^{\oplus s_p})\in  (\Cc\bxob \Vv)( (\ob{i}, \bullet), (\ob{j}, \bullet))$
\begin{equation*}\begin{split}
\left(\run_{\bfM(\ob{j})}((\bfM\tenrep\bfT_\Vv )(\bigoplus_p(\F_p,\mathbbm{1}_\bullet^{\oplus s_p})))\right)( X,\bullet_\Bbbk^{\oplus n}) &= \run_{\bfM(\ob{j})}\bigoplus_p( \bfM(\F_p)X,(\bullet_\Bbbk^{\oplus n})^{\oplus s_p}) \\ &=\bigoplus_p ((\bfM(\F_p)X)^{\oplus n})^{\oplus s_p}
\end{split}\end{equation*}
while
\begin{equation*}\begin{split}
\left( \bfM(\runtwo_\Cc (\bigoplus_p(\F_p,\mathbbm{1}_\bullet^{\oplus s_p})))
\run_{\bfM(\ob{i})} \right)( X,\bullet_\Bbbk^{\oplus n}) &= \bfM (\bigoplus_p\F_p^{\oplus s_p}) (X^{\oplus n}) 
\end{split}\end{equation*}
so again we can define 
\begin{equation}\label{taur2}(\tau_\bfM^r )_{ (\F,\mathbbm{1}_\bullet^{\oplus s})} \colon 
\bfM(\runtwo_\Cc (\F,\mathbbm{1}_\bullet^{\oplus s}))\run_{\bfM(\ob{i})}    \to  
\run_{\bfM(\ob{j})}((\bfM\tenrep\bfT_\Vv )(\F,\mathbbm{1}_\bullet^{\oplus s}))
\end{equation}
as the natural transformation whose components are given by the natural isomorphisms between these, induced by additivity of $\bfM$ and each $\bfM(\F_p)$.



Finally, for an object $X\in \bfM(\ob{i})$, we have
$$( (\bfM\tenrep \bfT_\Vv)(\bigoplus_p(\F_p,\mathbbm{1}_\bullet^{\oplus s_p}))) \run_{\bfM(\ob{i})}^*(X) = 
((\bfM\tenrep \bfT_\Vv)(\bigoplus_p(\F_p,\mathbbm{1}_\bullet^{\oplus s_p}))) ( X,\bullet_\Bbbk)
 = \bigoplus_p(\bfM(\F_p)X,\bullet_\Bbbk^{\oplus s_p})$$ 
 while 
 $$\left(\run_{\bfM(\ob{j})}^*   \bfM \runtwo_\Cc (\bigoplus_p(\F_p,\mathbbm{1}_\bullet^{\oplus s_p}))\right)(X) =\run_{\bfM(\ob{j})}^* (\bfM(\bigoplus_p\F_p^{\oplus s_p})(X))=( \bfM(\bigoplus_p\F_p^{\oplus s_p})(X), \bullet_\Bbbk)
 $$ 
 so again defining 
\begin{equation}\label{taur*2}(\tau_\bfM^{r*} )_{ \bigoplus_i(\F_i,\mathbbm{1}_\bullet^{\oplus s_i})} \colon 
((\bfM\tenrep \bfT_\Vv)(\F,\mathbbm{1}_\bullet^{\oplus s})) \run_{\bfM(\ob{i})}^*
\to
\run_{\bfM(\ob{j})}^*   \bfM \runtwo_\Cc (\F,\mathbbm{1}_\bullet^{\oplus s})
\end{equation}
to have as its component on each $X$ the natural isomorphism (given by additivity of $\bfM$ and resummation) between these, yields a natural transformation.

It is easy to see that this data assembles to yield the following proposition.

\begin{lemma}
\begin{enumerate}[(a)]
\item For any $\bfM\in \Cc\amod$, there is an adjoint equivalence $(\tau_\bfM^\ell, \tau_\bfM^{\ell*} , \eta^{\tau_\bfM^\ell} , \varepsilon^{\tau_\bfM^\ell} )$, where 
\begin{itemize}
\item the strong natural transformations $$\tau_\bfM^\ell \colon \bfT_\Vv\tenrep \bfM \rightleftarrows \bfM \luntwo_\Cc \colon \tau_\bfM^{\ell*} $$ between strictly unital pseudofunctors $\Vv\bxob \Cc \to \A$ (strictly unital birepresentations of $\Cc\bxob \Vv$) have component $1$-cells in $\A$ given by 
$$(\tau_\bfM^\ell )_{(\bullet,\ob{i})} = \lun_{\bfM(\ob{i})} \qquad  \text{resp.\ } (\tau_\bfM^{\ell*} )_{(\bullet,\ob{i})} = \lun^{*}_{\bfM(\ob{i})}$$ and component $2$-cells are given in \eqref{tauell2} resp.\ \eqref{tauell*2};
\item the modifications $$\eta^{\tau_\bfM^\ell} \colon \mathbbm{1}_{\bfT_\Vv\tenrep \bfM} \to \tau_\bfM^{\ell*} \tau_\bfM^\ell \qquad
 \varepsilon^{\tau_\bfM^\ell}\colon \tau_\bfM^\ell \tau_\bfM^{\ell*} \to \mathbbm{1}_{\bfM \luntwo_\Cc}$$ providing the adjunction
have as component $2$-cells the natural transformations $\eta^\ell_{\bfM(\ob{i})}$ and $\varepsilon^\ell_{\bfM(\ob{i})}$ given in Proposition \ref{unitorone}.
\end{itemize}
\item For any $\bfM\in \Cc\amod$,  there is an adjoint equivalence  $(\tau_\bfM^r, \tau_\bfM^{r*} , \eta^{\tau_\bfM^r} , \varepsilon^{\tau_\bfM^r} )$, where
\begin{itemize}
\item the strong natural transformations 
$$\tau_\bfM^r \colon \bfM\tenrep\bfT_\Vv \rightleftarrows \bfM \runtwo_\Cc \colon \tau_\bfM^{r*}$$
between strictly unital pseudofunctors $\Cc\bxob \Vv \to \A$ (strictly unital birepresentations of $\Cc\bxob \Vv$) have component $1$-cells in $\A$ given by 
$$(\tau_\bfM^r )_{(\ob{i},\bullet)} = \run_{\bfM(\ob{i})} \qquad  \text{resp.\ } (\tau_\bfM^{r*} )_{(\ob{i},\bullet)} = \run^*_{\bfM(\ob{i})}$$
and component $2$-cells in $\A$ given by \eqref{taur2} resp.\ \eqref{taur*2};
\item the modifications 
$$\eta^{\tau_\bfM^r} \colon \mathbbm{1}_{\bfM\tenrep \bfT_\Vv} \to \tau_\bfM^{r*}\tau_\bfM^r \qquad
 \varepsilon^{\tau_\bfM^r}\colon \tau_\bfM^r \tau_\bfM^{r*} \to \mathbbm{1}_{\bfM \runtwo_\Cc}$$ 
 providing the adjunction
have as component $2$-cells the natural transformations $\eta^r_{\bfM(\ob{i})}$ and $\varepsilon^r_{\bfM(\ob{i})}$ given in Proposition \ref{unitorone}.
\end{itemize}
\end{enumerate}
\end{lemma}

In the following, we will often use pasting diagrams, where the $0$-cells are locally additive $\Bbbk$-linear $2$-categories, the $1$-cells are strictly unital pseudofunctors, and the $2$-cells are strong natural transformations. A priori, such pasting diagrams are not well-defined since the interchange law only holds up to invertible modification for strong natural transformations. However, in all instances where we depict a resulting strong natural transformation as a pasting diagram, there will either be a unique way of composing, or it will be clear from context, in which order the cells need to be composed.


\begin{defn}\label{unitorthreeM}
\begin{enumerate}[(a)]
\item\label{unitorthreeMleft} For any $\bfM\in \Cc\amod$, we define $1$-cells  $$\lunthree_{\bfM}\colon  \bfT_\Cc\tenrepC_0\bfM \to\bfM \quad \text{and}\quad \lunthree_{\bfM}^*\colon  \bfM \to\bfT_\Cc\tenrepC_0\bfM$$ in $\Cc\amod$ by the pasting diagram
\begin{equation}\label{upasting1}\xymatrix{
\Cc\ar@{=}[rrrrdd] 
\ar^{\Delta}[rr] &&\Cc\bxob\Cc \ar^{\Theta\bxone\mathbbm{1}_{\Cc}}[rr] && \Vv\bxob\Cc  \ar^{\luntwo_\Cc}[dd]\ar^{\bfT_{\Vv}\bxone \bfM}[rr]\ar@2{->}@<5pt>[dl]^{\theta^\ell}&& \A\bxob \A\ar^{\what}[ddd]
\ar@2{->}@<5pt>[ddl]^{\tau^\ell_\bfM}\\
&&&\ar@2{->}@<5pt>[ur]^{\theta^{\ell*}}&&&\\
&&&&\Cc\ar_{\bfM}[drr]&\ar@2{->}@<5pt>[uur]^{\tau^{\ell*}_\bfM}&\\
&&&&&&\A.
}\end{equation}

\item\label{unitorthreeMright} For any $\bfM\in \Cc\amod$, we define $1$-cells $$\runthree_{\bfM}\colon \bfM \tenrepC_0 \bfT_\Cc \to\bfM  \quad \text{and}\quad \runthree^*_{\bfM}\colon \bfM  \to\bfM \tenrepC_0 \bfT_\Cc $$ in $\Cc\amod$ by the pasting diagram
\begin{equation}\label{upasting2}\xymatrix{
\Cc\ar@{=}[rrrrdd] 
\ar^{\Delta}[rr] &&\Cc\bxob\Cc \ar^{\mathbbm{1}_{\Cc}\bxone\Theta}[rr] && \Cc\bxob\Vv  \ar^{\runtwo_\Cc}[dd]\ar^{\bfM\bxone \bfT_{\Vv}}[rr]\ar@2{->}@<5pt>[dl]^{\theta^r}&& \A\bxob \A\ar^{\what}[ddd]
\ar@2{->}@<5pt>[ddl]^{\tau^r_\bfM}\\
&&&\ar@2{->}@<5pt>[ur]^{\theta^{r*}}&&&\\
&&&&\Cc\ar_{\bfM}[drr]&\ar@2{->}@<5pt>[uur]^{\tau^{r*}_\bfM}&\\
&&&&&&\A.
}\end{equation}

\end{enumerate}
\end{defn}

By definition, this means that $\lunthree_\bfM$ is the composition of the post-whiskering of $\theta_\ell$ with $\bfM$ and the pre-whiskering of $\tau^\ell_\bfM$ with $(\Theta\bxone\mathbbm{1}_{\Cc})\Delta$, see \cite[Definition 11.1.1]{JY}, and $\lunthree^*_\bfM, \runthree_\bfM, \runthree^*_\bfM$ are defined similarly.

Explicitly, for $\ob{i}\in \Cc$, the $1$-cell component of the strong natural transformations $\lunthree_\bfM, \lunthree^*_\bfM, \runthree_\bfM, \runthree^*_\bfM$ are given by the following functors:
\begin{align*}
(\lunthree_{\bfM})_{\ob{i}} = \bfM(\theta^{\ell}_{\ob{i}}) \lun_{\bfM(\ob{i}_r)}&\colon \Add\Bbbk\monob \bfM(\ob{i}_r) \to \bfM(\ob{i}), \quad \bigoplus_s(\bullet_\Bbbk^{\oplus n_s},X_s)\mapsto \bfM(\theta^{\ell}_{\ob{i}}) (\bigoplus_s X_s^{\oplus n_s})\\
(\runthree_{\bfM})_{\ob{i}} = \bfM(\theta^{r}_{\ob{i}}) \run_{\bfM(\ob{i}_{\ell})}&\colon  \bfM(\ob{i}_\ell)\monob \Add\Bbbk\to \bfM(\ob{i}), \quad \bigoplus_s (X_s, \bullet_\Bbbk^{\oplus n_s})\mapsto \bfM(\theta^{r}_{\ob{i}}) (\bigoplus_s X_s^{\oplus n_s})\\
(\lunthree_{\bfM}^*)_{\ob{i}} = \lun^*_{\bfM(\ob{i}_r)} \bfM((\theta^{\ell}_{\ob{i}})^*)&\colon \bfM(\ob{i})  \to \Add\Bbbk\monob \bfM(\ob{i}_r), \quad X\mapsto (\bullet_\Bbbk, \bfM(\theta^{\ell,*}_{\ob{i}}) (X))\\
(\runthree_{\bfM}^*)_{\ob{i}} = \run^*_{\bfM(\ob{i}_\ell)} \bfM((\theta^{r}_{\ob{i}})^*)&\colon \bfM(\ob{i})  \to\bfM(\ob{i}_\ell) \monob \Add\Bbbk, \quad X\mapsto ( \bfM(\theta^{r,*}_{\ob{i}}) (X),\bullet_\Bbbk).\\
\end{align*}

\begin{prop}\label{unitorthree}
\begin{enumerate}[(a)]
\item\label{unitorthreeleft} There is an adjoint equivalence $(\lunthree, \lunthree^*, \eta^{\lunthree}, \varepsilon^{\lunthree})$
$$\tenrepC( \One_{\Cc\amod} \times \mathbbm{1}_{\Cc\amod})\underset{\lunthree^*}{ \overset{\lunthree}{\rightleftarrows}}\mathbbm{1}_{\Cc\amod} $$ 
between $2$-functors $\Cc\amod\to \Cc\amod$.
\item\label{unitorthreeright}  There is an adjoint equivalence $(\runthree, \runthree^*, \eta^{\runthree}, \varepsilon^{\runthree})$
$$ \tenrepC(\mathbbm{1}_{\Cc\amod}\times \One_{\Cc\amod})\underset{\runthree^*}{ \overset{\runthree}{\rightleftarrows}} \mathbbm{1}_{\Cc\amod} $$
between $2$-functors $\Cc\amod\to \Cc\amod$.
\end{enumerate}
\end{prop}

\proof
\eqref{unitorthreeleft} By definition, such an adjoint equivalence consists of strong natural transformations $\lunthree, \lunthree^*$ and modifications $\eta^{\lunthree}\colon\mathbbm{1}_{ \mathbbm{1}_{\Cc\amod}}\to \lunthree\,\lunthree^*$, $\varepsilon^{\lunthree} \colon  \lunthree^*\lunthree\to  \mathbbm{1}_{\tenrepC( \One_{\Cc\amod} \times \mathbbm{1}_{\Cc\amod})}$ satisfying the triangle equality. 

Such strong natural transformations consist of, for each object $\bfM$ in $\Cc\amod$, a $1$-cell in $\Cc\amod$, and these are given by  the $1$-cells $\lunthree_{\bfM},\lunthree^*_{\bfM}$ in $\Cc\amod$. 
Moreover, for every $1$-cell $\Phi\colon \bfM\to \bfN$ in $\Cc\amod$, we need modifications 
$$\lunthree_{\Phi}\colon \Phi \lunthree_{\bfM}\to\lunthree_{\bfN}  (\mathbbm{1}_{\bfT_\Cc}\tenrepC_1 \Phi) $$
$$\lunthree^*_{\Phi}\colon  (\mathbbm{1}_{\bfT_\Cc}\tenrepC_1 \Phi)   \lunthree_{\bfM}^* \to \lunthree_{\bfN}^* \Phi$$
compatible with compositions of such $1$-cells. The first of these consists of, for each object $\ob{i}\in \Cc$, a natural transformation 
$(\lunthree_{\Phi})_{\ob{i}}\colon  \Phi_{\ob{i}} (\lunthree_{\bfM})_{\ob{i}}\to(\lunthree_{\bfN})_{\ob{i}}  (\mathbbm{1}_{\bfT_\Cc}\tenrepC_1 \Phi)_{\ob{i}}$. Unravelling the definitions, for an object $X\in \bfM(\ob{i}_r)$
$$(\Phi_{\ob{i}} (\lunthree_{\bfM})_{\ob{i}})(\bullet_{\Bbbk}^{\oplus n},X)
 =  \Phi_{\ob{i}} \bfM(\theta^{\ell}_{\ob{i}}) (X^{\oplus n}) $$
 while
$$( (\lunthree_{\bfN})_{\ob{i}}  (\mathbbm{1}_{\bfT_\Cc}\tenrepC_1 \Phi)_{\ob{i}}  )(\bullet_{\Bbbk}^{\oplus n},X) =  (\lunthree_{\bfN})_{\ob{i}}(\bullet_{\Bbbk}^{\oplus n},\Phi_{\ob{i}_r} (X))
 = \bfN(\theta^{\ell}_{\ob{i}})(\Phi_{\ob{i}_r} (X)^{\oplus n})$$
 so 
$((\lunthree_{\Phi})_{\ob{i}})_{(\bullet_{\Bbbk}^{\oplus n},X)}$ is induced by $(\Phi_{\theta^{\ell}_{\ob{i}}})_{X^{\oplus n}}$.
 
Similarly, the second of these consists of, for each object $\ob{i}\in \Cc$, a natural transformation $$(\lunthree^*_{\Phi})_{\ob{i}}\colon  (\mathbbm{1}_{\bfT_\Cc}\tenrepC_1 \Phi)_{\ob{i}}   (\lunthree_{\bfM}^*)_{\ob{i}} \to (\lunthree_{\bfN}^*)_{\ob{i}} \Phi_{\ob{i}}.$$ On an object $X\in \bfM(\ob{i})$, we note that  
$$((\mathbbm{1}_{\bfT_\Cc}\tenrepC_1 \Phi)_{\ob{i}}   (\lunthree_{\bfM}^*)_{\ob{i}} )(X)=  (\bullet_{\Bbbk},\Phi_{\ob{i}_r} \bfM(\theta^{\ell*}_{\ob{i}}) (X))$$ while
$$((\lunthree_{\bfN}^*)_{\ob{i}} \Phi_{\ob{i}})(X) = (\bullet_{\Bbbk}, \bfN(\theta^{\ell*}_{\ob{i}}) \Phi_{\ob{i}}(X)) $$
so we can set 
$((\lunthree^*_{\Phi})_{\ob{i}})_X = \mathbbm{1}\otimes \Phi_{\theta^{\ell*}_{\ob{i}}}$. Naturality of the resulting transformations, the modification axiom, as well as compatibility with the composition of $1$-cells in $\Cc\amod$ (i.e.\ lax naturality) are then immediate.

Finally, the modifications $\eta^{\lunthree}$ and $\varepsilon^{\lunthree}$ consist of  $2$-cells $\eta^{\lunthree}_\bfM$ and $\varepsilon^{\lunthree}_\bfM$ in $\Cc\amod$, i.e.\ modifications, for each $\bfM\in \Cc\amod$. Explicitly, noting that, on an object $\ob{i}\in \Cc$,
$$(\lunthree_\bfM)_{\ob{i}}(\lunthree^*_\bfM)_{\ob{i}} = \bfM(\theta^{\ell}_{\ob{i}}) \lun_{\bfM(\ob{i}_r)}\lun^*_{\bfM(\ob{i}_r)} \bfM(\theta^{\ell*}_{\ob{i}})$$
and
$$(\lunthree^*_\bfM)_{\ob{i}}(\lunthree_\bfM)_{\ob{i}} = \lun^*_{\bfM(\ob{i}_r)} \bfM(\theta^{\ell*}_{\ob{i}})\bfM(\theta^{\ell}_{\ob{i}}) \lun_{\bfM(\ob{i}_r)}$$
we see that 
$(\eta^{\lunthree}_\bfM)_{\ob{i}}$ and 
$(\varepsilon^{\lunthree}_\bfM)_{\ob{i}}  $ are assembled from $\eta^\lun_{\bfM(\ob{i}_r)}$, resp.\ $\varepsilon^\lun_{\bfM(\ob{i}_r)}$, and the adjunction $2$-morphisms for the adjoint pair $(\theta^{\ell}_{\ob{i}},\theta^{\ell*}_{\ob{i}})$ of $1$-morphisms in $\Cc$, to which we apply $\bfM$ (using the compositionality constraint $\bfM^2$ of the latter). Again, the facts that these assemble to modifications and satisfy the triangle identities follow from the properties of the ingredients.

\eqref{unitorthreeright} Mutatis mutandis.
\endproof

\subsection{Associators for tensor products of strictly unital birepresentations}

Analogously to the unitors, we construct an intermediate adjoint equivalence, and use these to define the data of the associators of $\what$, via the pasting diagram in Definition~\ref{assthreelocal}. For this, recall the construction of $\asstwo_{\Aa,\Bb,\Cc}$ from Proposition~\ref{ass2fun}, and $\ass_{\oneA,\oneB,\oneC}, \ass^*_{\oneA,\oneB,\oneC}$ from \eqref{pa:associator}.

\begin{prop}\label{alpha}
There is an adjoint equivalence $(\alpha,\alpha^*, \varepsilon^\alpha, \eta^{\alpha})$ 
$$\what(\what\bxone\mathbbm{1}_{\A})  \underset{\alpha^*}{ \overset{\alpha}{\rightleftarrows}}   \what(\mathbbm{1}_{\A}\bxone\what)  \asstwo_{\A,\A,\A}$$
between $2$-functors $(\A\bxob\A)\bxob\A \to \A$.
\end{prop}

\proof
We first need to construct the strong natural transformations $\alpha,\alpha^*$. Given that, for an object $ ((\oneB,\oneC),\oneD)$ in $(\A\bxob\A)\bxob\A$, we have 
\begin{align*}\what(\what\bxone\mathbbm{1}_{\A}) ((\oneB,\oneC),\oneD)  &= (\oneB\monob\oneC)\monob\oneD  \qquad \text{and}\\ \what(\mathbbm{1}_{\A}\bxone\what)  \asstwo_{\A,\A,\A}  (\oneB,(\oneC,\oneD) ) &= \oneB\monob(\oneC\monob\oneD),\end{align*}
we can choose the $1$-cell components of $\alpha,\alpha^*$ on such an object to be $\ass_{\oneB,\oneC,\oneD}, \ass^*_{\oneB,\oneC,\oneD}$, respectively.

Let $\F= (\bigoplus_p (\bigoplus_q (\F_{1pq},\F_{2pq}),\F_{3p}))$ be a $1$-morphism from $((\oneB,\oneC),\oneD) \to((\oneB',\oneC'),\oneD')$ in the $2$-category $(\A\bxob\A)\bxob\A$. Noting that
\begin{equation*}\begin{split}
\what(\what\bxone\mathbbm{1}_{\A})(\F) &=  \bigoplus_p (\bigoplus_q (\F_{1pq}\monmor\F_{2pq}) \monmor \F_{3p}) \qquad \text{while}\\
\what(\mathbbm{1}_{\A}\bxone\what)  \asstwo_{\A,\A,\A} (\F) 
&= \what(\mathbbm{1}_{\A}\bxone\what) (\bigoplus_p \bigoplus_q (\F_{1pq},(\F_{2pq},\F_{3p})))\\
&=\bigoplus_p \bigoplus_q (\F_{1pq}\monmor(\F_{2pq}\monmor\F_{3p}))
\end{split}\end{equation*}
for an object $X=\bigoplus_i(\bigoplus_j(X_{ij},Y_{ij}),Z_i)$ in  $ (\oneB\monob\oneC)\monob\oneD$, we compute 
\begin{equation*}\begin{split}
		\ass_{\oneB',\oneC',\oneD'}\left(\what(\what\bxone\mathbbm{1}_{\A})(\F)\right)(X)&=\ass_{\oneB',\oneC',\oneD'}\left( \bigoplus_p (\bigoplus_q (\F_{1pq}\monmor\F_{2pq}) \monmor \F_{3p})\right)(X)\\
		&=\ass_{\oneB',\oneC',\oneD'} ( \bigoplus_p \bigoplus_i (\bigoplus_q \bigoplus_j (\F_{1pq}X_{ij},\F_{2pq} Y_{ij}), \F_{3p}Z_i))\\
		&=\bigoplus_p\bigoplus_i\bigoplus_q\bigoplus_j (\F_{1pq}X_{ij},(\F_{2pq} Y_{ij}, \F_{3p}Z_{i}))
\end{split}\end{equation*}
and 
\begin{equation*}\begin{split}
		\left(\what(\mathbbm{1}_{\A}\bxone\what)  \asstwo_{\A,\A,\A} (\F) \right)\ass_{\oneB,\oneC,\oneD} (X)&=\left(\what(\mathbbm{1}_{\A}\bxone\what)  \asstwo_{\A,\A,\A} (\F) \right) (\bigoplus_i \bigoplus_j (X_{ij},(Y_{ij},Z_{i})))\\
		&= \bigoplus_p \bigoplus_q (\F_{1pq}\monmor(\F_{2pq}\monmor\F_{3p})) (\bigoplus_i \bigoplus_j (X_{ij},(Y_{ij},Z_{i})))\\
		&=\bigoplus_p\bigoplus_q \bigoplus_i\bigoplus_j(\F_{1pq}X_{ij},(\F_{2pq} Y_{ij}, \F_{3p}Z_{i})).
\end{split}\end{equation*}

Since it also follows from the universal property of the coproduct that the map reordering the summation commutes with all morphisms in $B\monob(C\monob D)$,  we may choose  the components of the $2$-cells associated to $\alpha$ to be this reordering. It is straightforward to check that this choice satisfies lax naturality.

Given the $1$-cell components of $\alpha, \alpha^*$, it is easy to see that, on an object $ ((\oneB,\oneC),\oneD)$ in $(\A\bxob\A)\bxob\A$,  the $2$-cells $ \varepsilon^\alpha_{\oneB,\oneC,\oneD} , \eta^{\alpha}_{\oneB,\oneC,\oneD} $ defining the adjunction modifications $ \varepsilon^\alpha, \eta^{\alpha}$ can be chosen as the natural transformations $ \varepsilon^\ass_{\oneB,\oneC,\oneD} , \eta^{\ass}_{\oneB,\oneC,\oneD} $. That these satisfy the triangle identities reduces to checking that $\ass$, $\ass^*$, $ \varepsilon^\ass_{\oneB,\oneC,\oneD}$, and $ \eta^{\ass}_{\oneB,\oneC,\oneD} $ satisfy the triangle identities.
\endproof

\begin{defn}\label{assthreelocal} For $\bfL,\bfM,\bfN\in \Cc\amod$, the pasting diagram
$$
\xymatrix{
&&\Cc\ar_{\Delta}[dll]\ar^{\Delta}[drr] &&\\
\Cc\bxob\Cc \ar_{\Delta\bxone  \mathbbm{1}_{\Cc}}[dd] &&&\ar@2{->}@<5pt>^{\delta^*}[dll]&\Cc\bxob\Cc \ar^{ \mathbbm{1}_{\Cc}\bxone \Delta}[dd]\\
&\ar@2{->}@<5pt>^{\delta}[urr]&&&\\
(\Cc\bxob\Cc)\bxob\Cc\ar^{\asstwo_{\Cc,\Cc,\Cc}}[rrrr]\ar_{(\bfL\bxone \bfM)\bxone \bfN}[dd]&&&&\Cc\bxob(\Cc\bxob\Cc)\ar^{\bfL\bxone (\bfM\bxone \bfN)}[dd]\\
&&&&\\
(\A\bxob\A)\bxob\A\ar^{\asstwo_{\A,\A,\A}}[rrrr]\ar_{\what(\what\bxone\mathbbm{1}_{\A})}[ddddrr]&&&&\A\bxob(\A\bxob\A)\ar^{ \what(\mathbbm{1}_{\A}\bxone\what) }[ddddll]\\
&&&\ar@2{->}@<5pt>^{\alpha^*}[dll]&\\
&\ar@2{->}@<5pt>^{\alpha}[urr]&&&\\
&&&&\\
&&\A&&
}
$$
defines $1$-cells $$\assthree_{\bfL,\bfM,\bfN}\colon \;\;(\bfL\tenrepC_0 \bfM)\tenrepC_0 \bfN \rightleftarrows \bfL\tenrepC_0 (\bfM\tenrepC_0 \bfN)\;\;\colon \assthree^*_{\bfL,\bfM,\bfN}$$
in $\Cc\amod$.
\end{defn}

Explicitly, on an object $\ob{i}\in \Cc$, the functor  $(\assthree_{\bfL,\bfM,\bfN})_{\ob{i}}$ is given by the functor
{\small \begin{align*}(\bfL\tenrep_0 (\bfM\tenrep_0 \bfN))(\delta_{\ob{i}}) \ass_{\bfL(\ob{i}_{\ell\ell}), \bfM(\ob{i}_{\ell r}),\bfN(\ob{i}_{r})} \colon (\bfL(\ob{i}_{\ell\ell})\monob \bfM(\ob{i}_{\ell r}))\monob \bfN(\ob{i}_{r}) &\to \bfL(\ob{i}_{\ell})\monob( \bfM(\ob{i}_{r \ell })\monob \bfN(\ob{i}_{rr}))
\end{align*}
}
while $(\assthree^*_{\bfL,\bfM,\bfN})_{\ob{i}}$ is given by 
$$\ass^*_{\bfL(\ob{i}_{\ell\ell}), \bfM(\ob{i}_{\ell r}),\bfN(\ob{i}_{r})} (\bfL\tenrep_0 (\bfM \tenrep_0 \bfN))(\delta^*_{\ob{i}})\colon \bfL(\ob{i}_{\ell})\monob( \bfM(\ob{i}_{r \ell })\monob \bfN(\ob{i}_{rr})) \to  (\bfL(\ob{i}_{\ell\ell})\monob \bfM(\ob{i}_{\ell r}))\monob \bfN(\ob{i}_{r}). $$

\begin{prop}\label{assthree}
There is an adjoint equivalence $(\assthree, \assthree^*, \eta^{\assthree}, \varepsilon^{\assthree})$
$$ \tenrepC(\tenrepC \times \mathbbm{1}_{\Cc\amod})\underset{\assthree^*}{ \overset{\assthree}{\rightleftarrows}} \tenrepC(\mathbbm{1}_{\Cc\amod} \times \tenrepC)$$ between $2$-functors $\Cc\amod\times\Cc\amod\times\Cc\amod\to \Cc\amod$.
\end{prop}

\proof 
 By definition, such an adjoint equivalence consists of strong natural transformations $\assthree, \assthree^*$ and modifications 
 $\eta^{\assthree}\colon \mathbbm{1}_{\tenrepC(\tenrepC \times \mathbbm{1}_{\Cc\amod}) }\to\assthree^*\,\assthree$, $\varepsilon^{\assthree} \colon  \assthree\assthree^*\to  \mathbbm{1}_{\tenrepC(\mathbbm{1}_{\Cc\amod} \times \tenrepC) }$ satisfying the triangle equality.
 
 The $1$-cells forming part of the data of the required strong natural transformations are given by, for each triple $\bfL,\bfM,\bfN$ of objects in $\Cc\amod$, the $1$-cells  $\assthree_{\bfL,\bfM,\bfN}, \assthree^*_{\bfL,\bfM,\bfN}$, respectively. To define the necessary $2$-cells, let $\Phi\colon \bfL\to\bfL', \Psi\colon \bfM\to\bfM', \Sigma\colon \bfN\to\bfN'$ be $1$-cells in $\Cc\amod$. We need modifications 
  $$\assthree_{\Phi,\Psi,\Sigma} \colon  (\Phi \tenrepC_1(\Psi\tenrepC_1\Sigma))  \assthree_{\bfL,\bfM,\bfN}\to \assthree_{\bfL',\bfM',\bfN'} ((\Phi \tenrepC_1\Psi)\tenrepC_1\Sigma)$$
and 
$$\assthree^*_{\Phi,\Psi,\Sigma} \colon ((\Phi \tenrepC_1\Psi)\tenrepC_1\Sigma)  \assthree^*_{\bfL,\bfM,\bfN} \to \assthree^*_{\bfL',\bfM',\bfN'} (\Phi \tenrepC_1(\Psi\tenrepC_1\Sigma)) .$$
On an object $\ob{i}\in \Cc$, 
$$
\left((\Phi \tenrepC_1(\Psi\tenrepC_1\Sigma))  \assthree_{\bfL,\bfM,\bfN}\right)_{\ob{i}}
=(\Phi_{\ob{i}_{\ell}} \monmor(\Psi_{\ob{i}_{r\ell }}\monmor\Sigma_{\ob{i}_{rr}}))
(\bfL\tenrep_0 (\bfM\tenrep_0 \bfN))(\delta_{\ob{i}}) \ass_{\bfL(\ob{i}_{\ell\ell}), \bfM(\ob{i}_{\ell r}),\bfN(\ob{i}_{r})}.
$$
while 
$$\left( \assthree_{\bfL',\bfM',\bfN'} ((\Phi \tenrepC_1\Psi)\tenrepC_1\Sigma)\right)_{\ob{i}} = 
(\bfL'\tenrep_0 (\bfM'\tenrep_0 \bfN'))(\delta_{\ob{i}}) \ass_{\bfL'(\ob{i}_{\ell\ell}), \bfM'(\ob{i}_{\ell r}),\bfN'(\ob{i}_{r})} 
((\Phi_{\ob{i}_{\ell\ell}} \monmor\Psi_{\ob{i}_{\ell r}}) \monmor\Sigma_{\ob{i}_{r}}) $$

Similarly,
$$
\left(((\Phi \tenrepC_1\Psi)\tenrepC_1\Sigma)  \assthree^*_{\bfL,\bfM,\bfN} \right)_{\ob{i}} =
((\Phi_{\ob{i}_{\ell\ell}} \monmor\Psi_{\ob{i}_{\ell r}})\monmor\Sigma_{\ob{i}_{r}})
 \ass^*_{\bfL(\ob{i}_{\ell\ell}), \bfM(\ob{i}_{\ell r}),\bfN(\ob{i}_{r})} (\bfL\tenrep_0 (\bfM \tenrep_0 \bfN))(\delta^*_{\ob{i}}) 
$$
while
$$ 
\left(\assthree^*_{\bfL',\bfM',\bfN'} (\Phi \tenrepC_1(\Psi\tenrepC_1\Sigma)) \right)_{\ob{i}} 
=
\ass^*_{\bfL'(\ob{i}_{\ell\ell}), \bfM'(\ob{i}_{\ell r}),\bfN'(\ob{i}_{r})} (\bfL'\tenrep_0 (\bfM' \tenrep_0 \bfN'))(\delta^*_{\ob{i}})  (\Phi_{\ob{i}_{\ell}} \monmor(\Psi_{\ob{i}_{r\ell}}\monmor\Sigma_{\ob{i}_{rr}})).
$$


Since, as argued in Proposition \ref{assone},  $\ass$ is a strict natural transformation, we have
$$ 
\left(((\Phi \tenrepC_1\Psi)\tenrepC_1\Sigma)  \assthree^*_{\bfL,\bfM,\bfN} \right)_{\ob{i}} =
\ass^*_{\bfL(\ob{i}_{\ell\ell}), \bfM(\ob{i}_{\ell r}),\bfN(\ob{i}_{r})}
(\Phi_{\ob{i}_{\ell\ell}} \monmor(\Psi_{\ob{i}_{\ell r}}\monmor\Sigma_{\ob{i}_{r}}))
  (\bfL\tenrep_0 (\bfM \tenrep_0 \bfN))(\delta^*_{\ob{i}}) 
$$
and
$$\left( \assthree_{\bfL',\bfM',\bfN'} ((\Phi \tenrepC_1\Psi)\tenrepC_1\Sigma)\right)_{\ob{i}} = (\bfL'\tenrep_0 (\bfM'\tenrep_0 \bfN'))(\delta_{\ob{i}})
(\Phi_{\ob{i}_{\ell\ell}} \monmor(\Psi_{\ob{i}_{\ell r}} \monmor\Sigma_{\ob{i}_{r}})) \ass_{\bfL(\ob{i}_{\ell\ell}), \bfM(\ob{i}_{\ell r}),\bfN(\ob{i}_{r})}$$
and hence $\assthree_{\Phi,\Psi,\Sigma},\assthree^*_{\Phi,\Psi,\Sigma}$ are simply induced from the commutation data associated to the $1$-cells in $\Cc\amod$, providing isomorphisms
$$
(\Phi_{\ob{i}_{\ell}} \monmor(\Psi_{\ob{i}_{r\ell }}\monmor\Sigma_{\ob{i}_{rr}}))
(\bfL\tenrep_0 (\bfM\tenrep_0 \bfN))(\delta_{\ob{i}}) \to (\bfL'\tenrep_0 (\bfM'\tenrep_0 \bfN'))(\delta_{\ob{i}})
 (\Phi_{\ob{i}_{\ell\ell}} \monmor(\Psi_{\ob{i}_{\ell r}} \monmor\Sigma_{\ob{i}_{r}}))
$$
and
$$ (\bfL'\tenrep_0 (\bfM' \tenrep_0 \bfN'))(\delta^*_{\ob{i}})  (\Phi_{\ob{i}_{\ell}} \monmor(\Psi_{\ob{i}_{r\ell}}\monmor\Sigma_{\ob{i}_{rr}})) \to(\Phi_{\ob{i}_{\ell\ell}} \monmor(\Psi_{\ob{i}_{\ell r}}\monmor\Sigma_{\ob{i}_{r}}))
  (\bfL\tenrep_0 (\bfM \tenrep_0 \bfN))(\delta^*_{\ob{i}}) . $$
It is straightforward to check that this is compatible with composition of $1$-cells in $\Cc\amod$ and hence the $1$-cells $\assthree_{\bfL,\bfM,\bfN}, \assthree^*_{\bfL,\bfM,\bfN}$ and $2$-cells $\assthree_{\Phi,\Psi,\Sigma},\assthree^*_{\Phi,\Psi,\Sigma}$ in $\Cc\amod$ assemble to strong natural transformations $\assthree, \assthree^*$ as required.

The modification $\eta^{\assthree}$ consists, for each triple $\bfL,\bfM,\bfN$ of objects in $\Cc\amod$, a $2$-cell 
$$\eta^{\assthree}_{\bfL,\bfM,\bfN}  \colon \mathbbm{1}_{(\bfL\tenrepC_0 \bfM)\tenrepC_0 \bfN} \to \assthree^*_{\bfL,\bfM,\bfN}\assthree_{\bfL,\bfM,\bfN} $$
in $\Cc\amod$. Such a $2$-cell is a modification between $1$-cells in $\Cc\amod$, which is given on an object $\ob{i}\in \Cc$ by the natural transformation
$$  \mathbbm{1}_{(\bfL(\ob{i}_{\ell}) \monob\bfM(\ob{i}_{r\ell}))\monob \bfN(\ob{i}_{rr})} \to \ass^*_{\bfL(\ob{i}_{\ell\ell}), \bfM(\ob{i}_{\ell r}),\bfN(\ob{i}_{r})}(\bfL\tenrep_0 (\bfM\tenrep_0 \bfN))(\delta^*_{\ob{i}}) (\bfL\tenrep_0 (\bfM \tenrep_0 \bfN))(\delta_{\ob{i}})\ass_{\bfL(\ob{i}_{\ell\ell}), \bfM(\ob{i}_{\ell r}),\bfN(\ob{i}_{r})}, $$
which is obtained by appropriately composing 
 $\eta^\ass_{\bfL(\ob{i}_{\ell\ell}), \bfM(\ob{i}_{\ell r}),\bfN(\ob{i}_{r})}$ from Proposition \ref{assone}, with $(\bfL\tenrep_0 (\bfM\tenrep_0 \bfN))((\eta^\delta)_\ob{i})$ where $\eta^\delta$ is the unit of the adjunction $(\delta_{\ob{i}},\delta^*_{\ob{i}})$.

Similarly, the modification $\varepsilon^{\assthree}$ consists, for each triple $\bfL,\bfM,\bfN$ of objects in $\Cc\amod$, a $2$-cell 
$$\varepsilon^{\assthree}_{\bfL,\bfM,\bfN}\colon  \assthree_{\bfL,\bfM,\bfN}\assthree^*_{\bfL,\bfM,\bfN}   \to\mathbbm{1}_{\bfL\tenrepC_0 (\bfM\tenrepC_0 \bfN)} $$
 in $\Cc\amod$, which is given, on an object $\ob{i}\in \Cc$, by the natural transformation
$$
(\bfL\tenrepC_0 (\bfM\tenrepC_0 \bfN))
(\delta_{\ob{i}})
 \ass_{\bfL(\ob{i}_{\ell\ell}), \bfM(\ob{i}_{\ell r}),\bfN(\ob{i}_{r})}
 \ass^*_{\bfL(\ob{i}_{\ell\ell}), \bfM(\ob{i}_{\ell r}),\bfN(\ob{i}_{r})}
(\bfL\tenrepC_0 (\bfM \tenrepC_0 \bfN))(\delta^*_{\ob{i}}) 
 \to \mathbbm{1}_{\bfL(\ob{i}_{\ell})\monob(\bfM(\ob{i}_{r\ell})\monob \bfN(\ob{i}_{rr}))}
$$
 induced similarly to the previous paragraph by $\varepsilon^\ass$ from Proposition \ref{assone} and the counit $\varepsilon^\delta$ of the adjunction $(\delta_{\ob{i}},\delta^*_{\ob{i}})$.
 

That both triangle identities are satisfied componentwise follows from the fact that they are satisfied for $\eta^\ass$ and $\varepsilon^\ass$, as well as for $\eta^\delta$, and $\varepsilon^\delta$.  
\endproof

\subsection{Pentagonator and $2$-unitors}

In order to obtain a monoidal structure on the $2$-category $\Cc\amod$, we also require 
\begin{itemize}
\item an invertible modification $\muthree$ in $\Bb icat^{ps}(\Cc\amod\times \Cc\amod,\Cc\amod)$, given by invertible modifications 
$$(\mathbbm{1}_{\bfM}\tenrepC_1 \lunthree_{\bfN})\, \assthree_{\bfM,\bfT_\Cc, \bfN} \,(\runthree^*_{\bfM}\tenrepC_1\mathbbm{1}_{\bfN} )\xrightarrow{\muthree_{\bfM,\bfN}}\mathbbm{1}_{\bfM\tenrepC_0\bfN}$$
\item an invertible modification $\lambdathree$ in $\Bb icat^{ps}(\Cc\amod\times \Cc\amod,\Cc\amod)$ given by modifications 
$$\lunthree_{\bfM}\tenrepC_1 \mathbbm{1}_{\bfN}\xrightarrow{\lambdathree_{\bfM,\bfN}}\lunthree_{\bfM\tenrepC_0\bfN}\, \assthree_{\bfT_\Cc,\bfM,\bfN } $$
\item an invertible modification $\rhothree$ in $\Bb icat^{ps}(\Cc\amod\times \Cc\amod,\Cc\amod)$ given by modifications 
$$\mathbbm{1}_{\bfM}\tenrepC_1 \runthree^*_{\bfN}
 \xrightarrow{\rhothree_{\bfM,\bfN}} \assthree_{\bfM,\bfN,\bfT_\Cc } \,\runthree^*_{\bfM\tenrepC_0\bfN},$$
 \item an invertible modification $\pithree$ in $\Bb icat^{ps}(\Cc\amod^{\times 4},\Cc\amod)$ given by modifications 
 $$(\mathbbm{1}_{\bfK}\tenrepC_1 \assthree_{ \bfL, \bfM,\bfN})\assthree_{\bfK, \bfL\tenrepC_0 \bfM,\bfN}(\assthree_{\bfK, \bfL, \bfM}\tenrepC_1\mathbbm{1}_{\bfN})
 \xrightarrow{\pithree_{\bfK, \bfL,\bfM,\bfN}} \assthree_{\bfK, \bfL,\bfM\tenrepC_0\bfN } \,\assthree_{\bfK\tenrepC_0 \bfL,\bfM,\bfN}.$$
\end{itemize}


We will now consider the $1$-cells in the domain of $\muthree_{\bfM,\bfN}$ and the domains and codomains of $\lambdathree_{\bfM,\bfN}$ and $\rhothree_{\bfM,\bfN}$. In each case, we will split their pasting diagrams into a part only involving the $2$-category and a part involving the birepresentations $\bfM$ and $\bfN$, and show that the parts involving the birepresentations $\bfM$ and $\bfN$ may be assumed to be identities.

In each case, we will need to swap the order of certain strong transformations, whose non-identity parts are separated by a $2$-functor.

\begin{lemma}\label{swaptwocells}
Suppose we have $2$-categories, pseudofunctors and strong natural transformations as follows:
\[\begin{tikzcd}[ampersand replacement=\&]
	\Aa \& \Bb \& \Cc \& \Dd
	\arrow[""{name=0, anchor=center, inner sep=0}, "\F", curve={height=-18pt}, from=1-1, to=1-2]
	\arrow[""{name=1, anchor=center, inner sep=0}, "{\F'}"', curve={height=18pt}, from=1-1, to=1-2]
	\arrow["\G"', from=1-2, to=1-3]
	\arrow[""{name=2, anchor=center, inner sep=0}, "\H", curve={height=-18pt}, from=1-3, to=1-4]
	\arrow[""{name=3, anchor=center, inner sep=0}, "{\H'}"', curve={height=18pt}, from=1-3, to=1-4]
	\arrow["\alpha", between={0.2}{0.8}, Rightarrow, from=0, to=1]
	\arrow["\beta", between={0.2}{0.8}, Rightarrow, from=2, to=3]
\end{tikzcd}\]
 where we additionally assume that $\H$ is a $2$-functor.
There is a modification 
\[
((\H'\G)\star\alpha)
(\beta  \star (\G\F)) 
\to
(\beta  \star (\G\F'))
((\H\G)\star\alpha)
\]
whose $2$-cell component at an object $\ob{i}\in \Aa$  is given by 
 $\beta_{G(\alpha_\ob{i})}$.
\end{lemma}

\proof
The composition
\[
((\H'\G)\star\alpha)
(\beta  \star (\G\F))
\]
has component at an object $\ob{i}\in \Aa$ given by 
$\H'\G(\alpha_{\ob{i}}) \beta_{\G\F(\ob{i})}$
whilst the corresponding component of the composition 
\[
(\beta  \star (\G\F'))
((\H\G)\star\alpha)
\]
is given by 
$\beta_{\G\F'(\ob{i})}\H\G(\alpha_{\ob{i}})$.
It is straightforward to check that the morphisms \[\beta_{G(\alpha_\ob{i})}\colon H'G(\alpha_\ob{i})\beta_{GF(\ob{i})}\Rightarrow \beta_{GF'(\ob{i})}H G(\alpha_\ob{i}),\] for each object $\ob{i}$ in $\Aa$, yield the modifiction between the two choices of composition.
\endproof

\subsubsection{The $2$-unitor $\muthree$} (Complete up to a few minor corrections and marked to dos.)
First we consider the domain $(\mathbbm{1}_{\bfM}\tenrepC_1 \lunthree_{\bfN})\, \assthree_{\bfM,\bfT_\Cc, \bfN} \,(\runthree^*_{\bfM}\tenrepC_1\mathbbm{1}_{\bfN} )$ of $\muthree_{\bfM,\bfN}$. We will show that this can be represented by the following pasting diagram. This diagram does not define the strong natural transformation uniquely, since it takes place in a tricategory where the interchange law does not hold strictly, although the visualisation via the pasting diagram will nevertheless be helpful. 

\begin{align}\label{middleunitor2cell}\scalebox{0.8}{
\begin{tikzcd}[ampersand replacement=\&]
	\&\&\& \Cc \&\&\& \\
	\\
	\& {\Cc\bxone \Cc} \&\&\&\& {\Cc\bxone \Cc} \\
	\\
	\&\& {(\Cc\bxob \Cc )\bxob \Cc} \&\& {\Cc\bxob(\Cc\bxob \Cc )} \\
	\\
	{\Cc\bxob \Cc} \&\& {(\Cc\bxob \Vv )\bxob \Cc} \&\& {\Cc\bxob(\Vv\bxob \Cc )} \&\& {\mathcal{C}\boxtimes_0 \mathcal{C}} \\
	{\Cc\bxob \Cc} \& {\Cc\bxob(\Vv\bxob \Cc )} \& {(\Cc\bxob \Vv )\bxob \Cc} \\
	{\Cc\bxob \Cc} \&\& {(\Cc\bxob \Vv )\bxob \Cc} \\
	\\
	\&\& {(\A\bxob  \A)\bxob \A} \&\& {\A\bxob(\A\bxob \A)} \\
	\\
	\& {\A\bxob \A} \&\&\&\& {\A\bxob \A} \\
	\\
	\&\&\& \A
	\arrow["\Delta"', from=1-4, to=3-2]
	\arrow["\Delta", from=1-4, to=3-6]
	\arrow["\delta", between={0.4}{0.6}, Rightarrow, from=3-2, to=3-6]
	\arrow["{\Delta\bxone\mathbbm{1}_\Cc}", from=3-2, to=5-3]
	\arrow[""{name=0, anchor=center, inner sep=0}, equals, from=3-2, to=7-1]
	\arrow["{\mathbbm{1}_\Cc\bxone\Delta}"', from=3-6, to=5-5]
	\arrow[""{name=1, anchor=center, inner sep=0}, equals, from=3-6, to=7-7]
	\arrow["{\asstwo_{\Cc,\Cc,\Cc}}", from=5-3, to=5-5]
	\arrow[""{name=2, anchor=center, inner sep=0}, "{(\mathbbm{1}_\Cc\bxone\Theta)\bxone \mathbbm{1}_\Cc}"', from=5-3, to=7-3]
	\arrow[""{name=3, anchor=center, inner sep=0}, "{\mathbbm{1}_\Cc\bxone(\Theta\bxone \mathbbm{1}_\Cc)}", from=5-5, to=7-5]
	\arrow[equals, from=7-1, to=8-1]
	\arrow[between={0.4}{0.6}, equals, from=7-3, to=5-5]
	\arrow["{\runtwo_\Cc\bxone \mathbbm{1}_\Cc}", from=7-3, to=7-1]
	\arrow["{\asstwo_{\Cc,\Vv,\Cc}}", from=7-3, to=7-5]
	\arrow["{\tilde\mutwo_{\Cc,\Cc}}"', between={0.4}{0.6}, Rightarrow, dashed, from=7-3, to=8-1]
	\arrow[equals, from=7-3, to=8-3]
	\arrow["{\mathbbm{1}_\Cc\bxone\luntwo_\Cc}"', from=7-5, to=7-7]
	\arrow["{\bfM\bxone(\bfT_\Vv \bxone \bfN)}", from=7-5, to=11-5]
	\arrow["{\bfM\bxone\bfN}", from=7-7, to=13-6]
	\arrow[equals, from=8-1, to=9-1]
	\arrow["{\mathbbm{1}_\Cc\bxone\luntwo_\Cc}"', dashed, from=8-2, to=8-1]
	\arrow["{\asstwo_{\Cc,\Vv,\Cc}}"', dashed, from=8-3, to=8-2]
	\arrow["{(\tilde\mutwo_{\Cc,\Cc})^{-1}}"', between={0.4}{0.6}, Rightarrow, dashed, from=8-3, to=9-1]
	\arrow[equals, from=8-3, to=9-3]
	\arrow["{\tau^{r*}_\bfM\bxtwo \mathbbm{1}_{\bfN}}"'{pos=0.7}, between={0.4}{0.6}, Rightarrow, from=9-1, to=11-3]
	\arrow["{\bfM\bxone\bfN}"', from=9-1, to=13-2]
	\arrow["{\runtwo_\Cc\bxone \mathbbm{1}_\Cc}", from=9-3, to=9-1]
	\arrow["{(\bfM\bxone\bfT_\Vv )\bxone \bfN}"', from=9-3, to=11-3]
	\arrow[between={0.4}{0.6}, equals, from=11-3, to=7-5]
	\arrow["{\asstwo_{\A,\A,\A}}", from=11-3, to=11-5]
	\arrow["{\what\bxob \mathbbm{1}_{\A}}", from=11-3, to=13-2]
	\arrow["{\mathbbm{1}_{\bfM}\bxtwo \tau^\ell_\bfN}"', between={0.4}{0.6}, Rightarrow, from=11-5, to=7-7]
	\arrow["{\mathbbm{1}_{\A}\bxone \what}", from=11-5, to=13-6]
	\arrow["\alpha", between={0.4}{0.6}, Rightarrow, from=13-2, to=13-6]
	\arrow["\what"', from=13-2, to=15-4]
	\arrow["\what", from=13-6, to=15-4]
	\arrow["{\theta^{r*}\bxtwo \mathbbm{1}_{\mathbbm{1}_\Cc}}"{pos=0.3}, between={0.4}{0.6}, Rightarrow, from=0, to=2]
	\arrow["{\mathbbm{1}_{\mathbbm{1}_\Cc}\bxtwo \theta^\ell}"{pos=0.7}, between={0.4}{0.6}, Rightarrow, from=3, to=1]
\end{tikzcd}
}
\end{align}

In what follows, we will examine the definition of the domain $2$-cell to identify the order of composition of the componenent $2$-cells. We will then use Lemma~\ref{swaptwocells} to show that the order of the strong natural transformations given by whiskerings of $\delta$ and $\mathbbm{1}_{\bfM}\bxtwo \tau^\ell_\bfN$ can be swapped. 
As we will show, the bottom half of the diagram, pre-whiskered with the functor $((\mathbbm{1}_\Cc\bxone\Theta)\bxone \mathbbm{1}_\Cc)(\Delta \bxone \mathbbm{1}_{\Cc})\Delta$, is an identity $2$-cell. What remains is then a diagram using only the comonoidal structure on $\Cc$, precisely the diagram in condition \eqref{middlecond} of the definition of a comonoidal structure, which we assumed to be equal to the identity. This means the domain and codomain of $\muthree_{\bfM,\bfN}$ are both identities, and thus we may choose $\muthree_{\bfM,\bfN}$ to be the identity modification.

We begin by unpacking the definition of the $2$-cell
$(\mathbbm{1}_{\bfM}\tenrepC_1 \lunthree_{\bfN})\, \assthree_{\bfM,\bfT_\Cc, \bfN} \,(\runthree^*_{\bfM}\tenrepC_1\mathbbm{1}_{\bfN} )$.
Below we will use an overline on a $2$-cell from the above diagram to denote that $2$-cell pre and post-whiskered with $\Delta$, $\Theta$, $\what$, the representations $\bfM$, $\bfN$ and $\bfT_\Vv$, and identities, to obtain a strong natural transformation between functors $\Cc$ to $\A$.

Applying the definition of $\tenrepC_1$, we have that \[\runthree^*_{\bfM}\tenrepC_1\mathbbm{1}_{\bfN}=
(\runthree^*_{\bfM}\tenrep_1\mathbbm{1}_{\bfN})\star\Delta = \what \star(\runthree^*_{\bfM}\bxtwo\mathbbm{1}_{\bfN})\star\Delta,\]
where, using Lemma.~\ref{le:whiskerassoc}, we may ignore the brackets. 
Using that \[\runthree^*_{\bfM}=\bigl(\tau_{\bfM}^{r *} \star ((\mathbbm{1}_\Cc \bxone \Theta)\Delta)\bigr)\bigl({\bfM} \star \theta^{r*}\bigr),\] and 
the definition of $\bxtwo$, we have 
\begin{align*}\runthree^*_{\bfM}\bxtwo\mathbbm{1}_{\bfN}  
&= \Bigl(\bigl(\tau_{\bfM}^{r*} \star ((\mathbbm{1}_\Cc \bxone \Theta)\Delta)\bigr)\bxtwo \mathbbm{1}_\bfN \Bigr) \Bigl(\bigl(\bfM \star \theta^{r*}\bigr)\bxtwo \mathbbm{1}_\bfN \Bigr)\\
&= \Bigl(\bigl(\tau_{\bfM}^{r*}  \bxtwo \mathbbm{1}_\bfN \bigr)\star \bigl(((\mathbbm{1}_\Cc\bxone\Theta)\bxone \mathbbm{1}_\Cc)(\Delta \bxone \mathbbm{1}_\Cc)\bigr)\Bigr) 
\Bigl(({\bfM\bxone \bfN}) \star \bigl(\theta^{r*}\bxtwo \mathrm{id}_{\mathbbm{1}_\Cc}\bigr) \Bigr)
\end{align*}
thus giving
\[\runthree^*_{\bfM}\tenrepC_1\mathbbm{1}_{\bfN}=
 {\what} \star
\Bigl(\bigl(\tau_{\bfM}^{r*}  \bxtwo \mathbbm{1}_\bfN \bigr)\star \bigl(((\mathbbm{1}_\Cc\bxone\Theta)\bxone \mathbbm{1}_\Cc)(\Delta \bxone \mathbbm{1}_\Cc)\bigr)\Bigr) 
\Bigl((\bfM\bxone \bfN) \star \bigl(\theta^{r*}\bxtwo \mathrm{id}_{\mathbbm{1}_\Cc}\bigr) \Bigr)
 \star \Delta .\]
Using Lemma~\ref{le:partialinterchange}, we can move ${\Delta}$ and  ${\what}$ inside the brackets to obtain
 \begin{multline*}\runthree^*_{\bfM}\tenrepC_1\mathbbm{1}_{\bfN}=
	\Bigl({\what}\star(\tau^{r*}_\bfM\bxtwo \mathbbm{1}_{\bfN}) \star\bigl(((\mathbbm{1}_\Cc\bxone\Theta)\bxone \mathbbm{1}_\Cc)(\Delta \bxone \mathbbm{1}_\Cc)\Delta\bigr)\Bigr)\\
 \Bigl((\bfM\tenrep_0\bfN) \star \bigl(\theta^{r*}\bxtwo \mathbbm{1}_{\mathbbm{1}_\Cc}\bigr) \star \Delta\Bigr)
.\end{multline*}

Let 
\begin{align*}
\overline{ \theta^{r*}\bxtwo \mathbbm{1}_{\mathbbm{1}_\Cc}}
&=  
(\bfM\tenrep_0\bfN) \star \bigl(\theta^{r*}\bxtwo \mathbbm{1}_{\mathbbm{1}_\Cc}\bigr) \star \Delta, \text{ and}\\
 \overline{\tau^{r*}_\bfM\bxtwo \mathbbm{1}_{\bfN}}
 &={\what}\star(\tau^{r*}_\bfM\bxtwo \mathbbm{1}_{\bfN}) \star\bigl(((\mathbbm{1}_\Cc\bxone\Theta)\bxone \mathbbm{1}_\Cc)(\Delta \bxone \mathbbm{1}_\Cc)\Delta\bigr).
\end{align*}
 
 Then we have that \[\runthree^*_{\bfM}\tenrepC_1\mathbbm{1}_{\bfN} = 
\overline{\tau^{r*}_\bfM\bxtwo \mathbbm{1}_{\bfN}} \;
\overline{\theta^{r*}\bxtwo \mathbbm{1}_{\mathbbm{1}_\Cc}}.\]

Recall the  invertible strong natural transformation $\tilde\mutwo_{\Cc,\Cc}$ from Section \ref{phimdef}.

Let 
\[
\overline{\tilde\mutwo_{\Cc,\Cc}} = \bigl(
	\bfM\tenrep_0\bfN
\bigr)\star \tilde\mutwo_{\Cc,\Cc}\star \bigl(
	((\mathbbm{1}_\Cc\bxone\Theta)\bxone \mathbbm{1}_\Cc)(\Delta\bxone\mathbbm{1}_{\Cc})\Delta
\bigr).
\]

Inserting $\overline{\tilde\mutwo_{\Cc,\Cc}}$ and its inverse into the expression for the domain of $\muthree_{\bfM,\bfN}$, we can write
\begin{align}\label{eq:domainmuthree}
\runthree^*_{\bfM}\tenrepC_1\mathbbm{1}_{\bfN} = 
\overline{\tau^{r*}_\bfM\bxtwo \mathbbm{1}_{\bfN}} \;{\overline{\tilde\mutwo_{\Cc,\Cc}}}^{-1}\;\overline{\tilde\mutwo_{\Cc,\Cc}}\;
\overline{\theta^{r*}\bxtwo \mathbbm{1}_{\mathbbm{1}_\Cc}}.
\end{align}
which is precisely the $2$-morphism represented by the left-hand side of the diagram \eqref{middleunitor2cell}.

Analogously, we have that 
\[\lunthree_{\bfN}=
\left({\bfN} \star \theta^\ell\right) \left(\tau_{\bfN}^\ell \star \bigl((\Theta \bxone \mathbbm{1}_\Cc )\Delta\bigr)\right),\]
and
\begin{align*}
	\mathbbm{1}_{\bfM} \bxtwo \lunthree_{\bfN}  
&= 
\Bigl(\mathbbm{1}_{\bfM}\bxtwo\left({\bfN} \star \theta^\ell\right) \Bigr)
\Bigl(\mathbbm{1}_{\bfM}\bxtwo\left(\tau_{\bfN}^\ell \star \bigl((\Theta \bxone \mathbbm{1}_{\Cc} )\Delta\bigr)\right) \Bigr) \\
&= \Bigl(({\bfM\bxone\bfN })\star\left( \mathrm{id}_{\mathbbm{1}_{\Cc}}\bxtwo\theta^\ell\right) \Bigr)
\Bigl(\left(\mathbbm{1}_{\bfM}\bxtwo\tau_{\bfN}^\ell \right)\star \bigl((\mathbbm{1}_\Cc\bxone(\Theta \bxone \mathbbm{1}_{\Cc} ))(\mathbbm{1}_{\Cc}\bxone\Delta)\bigr) \Bigr)
\end{align*}
and so
\begin{multline*}
	\mathbbm{1}_{\bfM} \tenrepC \lunthree_{\bfN}  
= \Bigl(\bigl(\bfM\tenrep\bfN \bigr)\star\left( \mathrm{id}_{\mathbbm{1}_{\Cc}}\bxtwo \theta^\ell\right)\star \Delta \Bigr)\\
\Bigl({\what}\star\left(\mathbbm{1}_{\bfM}\bxtwo\tau_{\bfN}^\ell \right)\star \bigl((\mathbbm{1}_\Cc\bxone(\Theta \bxone \mathbbm{1}_{\Cc} ))(\mathbbm{1}_{\Cc}\bxone\Delta)\Delta\bigr) \Bigr).
\end{multline*}

Similarly to above, we set
\begin{align*}
\overline{\mathrm{id}_{\mathbbm{1}_\Cc}\bxtwo \theta^\ell}
&=
\bigl(\bfM\tenrep\bfN \bigr)\star\left( \mathrm{id}_{\mathbbm{1}_{\Cc}}\bxtwo\theta^\ell\right)\star \Delta, \text{ and}\\
 \overline{\mathbbm{1}_{\bfM}\bxtwo \tau^\ell_\bfM}
&=
 {\what}\star\left(\mathbbm{1}_{\bfM}\bxtwo\tau_{\bfN}^\ell \right)\star \bigl((\mathbbm{1}_\Cc\bxone(\Theta \bxone \mathbbm{1}_{\Cc} ))(\mathbbm{1}_{\Cc}\bxone\Delta)\Delta\bigr) 
\end{align*}
such that 
\[\mathbbm{1}_{\bfM}\tenrepC_1 \lunthree_{\bfN} =\overline{\mathbbm{1}_{\bfM}\bxtwo \tau^\ell_\bfM} \; \overline{\mathbbm{1}_{\mathbbm{1}_\Cc}\bxtwo \theta^\ell}.\] 

Finally we have that
\[
\assthree_{\bfM,\bfT_\Cc, \bfN} =\Bigl(\bigl((\bfM \tenrep (\bfT_\Cc\tenrep \bfN))\bigr) \star \delta
\Bigr)
\Bigl(
\alpha\star\bigl(((\bfM \bxone \bfT_\Cc)\bxone \bfN)(\Delta\bxone \mathbbm{1}_\Cc)\Delta\bigr)
\Bigr).
\]
Thus, denoting
\begin{align}\label{eq:alphaoverline}
\overline{\alpha}
=
\alpha\star\bigl(((\bfM \bxone \bfT_\Cc)\bxone \bfN)(\Delta\bxone \mathbbm{1}_\Cc)\Delta\bigr)
\end{align} and 
\begin{align}\label{eq:deltaoverline}
\overline{\delta}
=
\bigl((\bfM \tenrep (\bfT_\Cc\tenrep \bfN))\bigr) \star \delta
\end{align} 
we have 
\begin{multline}
	(\mathbbm{1}_{\bfM}\tenrepC_1 \lunthree_{\bfN})\, \assthree_{\bfM,\bfT_\Cc, \bfN} \,(\runthree^*_{\bfM}\tenrepC_1\mathbbm{1}_{\bfN} ) \\
	= \overline{\mathbbm{1}_{\mathbbm{1}_\Cc}\bxtwo \theta^\ell}\; \overline{\mathbbm{1}_{\bfM}\bxtwo \tau^\ell_\bfN}   \; \overline{\delta}\; \overline{\alpha} \;\overline{\tau^{r*}_\bfM\bxtwo \mathbbm{1}_{\bfN}}\; {\overline{\tilde\mutwo_{\Cc,\Cc}}}^{-1}\;\overline{\tilde\mutwo_{\Cc,\Cc}}\;\overline{\theta^{r*}\bxtwo \mathbbm{1}_{\mathbbm{1}_\Cc}} \; .
\end{multline}
Recall that vertical composition of strong natural transformations is associative, so we do not need to add parentheses here.

We now show that 
\begin{multline}\label{eq:CommuteNatTrans}
	 \overline{\mathbbm{1}_{\mathbbm{1}_\Cc}\bxtwo \theta^\ell}\; \overline{\mathbbm{1}_{\bfM}\bxtwo \tau^\ell_\bfN}   \; \overline{\delta}\; \overline{\alpha} \;\overline{\tau^{r*}_\bfM\bxtwo \mathbbm{1}_{\bfN}}\;{\overline{\tilde\mutwo_{\Cc,\Cc}}}^{-1}\;\overline{\tilde\mutwo_{\Cc,\Cc}}\; \overline{\theta^{r*}\bxtwo \mathbbm{1}_{\mathbbm{1}_\Cc}}= \\
	\overline{\mathbbm{1}_{\mathbbm{1}_\Cc}\bxtwo \theta^\ell}\;\overline{\overline{\delta}}\;\overline{\overline{\mathbbm{1}_{\bfM}\bxtwo \tau^\ell_\bfN}}\; \overline{\alpha} \;\overline{\tau^{r*}_\bfM\bxtwo \mathbbm{1}_{\bfN}}\;{\overline{\tilde\mutwo_{\Cc,\Cc}}}^{-1}\;\overline{\tilde\mutwo_{\Cc,\Cc}}\; \overline{\theta^{r*}\bxtwo \mathbbm{1}_{\mathbbm{1}_\Cc}}
\end{multline}
where 
\[\overline{\overline{\delta}}= \bigl( (\bfM\tenrep_0\bfN)(\mathbbm{1}_\Cc\bxone \luntwo_\Cc) (  (\mathbbm{1}_\Cc\bxone \Theta) \bxone\mathbbm{1}_\Cc)\bigr) \star  \delta \]
and 
\[\overline{\overline{\mathbbm{1}_{\bfM}\bxtwo \tau^\ell_\bfN}} =\what \star  (\mathbbm{1}_{\bfM}\bxtwo \tau^\ell_\bfN ) \star   \bigl( (  (\mathbbm{1}_\Cc\bxone \Theta) \bxone\mathbbm{1}_\Cc)\asstwo_{\Cc,\Cc,\Cc}(\Delta\bxone \mathbbm{1}_\Cc)\Delta   \bigr)\]

In other words, we can swap the order of the strong natural transformations whose non-identity components are  ${\mathbbm{1}_{\bfM}\bxtwo \tau^\ell_\bfN}$ and ${\delta}$.

Letting $\alpha = \delta$, $\beta = \what\star(\mathbbm{1}_{\bfM}\bxtwo \tau^\ell_\bfN)$ and $G = \mathbbm{1}_\Cc\bxone(\Theta\bxone \mathbbm{1}_\Cc)$ in Lemma \ref{swaptwocells}, a modification 
\[
\Psi\colon 
\overline{\mathbbm{1}_{\bfM}\bxtwo \tau^\ell_\bfN}   \; \overline{\delta}\to\overline{\overline{\delta}}\;\overline{\overline{\mathbbm{1}_{\bfM}\bxtwo \tau^\ell_\bfN}}\]  is given by the natural transformations, for each $\ob{i}\in \Cc$, 
\begin{align}\label{eq:Modif}
\Psi_{\ob{i}} = \bigl(\what\star(\mathbbm{1}_{\bfM}\bxtwo \tau^\ell_\bfN)\bigr)
_{(\mathbbm{1}_\Cc\bxone(\Theta\bxone \mathbbm{1}_\Cc))(\delta_\ob{i})}.
\end{align}
Looking at Definition~\ref{de:whiskering}, we have that 
\[
\Psi_\ob{i}=\what((\mathbbm{1}_{\bfM}\bxtwo \tau^\ell_\bfN)_{(\mathbbm{1}_\Cc\bxone(\Theta\bxone \mathbbm{1}_\Cc))(\delta_\ob{i})}).
\]
As discussed above, 
$\delta_\ob{i}= (\delta_{\ob{i},1},(\delta_{\ob{i},2},\delta_{\ob{i},3}))$, and
 $\Theta(\delta_{\ob{i},3})=\mathbbm{1}_{\bullet}$ by \eqref{Thetatoone}.
It follows that
\[
\Psi_\ob{i} =
\what((\mathbbm{1}_{\bfM}\bxtwo \tau^\ell_\bfN)_{(\delta_{\ob{i},1},(\mathbbm{1}_{\bullet},\delta_{\ob{i},3}))})
=\mathbbm{1}_{\bfM(\delta_{\ob{i},1})}\otimes_2 \tau^\ell_{(\mathbbm{1}_{\bullet},\delta_{\ob{i},3})}.
\]

Looking at the calulations before \eqref{tauell2}, we see that $(\tau^\ell_\bfN)_{(\mathbbm{1}_{\bullet},\delta_{\ob{i},3})}$ is a natural transformation, with component at an object $(\bullet_\kk^{\oplus n},X)$ given by the canonical resummation from $(\bfM(\delta_{\ob{i},3})X)^{\oplus n}$ to $\bfM(\delta_{\ob{i},3})(X^{\oplus n})$.

Now the modification from the LHS of \eqref{eq:CommuteNatTrans} to the RHS is the composition of $\Psi$ with identity modifications.
These identity modifications have components which are identity natural transformations of the functor which is the $1$-cell component of each strong natural transformation, so for each $\ob{i}\in \Cc$ the component natural transformation of the modification is the composition of natural transformations
\begin{align}\label{eq:MuModifNatTrans}
\mathbbm{1}_{(
\overline{\mathbbm{1}_{\mathbbm{1}_\Cc}\bxtwo \theta^\ell})_{\ob{i}}}\; \Psi_\ob{i}\; \mathbbm{1}_{(\overline{\alpha})_{\ob{i}}}\; \mathbbm{1}_{(\overline{\tau^{r*}_\bfM\bxtwo \mathbbm{1}_{\bfN}})_{\ob{i}}}\; \mathbbm{1}_{(\overline{\theta^{r*}\bxtwo \mathbbm{1}_{\mathbbm{1}_\Cc}})_{\ob{i}}}
\end{align}
This is a natural transformation between endofunctors of $\bfM(\ob{i}_l)\monob \bfN(\ob{i}_r)$,
so consider a general object  $(X,Y)\in \bfM(\ob{i}_l)\monob \bfN(\ob{i}_r)$.
The component of the natural transformation in \eqref{eq:MuModifNatTrans} at $(X,Y)$is 
\[
(\mathbbm{1}_{\Cc}\bxtwo \theta^\ell)_{\ob{i}}((\Psi_\ob{i})_{(\overline{\alpha})_{\ob{i}}(\overline{\tau^{r*}_\bfM\bxtwo \mathbbm{1}_{\bfN}})_{\ob{i}}(\overline{\theta^{r*}\bxtwo \mathbbm{1}_{\mathbbm{1}_\Cc}})_{\ob{i}}(X,Y)}).
\]
Unpacking ${(\overline{\alpha})_{\ob{i}}(\overline{\tau^{r*}_\bfM\bxtwo \mathbbm{1}_{\bfN}})_{\ob{i}}(\overline{\theta^{r*}\bxtwo \mathbbm{1}_{\mathbbm{1}_\Cc}})_{\ob{i}}(X,Y)}$ we have 
\begin{align*}
&= (\alpha)_{((\bfM(\ob{i}_\ell),\Add \kk),\bfN(\ob{i}_r))}
((\tau^{r*}_\bfM)_{(\ob{i}_\ell,\bullet)}\monmor (\mathbbm{1}_{\bfN})_{\ob{i}_r})
(\bfM(\theta^{r*}_{\ob{i}_\ell})\monmor{\bfN(\mathbbm{1}_{\mathbbm{1}_\Cc})}_{\ob{i}_r})(X,Y)\\
&= (\alpha)_{((\bfM(\ob{i}_\ell),\Add \kk),\bfN(\ob{i}_r))}
((\tau^{r*}_\bfM)_{(\ob{i}_\ell,\bullet)}\monmor (\mathbbm{1}_{\bfN})_{\ob{i}_r})
(\bfM(\theta^{r*}_{\ob{i}_\ell})X,Y)\\
& = (\alpha)_{((\bfM(\ob{i}_\ell),\Add \kk),\bfN(\ob{i}_r))}
((\bfM(\theta^{r*}_{\ob{i}_\ell})X,\bullet_\kk),Y)\\
& = (\bfM(\theta^{r*}_{\ob{i}_\ell})X,(\bullet_\kk,Y)).
\end{align*}
Looking at the definition of $\Psi_{\ob{i}}$ in \eqref{eq:Modif}, and the following discussion of the components of $(\tau^\ell_\bfN)_{(\mathbbm{1}_{\bullet},\delta_{\ob{i},3})}$, we see that the component of $\Psi_\ob{i}$ at an object of the form $(\bfM(\theta^{r*}_{\ob{i}_\ell})X,(\bullet_\kk,Y))$ is the identity. Thus, the natural transformation in \eqref{eq:MuModifNatTrans} is indeed an identity, and we have the claimed equality given in \eqref{eq:CommuteNatTrans}.

At this point we have shown \begin{multline}
	(\mathbbm{1}_{\bfM}\tenrepC_1 \lunthree_{\bfN})\, \assthree_{\bfM,\bfT_\Cc, \bfN} \,(\runthree^*_{\bfM}\tenrepC_1\mathbbm{1}_{\bfN} ) \\
	= \overline{\mathbbm{1}_{\mathbbm{1}_\Cc}\bxtwo \theta^\ell}   \; \overline{\overline{\delta}} \; \overline{\overline{\mathbbm{1}_{\bfM}\bxtwo \tau^\ell_\bfN}}\; \overline{\alpha} \;\overline{\tau^{r*}_\bfM\bxtwo \mathbbm{1}_{\bfN}}\;{\overline{\tilde\mutwo_{\Cc,\Cc}}}^{-1}\;\overline{\tilde\mutwo_{\Cc,\Cc}}\; \overline{\theta^{r*}\bxtwo \mathbbm{1}_{\mathbbm{1}_\Cc}} \; .
\end{multline}

We now show that $\overline{\overline{\mathbbm{1}_{\bfM}\bxtwo \tau^\ell_\bfN}}\; \overline{\alpha} \;\overline{\tau^{r*}_\bfM\bxtwo \mathbbm{1}_{\bfN}}\;{\overline{\tilde\mutwo_{\Cc,\Cc}}}^{-1}$ is an identity strong natural transformation.
Initially we ignore the pre-whiskering with $((\mathbbm{1}_\Cc\bxone\Theta)\bxone \mathbbm{1}_\Cc)(\Delta \bxone \mathbbm{1}_\Cc)\Delta$, so we consider the
composition
\begin{align}
\label{eq:HalfMu}
\left(\what\star(\mathbbm{1}_{\bfM}\bxtwo \tau^\ell_\bfN)\star\asstwo_{\Cc,\Vv,\Cc} \right)\left(\alpha \star ((\bfM\bxone \bfT_\Vv)\bxone \bfN)\right)\left(    \what\star (\tau^{r*}_\bfM\bxtwo \mathbbm{1}_{\bfN})\right)\left( (\bfM\tenrep_0\bfN)\star {\tilde\mutwo_{\Cc,\Cc}}^{-1} \right)
\end{align}
which is represented diagramatically by the bottom half of \eqref{middleunitor2cell}.
For an object $((\ob{i},\bullet),\ob{j})\in (\Cc\bxob \Vv)\bxob \Cc$, the $1$-cell component of \eqref{eq:HalfMu} is given by the endofunctor 
\begin{multline*}
	({\mathbbm{1}_{\bfM(\tt{i})}\otimes_1 (\tau^\ell_{\bfM})_{(\bullet,\tt{j})}})(\alpha_{\bfM(\tt{i}),\rm{Add}\mathbbm{k}, \bfN(\tt{j})})((\tau^{r*}_{\bfM})_{(\mathtt{i},\bullet)}\otimes_1 {\mathbbm{1}}_{\bfN(\tt{j})})\mathbbm{1}_{\bfM(\tt{i})\monob\bfN(\tt{j})}\\
	=
	(\mathbbm{1}_{\bfM(\tt{i})}\otimes_1 \ell_{\bfM(\ob{j})})(\alpha_{\bfM(\tt{i}),\rm{Add}\mathbbm{k}, \bfN(\tt{j})})(r_{\bfM(\ob{i})}^*\otimes_1 {\mathbbm{1}}_{\bfN(\tt{j})})
\end{multline*}
of $\bfM(\ob{i})\monob \bfN(\ob{j})$, which we now show is equal to the identity.
Letting $\bigoplus_i(X_i,Y_i)\in \mathbf{M}(\mathtt{i})\otimes_0 \mathbf{N}(\mathtt{j})$ be an object, we have
\begin{multline*}
	(\mathbbm{1}_{\bfM(\tt{i})}\otimes_1 \ell_{\bfM(\ob{j})})(\alpha_{\bfM(\tt{i}),\rm{Add}\mathbbm{k}, \bfN(\tt{j})})(r_{\bfM(\ob{i})}^*\otimes_1 {\mathbbm{1}}_{\bfN(\tt{j})})(\bigoplus_i(X_i,Y_i))\\
	\begin{aligned}
		&=
		(\mathbbm{1}_{\bfM(\tt{i})}\otimes_1 \ell_{\bfM(\ob{j})})(\alpha_{\bfM(\tt{i}),\rm{Add}\mathbbm{k}, \bfN(\tt{j})})(\bigoplus_i((X_i,\bullet_{\Bbbk}),Y_i))\\
		&=(\mathbbm{1}_{\bfM(\tt{i})}\otimes_1 \ell_{\bfM(\ob{j})})(\bigoplus_i(X_i,(\bullet_{\Bbbk},Y_i)))\\
		&=(\bigoplus_i(X_i,Y_i))
	\end{aligned}
\end{multline*}
and similarly for morphisms.

Thus the $1$-cell component of \eqref{eq:HalfMu} is indeed the identity.
It follows that the $1$-cell component of $\overline{\overline{\mathbbm{1}_{\bfM}\bxtwo \tau^\ell_\bfN}}\; \overline{\alpha} \;\overline{\tau^{r*}_\bfM\bxtwo \mathbbm{1}_{\bfN}}\;{\overline{\tilde\mutwo_{\Cc,\Cc}}}^{-1}$ is the identity functor for each object in $\Cc$.

To see that the $2$-cell component is also the identity, note first that the source and target pseudofunctors are equal. Since the $2$-cell data of all component natural transformations is given by canonical isomorphisms (or identities), it follows that $2$-cell component of the composition must be the identity. 

We conclude that
\[
(\mathbbm{1}_{\bfM}\tenrepC_1 \lunthree_{\bfN})\, \assthree_{\bfM,\bfT_\Cc, \bfN} \,(\runthree^*_{\bfM}\tenrepC_1\mathbbm{1}_{\bfN} ) \\
	= \overline{\mathbbm{1}_{\mathbbm{1}_\Cc}\bxtwo \theta^\ell}   \; \overline{\overline{\delta}} \; \overline{\theta^{r*}\bxtwo \mathbbm{1}_{\mathbbm{1}_\Cc}} \; 
\]
which is the identity due to the condition \eqref{middlecond} in the definition of a comonoidal structure.

\subsubsection{The $2$-unitor $\lambdathree$}

We will now consider the domain and codomain of the $2$-unitor $\lambdathree_{\bfM,\bfN}$. Similarly to the above approach for $\mu_{\bfM,\bfN}$ we will show that we can reorder the composition of strong natural transformations in the codomain of $\lambdathree_{\bfM,\bfN}$, which will allow us to consider the part of the composition containing all instances of $\bfM$ and $\bfN$ and show that these parts of the domain and codomain are equal. It follows that, given condition \ref{??}, the domain and codomain of $\lambdathree_{\bfM,\bfN}$ are equal, and we may choose it to be the identity.

Both the domain and codomain of $\lambdathree_{\bfM,\bfN}$ may be represented, non uniquely as discussed above, by the following diagram. The domain corresponds to the left hand side of the diagram, which is a strong natural transformation from the middle blue line to the left-most composition of $1$-cells, and the codomain corresponds to the right hand side of the diagram, which is a strong natural transformation from the middle blue line to the right-most composition of $1$-cells, where the left-most and right-most composition of $1$-cells are equal.

\begin{align}\label{lambdaunitor2cell}\scalebox{0.7}{
\begin{tikzcd}[ampersand replacement=\&]
	\&\&\& \Cc \&\&\&\& \\
	\&\& {\Cc\bxone \Cc} \&\& {\Cc\bxob \Cc} \\
	\\
	\&\& {(\Cc\bxob \Cc )\bxob \Cc} \&\& {\Cc\bxob(\Cc\bxob \Cc )} \& {\Vv\bxob \Cc} \&\& \Cc \\
	\\
	{\Cc\bxob \Cc} \&\& {(\Vv\bxob \Cc )\bxob \Cc} \&\& {\Vv\bxob(\Cc\bxob \Cc )} \&\&\& {\Cc\bxob \Cc} \\
	{\Cc\bxob \Cc} \& {\Vv\bxob(\Cc\bxob \Cc )} \& {(\Vv\bxob \Cc )\bxob \Cc} \\
	{\Cc\bxob \Cc} \&\& {(\Vv\bxob \Cc )\bxob \Cc} \\
	\&\& {(\A\bxob  \A)\bxob \A} \&\& {\A\bxob(\A\bxob \A)} \\
	\\
	\\
	\&\& {\A\bxob \A} \&\& {\A\bxob\A} \&\&\& {\A\bxob \A} \\
	\&\&\& \A
	\arrow["\Delta"', color={rgb,255:red,92;green,92;blue,214}, from=1-4, to=2-3]
	\arrow["\Delta", from=1-4, to=2-5]
	\arrow[""{name=0, anchor=center, inner sep=0}, curve={height=-30pt}, equals, from=1-4, to=4-8]
	\arrow["\delta", between={0.4}{0.6}, Rightarrow, from=2-3, to=2-5]
	\arrow["{\Delta\bxone\mathbbm{1}_{\Cc}}", color={rgb,255:red,92;green,92;blue,214}, from=2-3, to=4-3]
	\arrow[""{name=1, anchor=center, inner sep=0}, curve={height=12pt}, equals, from=2-3, to=6-1]
	\arrow["{\mathbbm{1}_\Cc\bxone\Delta}"', from=2-5, to=4-5]
	\arrow[""{name=2, anchor=center, inner sep=0}, "{\Theta \bxone \mathbbm{1}_{\Cc}}", from=2-5, to=4-6]
	\arrow["{\asstwo_{\Cc,\Cc,\Cc}}", from=4-3, to=4-5]
	\arrow[""{name=3, anchor=center, inner sep=0}, "{(\Theta\bxone\mathbbm{1}_\Cc)\bxone \mathbbm{1}_\Cc}"', color={rgb,255:red,92;green,92;blue,214}, from=4-3, to=6-3]
	\arrow[between={0.4}{0.6}, equals, from=4-5, to=4-6]
	\arrow["{\Theta\bxone (\mathbbm{1}_\Cc\bxone\mathbbm{1}_\Cc)}"', from=4-5, to=6-5]
	\arrow["{\luntwo_\Cc}"', from=4-6, to=4-8]
	\arrow["{\mathbbm{1}_\Vv \bxone \Delta}", from=4-6, to=6-5]
	\arrow["\Delta", from=4-8, to=6-8]
	\arrow[equals, from=6-1, to=7-1]
	\arrow[between={0.4}{0.6}, equals, from=6-3, to=4-5]
	\arrow["{\luntwo_\Cc\bxone \mathbbm{1}_\Cc}"', from=6-3, to=6-1]
	\arrow["{\asstwo_{\Vv,\Cc,\Cc}}", from=6-3, to=6-5]
	\arrow[draw={rgb,255:red,92;green,92;blue,214}, equals, from=6-3, to=7-3]
	\arrow["{\luntwo_\Delta^{-1}}", between={0.4}{0.6}, Rightarrow, from=6-5, to=4-8]
	\arrow["{\luntwo_{\Cc \bxob \Cc}}"', from=6-5, to=6-8]
	\arrow["{\bfT_\Vv \bxone(\bfM\bxone\bfN)}", from=6-5, to=9-5]
	\arrow[""{name=4, anchor=center, inner sep=0}, "{\bfM\bxone\bfN}", from=6-8, to=12-8]
	\arrow[equals, from=7-1, to=8-1]
	\arrow["{\luntwo_{\Cc \bxob \Cc}}"', dashed, from=7-2, to=7-1]
	\arrow["{\lambdatwo_{\Cc,\Cc}^{-1}}"', between={0.4}{0.6}, Rightarrow, dashed, from=7-3, to=6-1]
	\arrow["{\asstwo_{\Vv,\Cc,\Cc}}"', dashed, from=7-3, to=7-2]
	\arrow[draw={rgb,255:red,92;green,92;blue,214}, equals, from=7-3, to=8-3]
	\arrow["{\tau^{\ell}_{\bfM}\bxtwo \mathbbm{1}_{\bfN}}"'{pos=0.7}, between={0.4}{0.6}, equals, 2tail reversed, from=8-1, to=9-3]
	\arrow["{\bfM\bxone\bfN}"', curve={height=12pt}, from=8-1, to=12-3]
	\arrow["{\lambdatwo_{\Cc,\Cc}}", between={0.4}{0.6}, Rightarrow, dashed, from=8-3, to=7-1]
	\arrow["{\luntwo_\Cc\bxone \mathbbm{1}_\Cc}", from=8-3, to=8-1]
	\arrow["{(\bfT_\Vv\bxone\bfM )\bxone \bfN}"', color={rgb,255:red,92;green,92;blue,214}, from=8-3, to=9-3]
	\arrow[between={0.4}{0.6}, equals, from=9-3, to=6-5]
	\arrow["{\asstwo_{\A,\A,\A}}", from=9-3, to=9-5]
	\arrow["{\what\bxone \mathbbm{1}_{\A}}", color={rgb,255:red,92;green,92;blue,214}, from=9-3, to=12-3]
	\arrow[""{name=5, anchor=center, inner sep=0}, "{\mathbbm{1}_{\A}\bxone \what}", from=9-5, to=12-5]
	\arrow["\alpha", between={0.4}{0.6}, Rightarrow, from=12-3, to=9-5]
	\arrow["\what"', color={rgb,255:red,92;green,92;blue,214}, from=12-3, to=13-4]
	\arrow["\what", from=12-5, to=13-4]
	\arrow["\what", curve={height=-12pt}, from=12-8, to=13-4]
	\arrow["{\theta^{\ell}\bxtwo \mathbbm{1}_{\mathbbm{1}_{\Cc}}}"{pos=0.3}, between={0.4}{0.6}, equals, 2tail reversed, from=1, to=3]
	\arrow["{\theta^{\ell}}", between={0.4}{0.6}, Rightarrow, from=2, to=0]
	\arrow["{\tau^\ell_{\bfM\tenrep_0\bfN}}", between={0.4}{0.6}, Rightarrow, from=5, to=4]
\end{tikzcd}
}
\end{align}

We being by unpacking the definition of the domain $\lunthree_{\bfM}\tenrepC_1 \mathbbm{1}_{\bfN}$ as above.

Applying the definition of $\tenrepC_1$, we have 
\[
\lunthree_{\bfM}\tenrepC_1 \mathbbm{1}_{\bfN}
=(\lunthree_{\bfM}\tenrep_1 \mathbbm{1}_{\bfN})\star\Delta
=\what\star(\lunthree_{\bfM}\bxtwo\mathbbm{1}_{\bfN})\star\Delta.\]

Using that \[\lunthree_{\bfM}=
\left(
	\bfM \star \theta^\ell
\right)\left(\tau^\ell_\bfM \star ((\Theta\bxone\mathbbm{1}_\Cc  )\Delta)
\right)
\]
we have 
\begin{align*}
\lunthree_{\bfM}\bxtwo\mathbbm{1}_{\bfN}
&=
\Bigl(
	\bigl(
		\bfM \star \theta^\ell
	\bigr)\bxtwo \mathbbm{1}_{\bfN}
\Bigr)
\Bigl(
	\bigl(\tau^\ell_\bfM \star ((\Theta\bxone\mathbbm{1}_\Cc  )\Delta) 
	\bigr) \bxtwo \mathbbm{1}_{\bfN}
\Bigr)\\
&=
\Bigl(
	({\bfM\bxone \bfN}) \star \bigl(\theta^\ell\bxtwo \mathbbm{1}_{\mathbbm{1}_\Cc}\bigr)
\Bigr)
\Bigl(
	\bigl(
		\tau^{\ell}_\bfM \bxtwo \mathbbm{1}_\bfN
	\bigr)\star 
	\bigl(((\Theta\bxone\mathbbm{1}_\Cc)\bxone \mathbbm{1}_\Cc)(\Delta \bxone \mathbbm{1}_\Cc)\bigr)
\Bigr)
\end{align*}
thus obtaining 
\[
\lunthree_{\bfM}\tenrepC_1 \mathbbm{1}_{\bfN}=
\what\star
\Bigl(
	({\bfM\bxone \bfN}) \star \bigl(\theta^\ell\bxtwo \mathbbm{1}_{\mathbbm{1}_\Cc}\bigr)
\Bigr)
\Bigl(
	\bigl(
		\tau^{\ell}_\bfM \bxtwo \mathbbm{1}_\bfN
	\bigr)\star 
	\bigl(((\Theta\bxone\mathbbm{1}_\Cc)\bxone \mathbbm{1}_\Cc)(\Delta \bxone \mathbbm{1}_\Cc)\bigr)
\Bigr)
\star\Delta.
\]

Again, using Lemma~\ref{le:partialinterchange}, we can move ${\Delta}$ and  ${\what}$ inside the brackets to obtain
 \begin{multline*}
\lunthree_{\bfM}\tenrepC_1 \mathbbm{1}_{\bfN}=
\Bigl(
	(\bfM\tenrep_0 \bfN) \star \bigl(\theta^\ell\bxtwo \mathbbm{1}_{\mathbbm{1}_\Cc}\bigr)\star \Delta
\Bigr)\\
\Bigl(
	\what \star
	\bigl(
		\tau^{\ell}_\bfM \bxtwo \mathbbm{1}_\bfN
	\bigr)\star 
	\bigl(((\Theta\bxone\mathbbm{1}_\Cc)\bxone \mathbbm{1}_\Cc)(\Delta \bxone \mathbbm{1}_\Cc)
	\Delta
	\bigr)
\Bigr)
.\end{multline*}

Let 
\begin{align*}
\overline{\theta^\ell\bxtwo \mathbbm{1}_{\mathbbm{1}_\Cc}}
&=  
(\bfM\tenrep_0 \bfN) \star \bigl(\theta^\ell\bxtwo \mathbbm{1}_{\mathbbm{1}_\Cc}\bigr)\star \Delta, \text{ and}\\
 \overline{\tau^\ell\bxtwo \mathbbm{1}_{\bfN}}
 &=\what \star
	\bigl(
		\tau^{\ell}_\bfM \bxtwo \mathbbm{1}_\bfN
	\bigr)\star 
	\bigl(((\Theta\bxone\mathbbm{1}_\Cc)\bxone \mathbbm{1}_\Cc)(\Delta \bxone \mathbbm{1}_\Cc)
	\Delta
	\bigr).
\end{align*}
 
 Then we have \[\lunthree_{\bfM}\tenrepC_1\mathbbm{1}_{\bfN} = 
\overline{\theta^\ell\bxtwo \mathbbm{1}_{\mathbbm{1}_\Cc}}\;
\overline{\tau^\ell\bxtwo \mathbbm{1}_{\bfN}} .\]

Recall the definition of $\lambdatwo_{\Cc,\Cc}$ from Proposition \ref{2unitortwo}.
Let 
\[
\overline{\lambdatwo_{\Cc,\Cc}} = \bigl(
	\bfM\tenrep_0\bfN
\bigr)\star \lambdatwo_{\Cc,\Cc} \star \bigl(
	((\Theta\bxone\mathbbm{1}_\Cc)\bxone \mathbbm{1}_\Cc)(\Delta\bxone\mathbbm{1}_{\Cc})\Delta
\bigr)
\]

Inserting $\overline{\lambdatwo_{\Cc,\Cc}}$ and its inverse into the expression for the domain of $\lambdathree_{\bfM,\bfN}$, we can write
\begin{align}\label{eq:domainlambdathree}
\lunthree_{\bfM}\tenrepC_1\mathbbm{1}_{\bfN} = 
\overline{\theta^{\ell}\bxtwo \mathbbm{1}_{\mathbbm{1}_\Cc}}\;
 \overline{\lambdatwo_{\Cc,\Cc}}^{-1}\overline{\lambdatwo_{\Cc,\Cc}}\;
 \overline{\tau^\ell\bxtwo \mathbbm{1}_{\bfN}}
\end{align}
which is precisely the $2$-morphism represented by the left-hand side of the diagram \eqref{lambdaunitor2cell}.

Now considering the codomain of $\lambdathree_{\bfM,\bfN}$, note first that
\begin{equation}\label{eq:elltenrep}\lunthree_{\bfM\tenrepC_0\bfN}
=\bigl((\bfM \tenrepC_0\bfN)\star \theta^\ell\bigr) \bigl(\tau^{\ell}_{\bfM \tenrepC_0\bfN}\star((\Theta\bxone\mathbbm{1}_\Cc)\Delta)\bigr).
\end{equation}

It will be useful to split $\tau^{r}_{\bfM \tenrepC_0\bfN}$ into two parts by inserting the (inverse of the) strong transformation $\luntwo_\Delta$ from Proposition \ref{unitortwocommute}.
This is illustrated by the following diagram, where we will prove the equality below. 
\begin{multline}
\begin{tikzcd}[ampersand replacement=\&]
	{\Vv \bxob \Cc} \&\& {\A\bxob \A} \\
	\\
	\Cc \&\& \A
	\arrow["{\bfT_\Vv \bxone ({\bfM\tenrepC_0\bfN})}", from=1-1, to=1-3]
	\arrow["{\luntwo_\Cc}"', from=1-1, to=3-1]
	\arrow["{\tau^\ell_{\bfM\tenrepC_0\bfN}}", between={0.4}{0.6}, Rightarrow, from=1-3, to=3-1]
	\arrow["\what", from=1-3, to=3-3]
	\arrow["{{\bfM\tenrepC_0\bfN}}"', from=3-1, to=3-3]
\end{tikzcd}
=\\
\begin{tikzcd}[ampersand replacement=\&]
	{\Vv \bxob \Cc} \&\& {\Vv\bxob(\Cc\bxone\Cc)} \&\& {\A\bxob \A} \\
	\\
	\Cc \&\& {\Cc\bxone\Cc} \&\& \A
	\arrow["{\mathbbm{1}_\Vv\bxone \Delta}", from=1-1, to=1-3]
	\arrow["{\luntwo_\Cc}"', from=1-1, to=3-1]
	\arrow["{\bfT_\Vv \bxone ({\bfM\tenrep_0\bfN})}", from=1-3, to=1-5]
	\arrow["{\luntwo_\Delta^{-1}}"', between={0.4}{0.6}, Rightarrow, from=1-3, to=3-1]
	\arrow["{\luntwo_{\Cc\bxob\Cc}}"', from=1-3, to=3-3]
	\arrow["{\tau^\ell_{\bfM\tenrep_0\bfN}}"', between={0.4}{0.6}, Rightarrow, from=1-5, to=3-3]
	\arrow["\what", from=1-5, to=3-5]
	\arrow["\Delta"', from=3-1, to=3-3]
	\arrow["{\bfM\tenrep_0\bfN}"', from=3-3, to=3-5]
\end{tikzcd}
\end{multline}

\begin{lemma}
We have 
\[\tau^l_{\bfM\tenrepC_0\bfN} =  
\left((\bfM\tenrep_0\bfN)\star \luntwo_\Delta^{-1}\right)
\left(\tau^\ell_{\bfM\tenrep_0\bfN} \star (\mathbbm{1}_{\Vv}\bxone\Delta)\right)
\] as strong transformations from $(\bfM\tenrepC_0\bfN)\tenrep_0 \bfT_\Vv$ to $(\bfM\tenrepC_0\bfN)\luntwo_\Cc$.
\end{lemma}

\proof
The $1$-cell components at an object $(\bullet, \ob{i})$ in both cases given by $\lun_{\bfM(\ob{i}_\ell) \monob \bfN(\ob{i}_{r})}$. Thus the source and target of the $2$-cell components of both strong transformations are the same, and since both are given by canonical resummations and additivity data of the functors involved, they coincide.
\endproof

Let
\[
\overline{\luntwo_\Delta^{-1}} = \left(\bfM\tenrep_0\bfN\right)\star \luntwo_\Delta^{-1} \star \left((\Theta\bxone\mathbbm{1}_\Cc)\Delta\right) ,
\]
\[
\overline{\tau^\ell_{\bfM\tenrep_0\bfN}} = 
\tau^\ell_{\bfM\tenrep_0\bfN} \star (\mathbbm{1}_{\Vv}\bxone\Delta)((\Theta\bxone\mathbbm{1}_\Cc)\Delta)
\]
and
\[
\overline{\theta^{\ell}} = 
\bigl((\bfM \tenrepC_0\bfN)\star \theta^\ell\bigr).
\]
Thus proceeding from \eqref{eq:elltenrep} we have that
\begin{align}\label{eq:elltenrep2}
	\lunthree_{\bfM\tenrepC_0\bfN} = \overline{\theta^\ell}\,\overline{\luntwo_\Delta^{-1}} \, \overline{\tau^\ell_{\bfM\tenrep_0\bfN}},
\end{align}
and hence the codomain of $\lambdathree_{\bfM,\bfN}$ is
\begin{align}
\lunthree_{\bfM\tenrepC_0\bfN}\, \assthree_{\bfT_\Cc,\bfM,\bfN }  = 
\overline{\theta^\ell}\, \overline{\luntwo_\Delta^{-1}} \, \overline{\tau^\ell_{\bfM\tenrep_0\bfN}} \, \overline{\delta} \, \overline{\alpha}
\end{align}
where 
\begin{align}\label{eq:deltaoverline2}
\overline{\delta} = \bigl(\bfT_\Cc \tenrep_0 (\bfM\tenrep_0 \bfN)\bigr) \star \delta
\end{align} 
 and 
\begin{align}\label{eq:alphaoverline2}
\overline{\alpha} = \alpha\star\bigl((( \bfT_\Cc \bxone\bfM)\bxone \bfN)(\Delta\bxone \mathbbm{1}_\Cc)\Delta\bigr)
\end{align} 
in analogy to \eqref{eq:alphaoverline} and \eqref{eq:deltaoverline}, but now with the roles of $\bfT_\Cc$ and $\bfM$ swapped in the whiskerings.

We now claim that we can swap the strong natural transformations involving $\delta$ and $\tau^\ell_{\bfM\tenrep_0\bfN}$ on the right hand side of the desired equality. More precisely, we claim that
\begin{equation}\label{eq:deltatauswapinlambda}
\overline{\tau^\ell_{\bfM\tenrep_0\bfN}} \, \overline{\delta} =  \overline{\overline{\delta}}\,\overline{\overline{\tau^\ell_{\bfM\tenrep_0\bfN}}}
\end{equation}
where 
\[
\overline{\overline{\delta}} = \bigl(  (\bfM\tenrep_0\bfN) \luntwo_{\Cc\bxob\Cc} (\Theta\bxone(\mathbbm{1}_\Cc\bxone \mathbbm{1}_\Cc)) \bigr)\star \delta
\]
and
\begin{equation*}\begin{split}
\overline{\overline{\tau^\ell_{\bfM\tenrep_0\bfN}}} &= \tau^\ell_{\bfM\tenrep_0\bfN} \star\bigl( (\Theta\bxone(\mathbbm{1}_\Cc\bxone \mathbbm{1}_\Cc))  \asstwo_{\Cc,\Cc,\Cc} (\Delta \bxone \mathbbm{1}_\Cc) \Delta  \bigr)\\
&=\tau^\ell_{\bfM\tenrep_0\bfN} \star\bigl(  \asstwo_{\Vv,\Cc,\Cc} (\Theta\bxone(\mathbbm{1}_\Cc\bxone \mathbbm{1}_\Cc)) (\Delta \bxone \mathbbm{1}_\Cc) \Delta  \bigr).
\end{split}\end{equation*}
Indeed, by Lemma \ref{swaptwocells}, the invertible modification between the two sides of \eqref{eq:deltatauswapinlambda}, has $2$-cell components given by the natural transformations
\[
(\tau^\ell_{\bfM\tenrep_0\bfN})_{ (\Theta\bxone(\mathbbm{1}_\Cc\bxone \mathbbm{1}_\Cc))  (\delta_{\ob{i}})} = (\tau^\ell_{\bfM\tenrep_0\bfN})_{(\mathbbm{1}_\bullet, (\delta_{\ob{i},2},\delta_{\ob{i},3})) }.
\]
To check that these are the identity natural transformations, we evaluate them on $\bigoplus_s( \bullet_\Bbbk^{\oplus a_s}, \bigoplus_t (X_{st}, Y_{st}))$ in $\Add \Bbbk \monob  ( \bfM(\ob{i}_{\ell r})\monob\bfN(\ob{i}_r))$. By \eqref{tauell2}, this is the canonical resummation from

\begin{equation*}\begin{split}
\bigl(\lun_{\bfM(\ob{i}_{r \ell})\monob\bfN(\ob{i}_{rr})}&((\bfT_\Vv\tenrep_0(\bfM\tenrep_0\bfN))(\mathbbm{1}_\bullet, (\delta_{\ob{i},2},\delta_{\ob{i},3}) ))\bigr)\left(\bigoplus_s( \bullet_\Bbbk^{\oplus a_s}, \bigoplus_t (X_{st}, Y_{st}))\right)\\
&= \lun_{\bfM(\ob{i}_{r \ell})\monob\bfN(\ob{i}_{rr})}
\bigl(\bigoplus_s( \bullet_\Bbbk^{\oplus a_s}, \bigoplus_t (  \bfM(\delta_{\ob{i},2})X_{st},\bfN(\delta_{\ob{i},3}) Y_{st}))
\bigr) 
\\ 
&=\bigoplus_s(\bigoplus_t (  \bfM(\delta_{\ob{i},2})X_{st},\bfN(\delta_{\ob{i},3}) Y_{st}))^{\oplus a_s})
\end{split}\end{equation*}
to
\begin{equation*}\begin{split}
\bigl( (\bfM\tenrep_0\bfN)(\luntwo_\Cc (\mathbbm{1}_\bullet, (\delta_{\ob{i},2},\delta_{\ob{i},3}) ))&
\lun_{(\bfM\tenrep_0\bfN)(\ob{i}_{\ell r} , \ob{i}_{r}  )}\bigr) (\bigoplus_s( \bullet_\Bbbk^{\oplus a_s}, \bigoplus_t (X_{st}, Y_{st})))\\
& = (\bfM(\delta_{\ob{i},2}) \bxone\bfN(\delta_{\ob{i},3})) ( \bigoplus_s( \bigoplus_t (X_{st}, Y_{st})^{\oplus a_s}) )\\
&=( \bigoplus_s( \bigoplus_t (\bfM(\delta_{\ob{i},2})X_{st}, \bfN(\delta_{\ob{i},3})Y_{st})^{\oplus a_s}) )
\end{split}\end{equation*}
Given those expressions are equal, the isomorphism is indeed the identity and the equality in \eqref{eq:deltatauswapinlambda} holds.

Thus, the codomain of $\lambdathree_{\bfM,\bfN}$ is
\begin{align}
\lunthree_{\bfM\tenrepC_0\bfN}\, \assthree_{\bfT_\Cc,\bfM,\bfN }  = 
\overline{\theta^\ell}\, \overline{\luntwo_\Delta^{-1}} \, \overline{\overline{\delta}}\,\overline{\overline{\tau^\ell_{\bfM\tenrep_0\bfN}}} \, \overline{\alpha}.
\end{align}

To summarise, we would like to show that

\[\overline{\theta^{\ell}\bxtwo \mathbbm{1}_{\mathbbm{1}_\Cc}}\;
 \overline{\lambdatwo_{\Cc,\Cc}}^{-1}\overline{\lambdatwo_{\Cc,\Cc}}\;
 \overline{\tau^\ell\bxtwo \mathbbm{1}_{\bfN}} =  \overline{\theta^\ell}\, \overline{\luntwo_\Delta^{-1}} \,\overline{\overline{\delta}}\,\overline{\overline{\tau^\ell_{\bfM\tenrep_0\bfN}}}\, \overline{\alpha}
.\]

We now show that 
\[
\overline{\overline{\tau^\ell_{\bfM\tenrep_0\bfN}}}\overline{\alpha} = \overline{\lambdatwo_{\Cc,\Cc}}\;
 \overline{\tau^\ell\bxtwo \mathbbm{1}_{\bfN}}.
\]

In fact, this follows by pre-whiskering with  $\bigl( (\Theta\bxone(\mathbbm{1}_\Cc\bxone \mathbbm{1}_\Cc))   (\Delta \bxone \mathbbm{1}_\Cc) \Delta  \bigr)$ from the following lemma.

\begin{lemma}
We have
\[\tau^\ell_{\bfM\tenrep_0\bfN} \star\asstwo_{\Vv,\Cc,\Cc}= \what\star (\lambdatwo_{\Cc,\Cc}(\tau^\ell\bxtwo \mathbbm{1}_{\bfN}))\]
as strong natural transformations 
\[
(\bfT_\Vv\tenrep_0 \bfM)\tenrep_0 \bfN) \to  (\bfM\tenrep_0 \bfN) \luntwo_{\Cc \bxob \Cc}\asstwo_{\Cc,\Cc,\Vv}  
\]
of pseudofunctors from $(\Vv\bxob \Cc)\bxob\Cc$ to $\A$.
\end{lemma}
\proof
In order to see this, we compute the $1$-cells associated to the two strong natural transformations, which are given by functors in $\A$. Explicitly, for
$\bigoplus_t(\bigoplus_s(\bullet^{\oplus a_{ts}},X_{ts}), Y_t) \in (\Add\, \Bbbk\monob \bfM(\ob{i}))\monob \bfN(\ob{j})$, we compute
\begin{equation*}
\begin{split}
(\tau^\ell_{\bfM\tenrep_0\bfN}& \star\asstwo_{\Vv,\Cc,\Cc}))_{((\bullet ,\ob{i}),\ob{j})}\bigl(\bigoplus_t(\bigoplus_s(\bullet^{\oplus a_{ts}},X_{ts}), Y_t)\bigr)\\
&=(\tau^\ell_{\bfM\tenrep_0\bfN})_{(\bullet ,(\ob{i},\ob{j}))}\bigl(\bigoplus_t\bigoplus_s(\bullet^{\oplus a_{ts}},(X_{ts}, Y_t))\bigr)\\
&=\bigoplus_t\bigoplus_s(X_{ts}, Y_t)^{\oplus a_{ts}}
\end{split}
\end{equation*}
while 
\begin{equation*}
\begin{split}
\bigl( (\what\star &\lambdatwo_{\Cc,\Cc})(\what\star (\tau^\ell\bxtwo \mathbbm{1}_{\bfN}))\bigr)_{((\bullet ,\ob{i}),\ob{j})}\bigl(\bigoplus_t(\bigoplus_s(\bullet^{\oplus a_{ts}},X_{ts}), Y_t)\bigr)\\
 &=(\what\star \lambdatwo_{\Cc,\Cc})_{((\bullet ,\ob{i}),\ob{j})} \bigl(\bigoplus_t(\bigoplus_s(X_{ts}^{\oplus a_{ts}}), Y_t)\bigr)\\
 &=\bigoplus_t\bigoplus_s(X_{ts}, Y_t)^{\oplus a_{ts}}
\end{split}
\end{equation*}
so the functors agree on objects and it is straightforward to see that they agree on morphisms as well. Hence the $1$-cell components of our strong natural transformations coincide. Thus the $2$-cell components of both have the same source and target, and since both are given by resummations, they coincide.
\endproof

It follows that $\lunthree_{\bfM}\tenrepC_1\mathbbm{1}_{\bfN}=\lunthree_{\bfM\tenrepC_0\bfN}\, \assthree_{\bfT_\Cc,\bfM,\bfN }$
if 
\[
\overline{\theta^{\ell}\bxtwo \mathbbm{1}_{\mathbbm{1}_\Cc}}\;
 \overline{\lambdatwo_{\Cc,\Cc}}^{-1}
=
\overline{\theta^\ell}\, \overline{\luntwo_\Delta^{-1}} \, \overline{\delta}
\]
which is precisely the post-whiskering with $\bfM\tenrep_0\bfN$ of the equality in condition \eqref{leftcond} in the definition of a comonoidal structure.

\subsubsection{The $2$-unitor $\rhothree$}

In this subsection, we will prove that $\rhothree$ can be chosen as the identity. In order to do this, we need to verify that for each pair of strictly unital birepresentations $\bfM$ and $\bfN$ of $\Cc$, the domain $\mathbbm{1}_{\bfM}\tenrepC_1 \runthree^*_{\bfN}$ and the codomain $ \assthree_{\bfM,\bfN,\bfT_\Cc } \,\runthree^*_{\bfM\tenrepC_0\bfN}$ of $\rhothree_{\bfM,\bfN}$ coincide. Each of these morphisms of birepresentations is a strong transformation with domain $\bfM\tenrepC_0\bfN$ and codomain $\bfM\tenrepC_0(\bfN\tenrepC_0 \bfT_\Cc)$.

Both may be represented, non-uniquely as discussed above, by the following diagram, where the left-most and right-most composition of $1$-cells are equal. The domain $\mathbbm{1}_{\bfM}\tenrepC_1 \runthree^*_{\bfN}$ corresponds to the right hand side of the diagram, which is a strong natural transformation from the right-most composition of $1$-cells to the middle blue line, and the codomain $ \assthree_{\bfM,\bfN,\bfT_\Cc } \,\runthree^*_{\bfM\tenrepC_0\bfN}$ corresponds to the left hand side of the diagram, which is a strong natural transformation from the right-most composition of $1$-cells to the middle blue line.

\begin{align}\label{rhounitor2cell}\scalebox{0.7}{
\begin{tikzcd}[ampersand replacement=\&]
	\&\&\& \Cc \&\&\& \\
	\&\& {\Cc\bxob\Cc} \&\& {\Cc\bxob\Cc} \\
	\\
	\Cc \& {\Cc\bxob\Vv} \& {(\Cc\bxob\Cc)\bxob\Cc} \&\& {\Cc\bxob(\Cc\bxob\Cc)} \\
	\\
	{\Cc\bxob\Cc} \&\& {(\Cc\bxob\Cc)\bxob\Vv} \&\& {\Cc\bxob(\Cc\bxob\Vv)} \&\& {\Cc\bxob\Cc} \\
	{\Cc\bxob\Cc} \& {\Cc\bxob(\Cc\bxob\Vv)} \& {(\Cc\bxob\Cc)\bxob\Vv} \\
	{\Cc\bxob\Cc} \&\& {(\Cc\bxob\Cc)\bxob\Vv} \\
	\\
	\&\& {(\A\bxob \A)\bxob \A} \&\& {\A\bxob(\A\bxob \A)} \\
	\\
	{\A\bxob\A} \&\& {\A \bxob \A} \&\& {\A \bxob \A} \\
	\&\&\& \A
	\arrow["\Delta"', from=1-4, to=2-3]
	\arrow["\Delta", color={rgb,255:red,92;green,92;blue,214}, from=1-4, to=2-5]
	\arrow[""{name=0, anchor=center, inner sep=0}, curve={height=30pt}, equals, from=1-4, to=4-1]
	\arrow[""{name=1, anchor=center, inner sep=0}, "{\mathbbm{1}_\Cc\bxone \Theta}", from=2-3, to=4-2]
	\arrow[""{name=2, anchor=center, inner sep=0}, "{\Delta\bxone \mathbbm{1}_{\Cc}}", from=2-3, to=4-3]
	\arrow["{\mathbbm{1}_{\Cc}\bxone \Delta}"', color={rgb,255:red,92;green,92;blue,214}, from=2-5, to=4-5]
	\arrow[""{name=3, anchor=center, inner sep=0}, curve={height=-18pt}, equals, from=2-5, to=6-7]
	\arrow["\Delta"', from=4-1, to=6-1]
	\arrow["{\runtwo_\Delta}", between={0.4}{0.6}, Rightarrow, from=4-1, to=6-3]
	\arrow["{\runtwo_\Cc}", from=4-2, to=4-1]
	\arrow[between={0.4}{0.6}, equals, from=4-2, to=4-3]
	\arrow["{\Delta\bxone \mathbbm{1}_{\Vv}}", from=4-2, to=6-3]
	\arrow["{\asstwo_{\Cc,\Cc,\Cc}}"', from=4-3, to=4-5]
	\arrow["{(\mathbbm{1}_{\Cc}\bxone\mathbbm{1}_{\Cc} )\bxone \Theta}", from=4-3, to=6-3]
	\arrow["{\mathbbm{1}_\Cc\bxone(\mathbbm{1}_\Cc\bxone \Theta)}"', color={rgb,255:red,92;green,92;blue,214}, from=4-5, to=6-5]
	\arrow[equals, from=6-1, to=7-1]
	\arrow[between={0.4}{0.6}, equals, from=6-3, to=4-5]
	\arrow["{\runtwo_{\Cc\bxob\Cc}}"', from=6-3, to=6-1]
	\arrow["{\asstwo_{\Cc,\Cc,\Vv}}"', from=6-3, to=6-5]
	\arrow[equals, from=6-3, to=7-3]
	\arrow["{\mathbbm{1}_\Cc\bxone \runtwo_\Cc}"', from=6-5, to=6-7]
	\arrow["{\bfM \bxone (\bfN \bxone \bfT_{\Vv})}"', color={rgb,255:red,92;green,92;blue,214}, from=6-5, to=10-5]
	\arrow["{\mathbbm{1}_\bfM \bxone {\tau^{r*}_\bfN}}"', between={0.4}{0.6}, Rightarrow, from=6-7, to=10-5]
	\arrow["{\bfM \bxone \bfN}", curve={height=-18pt}, from=6-7, to=12-5]
	\arrow[equals, from=7-1, to=8-1]
	\arrow[""{name=4, anchor=center, inner sep=0}, "{\mathbbm{1}_\Cc\bxone \runtwo_\Cc}", dashed, from=7-2, to=7-1]
	\arrow[""{name=5, anchor=center, inner sep=0}, "{\asstwo_{\Cc,\Cc,\Vv}}", dashed, from=7-3, to=7-2]
	\arrow[equals, from=7-3, to=8-3]
	\arrow[""{name=6, anchor=center, inner sep=0}, "{\bfM \bxone \bfN}"', from=8-1, to=12-1]
	\arrow["{\runtwo_{\Cc\bxob \Cc}}", from=8-3, to=8-1]
	\arrow["{(\bfM \bxone \bfN)\bxone \bfT_{\Vv}}", from=8-3, to=10-3]
	\arrow[between={0.4}{0.6}, equals, from=10-3, to=6-5]
	\arrow["{\asstwo_{\A,\A,\A}}"', from=10-3, to=10-5]
	\arrow[""{name=7, anchor=center, inner sep=0}, "{\what\bxone \mathbbm{1}_{\A}}", from=10-3, to=12-3]
	\arrow[""{name=8, anchor=center, inner sep=0}, "{\mathbbm{1}_{\A}\bxone \what}"', color={rgb,255:red,92;green,92;blue,214}, from=10-5, to=12-5]
	\arrow["\what"', curve={height=18pt}, from=12-1, to=13-4]
	\arrow["\what"', from=12-3, to=13-4]
	\arrow["\what", color={rgb,255:red,92;green,92;blue,214}, from=12-5, to=13-4]
	\arrow["{\theta^{r*}}"', between={0.4}{0.6}, Rightarrow, from=0, to=1]
	\arrow["\delta", between={0.4}{0.6}, Rightarrow, from=2, to=2-5]
	\arrow["{\mathbbm{1}_{\Cc}\bxone \theta^{r*}}", between={0.4}{0.6}, Rightarrow, from=3, to=6-5]
	\arrow["{\tilde\rhotwo_{\Cc,\Cc}}", between={0.4}{0.6}, Rightarrow, dashed, from=6-1, to=5]
	\arrow["{{\tilde\rhotwo_{\Cc,\Cc}}^{-1}}"', between={0.4}{0.6}, Rightarrow, dashed, from=4, to=8-3]
	\arrow["{\tau^{r*}_{\bfM\tenrep_0 \bfN}}"', between={0.4}{0.6}, Rightarrow, from=6, to=7]
	\arrow["\alpha"', between={0.4}{0.6}, Rightarrow, from=12-3, to=8]
\end{tikzcd}
}
\end{align}

We again start by unpacking the definition  of  $\mathbbm{1}_{\bfM}\tenrepC_1 \runthree^*_{\bfN}$. By definition,
\[\mathbbm{1}_{\bfM}\tenrepC_1 \runthree^*_{\bfN}=(\mathbbm{1}_{\bfM}\tenrepC \runthree^*_{\bfN})\star\Delta = \what\star (\mathbbm{1}_{\bfM}\bxtwo \runthree^*_{\bfN})\star\Delta.\]
Using that 
\[\runthree^*_{\bfN} = (\tau_\bfN^{r*} \star( (\mathbbm{1}_\Cc\bxone \Theta)\Delta)    ) (\bfN\star \theta^{r*})\]
we can see, in analogy to the computation for $\lunthree_{\bfM}\tenrepC_0\mathbbm{1}_\bfN$, that
\[ \mathbbm{1}_{\bfM}\tenrepC_1 \runthree^*_{\bfN} = \bigl( \what\star(\mathbbm{1}_\bfM\bxtwo \tau^{r*}_\bfN)\star ( (\mathbbm{1}_\Cc\bxone(\mathbbm{1}_\Cc \bxone \Theta))(\mathbbm{1}_\Cc\bxone \Delta)\Delta    )    \bigr)   \bigl( (\bfM\tenrep_0\bfN   )\star (\mathbbm{1}_{\mathbbm{1}_\Cc} \bxtwo \theta^{r*}) \star \Delta  \bigr). \]

Denote
\[\overline{\mathbbm{1}_\bfM\bxtwo \tau^{r*}_\bfN} =\bigl( \what\star(\mathbbm{1}_\bfM\bxtwo \tau^{r*}_\bfN)\star ( (\mathbbm{1}_\Cc\bxone(\mathbbm{1}_\Cc \bxone \Theta))(\mathbbm{1}_\Cc\bxone \Delta)\Delta    )    \bigr) \]
and 
\[  \overline{\mathbbm{1}_{\mathbbm{1}_\Cc} \bxtwo \theta^{r*}} = \bigl( (\bfM\tenrep_0\bfN   )\star (\mathbbm{1}_{\mathbbm{1}_\Cc} \bxtwo \theta^{r*}) \star \Delta  \bigr) \]
we thus have 
\[\mathbbm{1}_{\bfM}\tenrepC_1 \runthree^*_{\bfN} = \overline{\mathbbm{1}_\bfM\bxtwo \tau^{r*}_\bfN}\, \overline{\mathbbm{1}_{\mathbbm{1}_\Cc} \bxtwo \theta^{r*}} .\]

We now consider the codomain $ \assthree_{\bfM,\bfN,\bfT_\Cc } \,\runthree^*_{\bfM\tenrepC_0\bfN}$. Recall the strong transformations $\runtwo_\Delta$ from Proposition \ref{unitortwocommute}.

\begin{lemma}
We have $\tau^{r*}_{\bfM\tenrepC_0\bfN}= (\tau^{r*}_{\bfM\tenrep_0\bfN}\star (\Delta\bxone\mathbbm{1}_{\Vv}))((\bfM\tenrep_0\bfN)\star \kappa_r)
$ as strong transformations from $(\bfM\tenrepC_0\bfN) \runtwo_\Cc$ to $(\bfM\tenrepC_0\bfN)\tenrep_0 \bfT_\Vv$.
\end{lemma}

\proof
The $1$-cell components at an object $(\ob{i}, \bullet)$ in both cases given by $\run^*_{\bfM(\ob{i}_\ell) \monob \bfN(\ob{i}_{r})}$. Thus the source and target of the $2$-cell components of both strong transformations are the same, and since both are given by resummations, they coincide.
\endproof

Defining
\[\overline{\theta^{r*}} = ((\bfM\tenrep_0\bfN)\Delta) \star \theta^{r*},\]
\[ \overline{\runtwo_\Delta} =  (\bfM\tenrep_0\bfN) \star \runtwo_\Delta\star ((\mathbbm{1}_{\Cc}\bxone \Delta)\Delta) \]
and
\begin{equation*}\begin{split}  
\overline{\tau^{r*}_{\bfM\tenrep_0\bfN}} &=\tau^{r*}_{\bfM\tenrep_0\bfN}\star \bigl( (\Delta\bxone \mathbbm{1}_\Vv)(\mathbbm{1}_{\Cc}\bxone \Theta)\Delta \bigr) \\
&=\tau^{r*}_{\bfM\tenrep_0\bfN}\star \bigl(((\mathbbm{1}_{\Cc}\bxone\mathbbm{1}_{\Cc})\bxone \Theta) (\Delta\bxone \mathbbm{1}_\Cc)\Delta \bigr)
\end{split}\end{equation*}
we thus have
\[\runthree^*_{\bfM\tenrepC_0\bfN} = \overline{\tau^{r*}_{\bfM\tenrep_0\bfN}} \,\overline{\runtwo_\Delta} \, \overline{\theta^{r*}}. \]

Recall the invertible strong natural transformation $\tilde\rhotwo_{\Cc,\Cc}$ defined in Section \ref{phirdef}. Defining 
\[\overline{\tilde\rhotwo_{\Cc,\Cc}} = (\bfM\tenrep_0\bfN)\star\tilde\rhotwo_{\Cc,\Cc} \star \bigl( ((\mathbbm{1}_{\Cc}\bxone\mathbbm{1}_{\Cc})\bxone \Theta) (\Delta\bxone \mathbbm{1}_\Cc)\Delta  \bigr)\]
and inserting 

$\overline{\tilde\rhotwo_{\Cc,\Cc}} $ and $(\overline{\tilde\rhotwo_{\Cc,\Cc}} )^{-1}$ as in \eqref{rhounitor2cell}, we have 
\[\runthree^*_{\bfM\tenrepC_0\bfN} = \overline{\tau^{r*}_{\bfM\tenrep_0\bfN}} \,(\overline{\tilde\rhotwo_{\Cc,\Cc}} )^{-1}\,\overline{\tilde\rhotwo_{\Cc,\Cc}} \,\overline{\runtwo_\Delta} \, \overline{\theta^{r*}}. \]

As before, we have 
\[ \assthree_{\bfM,\bfN,\bfT_\Cc } = \overline{\delta}\,\overline{\alpha}\]
where now
\[
\overline{\delta} = ( \bfM\tenrep_0(\bfN\tenrep_0\bfT_\Cc) )\star \delta \]
and
\[\overline{\alpha} = \alpha\star \bigl( ( (\bfM\bxone\bfN)\bxone \bfT_\Cc  ) (\Delta\bxone \mathbbm{1}_\Cc)\Delta \bigr),\]
and thus
\[ \assthree_{\bfM,\bfN,\bfT_\Cc } \,\runthree^*_{\bfM\tenrepC_0\bfN} =   \overline{\delta}\,\overline{\alpha}\,\overline{\tau^{r*}_{\bfM\tenrep_0\bfN}} \,\overline{\tilde\rhotwo_{\Cc,\Cc}} ^{-1}\,\overline{\tilde\rhotwo_{\Cc,\Cc}} \,\overline{\runtwo_\Delta} \, \overline{\theta^{r*}}. \]

We therefore wish to show an equality
\begin{equation}\label{lambdathreeeq}
\overline{\mathbbm{1}_\bfM\bxtwo \tau^{r*}_\bfN}\, \overline{\mathbbm{1}_{\mathbbm{1}_\Cc} \bxtwo \theta^{r*}} 
=
\overline{\delta}\,\overline{\alpha}\,\overline{\tau^{r*}_{\bfM\tenrep_0\bfN}} \,\overline{\tilde\rhotwo_{\Cc,\Cc}} ^{-1}\,\overline{\tilde\rhotwo_{\Cc,\Cc}} \,\overline{\runtwo_\Delta} \, \overline{\theta^{r*}}.
\end{equation}

We now set 
\[\overline{\overline{\mathbbm{1}_\bfM\bxtwo \tau^{r*}_\bfN}} =\what\star(\mathbbm{1}_\bfM\bxtwo \tau^{r*}_\bfN)\star \bigl( \asstwo_{\Cc,\Cc,\Vv}((\mathbbm{1}_\Cc\bxone\mathbbm{1}_\Cc) \bxone \Theta)( \Delta\bxone\mathbbm{1}_\Cc)\Delta \bigr) \]
and claim that
\[\overline{\overline{\mathbbm{1}_\bfM\bxtwo \tau^{r*}_\bfN}} =\overline{\alpha}\,\overline{\tau^{r*}_{\bfM\tenrep_0\bfN}} \,\overline{\tilde\rhotwo_{\Cc,\Cc}} ^{-1}. \]
Indeed, this follows by pre-whiskering with $ ((\mathbbm{1}_\Cc\bxone\mathbbm{1}_\Cc) \bxone \Theta)( \Delta\bxone\mathbbm{1}_\Cc)\Delta  $
from the following lemma.

\begin{lemma}
We have
\[\what\star(\mathbbm{1}_\bfM\bxtwo \tau^{r*}_\bfN)\star \asstwo_{\Cc,\Cc,\Vv} =\bigl( \alpha\star  ((\bfM\bxone\bfN)\bxone \bfT_\Vv) \bigr)   \tau^{r*}_{\bfM\tenrep_0\bfN}\bigl(   (\bfM\tenrep_0\bfN)\star{\tilde\rhotwo_{\Cc,\Cc}}^{-1} \bigr)\]
as strong transformations from $(\bfM\tenrep_0\bfN)(\mathbbm{1}_{\Cc} \bxone \runtwo_\Cc)\asstwo_{\Cc,\Cc,\Vv}$ to $\bfM\tenrep_0(\bfN\tenrep_0\bfT_\Vv)\asstwo_{\Cc,\Cc,\Vv}$ of pseudofunctors from $(\Cc\bxob \Cc)\bxob\Vv$ to $\A$.
\end{lemma}

\proof
In order to see this, we compute the $1$-cells associated to the two strong natural transformations, which are given by functors in $\A$. More precisely, on an object 
$\bigoplus_t(X_t, Y_t) \in \bfM(\ob{i})\monob \bfN(\ob{j})$, we compute 
\begin{equation*}\begin{split}  
(\what\star(\mathbbm{1}_\bfM\bxtwo \tau^{r*}_\bfN)\star  \asstwo_{\Cc,\Cc,\Vv} )_{((\ob{i},\ob{j}),\bullet)}(\bigoplus_t(X_t, Y_t)) = \bigoplus_t(X_t, (Y_t, \bullet_\Bbbk))
\end{split}\end{equation*}
while 
\begin{equation*}\begin{split}  
\bigl(\bigl( \alpha\star  ((\bfM\bxone\bfN)\bxone \bfT_\Vv) \bigr)&   \tau^{r*}_{\bfM\tenrep_0\bfN}\bigl(   (\bfM\tenrep_0\bfN)\star{\tilde\rhotwo_{\Cc,\Cc}}^{-1} \bigr)\bigr)_{((\ob{i},\ob{j}),\bullet)}(\bigoplus_t(X_t, Y_t)) \\
&= \bigl(\bigl( \alpha\star  ((\bfM\bxone\bfN)\bxone \bfT_\Vv) \bigr)  \tau^{r*}_{\bfM\tenrep_0\bfN}\bigr)_{((\ob{i},\ob{j}),\bullet)} (\bigoplus_t(X_t, Y_t))\\
&=   \bigl( \alpha\star  ((\bfM\bxone\bfN)\bxone \bfT_\Vv) \bigr)_{((\ob{i},\ob{j}),\bullet)}     ( \bigoplus_t(X_t, Y_t), \bullet_\Bbbk)\\
&= \bigoplus_t(X_t, (Y_t, \bullet_\Bbbk))
\end{split}\end{equation*}
so the functors agree on objects and it is straightforward to see that they agree on morphisms as well. Hence the $1$-cell components of our strong natural transformations coincide. Thus the $2$-cell components of both have the same source and target, and since both are given by resummations, they coincide.
\endproof

We can therefore rewrite our desired equality \eqref{lambdathreeeq} as
\begin{equation}\label{lambdathreeeq2}
\overline{\mathbbm{1}_\bfM\bxtwo \tau^{r*}_\bfN}\, \overline{\mathbbm{1}_{\mathbbm{1}_\Cc} \bxtwo \theta^{r*}} 
=
\overline{\delta}\,\overline{\overline{\mathbbm{1}_\bfM\bxtwo \tau^{r*}_\bfN}}\,\overline{\tilde\rhotwo_{\Cc,\Cc}} \,\overline{\runtwo_\Delta} \, \overline{\theta^{r*}}.
\end{equation}

The next step is to swap $\overline{\delta}$ and $\overline{\overline{\mathbbm{1}_\bfM\bxtwo \tau^{r*}_\bfN}}$ to obtain
\begin{equation}\label{swapdeltatauinlambda}
\overline{\delta}\,\overline{\overline{\mathbbm{1}_\bfM\bxtwo \tau^{r*}_\bfN}}  =\overline{\mathbbm{1}_\bfM\bxtwo \tau^{r*}_\bfN} \,  \overline{\overline{\delta}} 
\end{equation}
where
\[\overline{\overline{\delta}} =\bigl( (\bfM\tenrep_0\bfN)(\mathbbm{1}_\Cc\bxone \runtwo_\Cc)\asstwo_{\Cc,\Cc,\Vv} (\mathbbm{1}_\Cc\bxone(\mathbbm{1}_\Cc \bxone \Theta))  \bigr)\star\delta.\]
To see this, we compute the modification between the two compositions using Lemma \ref{swaptwocells} and obtain the $1$-cell component
\[
(\mathbbm{1}_\bfM\bxtwo \tau^{r*}_\bfN)_{(\mathbbm{1}_\Cc\bxone(\mathbbm{1}_\Cc \bxone \Theta))(\delta_{\ob(i)})} = \mathbbm{1}_{\bfM(\delta_{\ob{i},1})} \otimes_2 
(\tau^{r*}_\bfN)_{(\delta_{\ob{i},2}, \mathbbm{1}_\bullet)}.
\]
Evaluated on any object $X\in \bfN(\ob{i}_{\ell r})$,  the map
\[\bigl((\tau^{r*}_\bfN)_{(\delta_{\ob{i},1}, \mathbbm{1}_\bullet)}\bigr)_X \colon (\bfN(\delta_{\ob{i},2})X, \bullet_\Bbbk) \to  (\bfN(\delta_{\ob{i},2})X, \bullet_\Bbbk)  \]
is the identity, thus the modification is the identity modification and the equality in \eqref{swapdeltatauinlambda} holds.

We deduce that \eqref{swapdeltatauinlambda} is equivalent to
\begin{equation}\label{lambdathreeeq3}
\overline{\mathbbm{1}_\bfM\bxtwo \tau^{r*}_\bfN}\, \overline{\mathbbm{1}_{\mathbbm{1}_\Cc} \bxtwo \theta^{r*}} 
=
\overline{\mathbbm{1}_\bfM\bxtwo \tau^{r*}_\bfN}\,\overline{\overline{\delta}}\,\overline{\tilde\rhotwo_{\Cc,\Cc}} \,\overline{\runtwo_\Delta} \, \overline{\theta^{r*}}.
\end{equation}
which is guaranteed by the equality
\begin{equation*}
 \overline{\mathbbm{1}_{\mathbbm{1}_\Cc} \bxtwo \theta^{r*}} 
=
\overline{\overline{\delta}}\,\overline{\tilde\rhotwo_{\Cc,\Cc}} \,\overline{\runtwo_\Delta} \, \overline{\theta^{r*}}.
\end{equation*}
obtained by post-whiskering the equality in condition
 \eqref{rightcond} in the definition of a comonoidal structure with $\bfM\tenrep_0\bfN$.

\subsubsection{The pentagonator $\pithree$} (Incomplete)
We now consider the domain  and codomain 
of $\pithree_{\bfK, \bfL,\bfM,\bfN}$.
As above we will give non unique diagrammatic representations of these natural transformations, which allow us to describe the intuitive idea of what follows, before explaining the details.
The domain
$(\mathbbm{1}_{\bfK}\tenrepC \assthree_{ \bfL, \bfM,\bfN})\assthree_{\bfK, \bfL\tenrepC \bfM,\bfN}(\assthree_{\bfK, \bfL, \bfM}\tenrepC\mathbbm{1}_{\bfN})$ can be represented, again non-uniquely, by the following diagram:

\begin{align}\label{pentadom}
\scalebox{0.55}{
\begin{tikzcd}[ampersand replacement=\&]
	\&\&\& \Cc \\
	\\
	\& {\Cc\bxob\Cc} \&\&\&\& {\Cc\bxob\Cc} \\
	\\
	{(\Cc\bxob\Cc)\bxob \Cc} \&\& {(\Cc\bxob\Cc)\bxob \Cc} \&\& {\Cc\bxob(\Cc\bxob\Cc)} \&\& {\Cc\bxob(\Cc\bxob\Cc)} \\
	\\
	{((\Cc\bxob\Cc)\bxob\Cc)\bxob \Cc} \&\& {(\Cc\bxob(\Cc\bxob\Cc))\bxob\Cc} \&\& {\Cc\bxob ((\Cc\bxob\Cc)\bxob\Cc)} \&\& {\Cc\bxob (\Cc\bxob(\Cc\bxob\Cc))} \\
	\\
	{((\A\bxob\A)\bxob\A)\bxob\A} \&\& {(\A\bxob(\A\bxob\A))\bxob\A} \&\& {\A\bxob((\A\bxob\A)\bxob\A)} \&\& {\A\bxob(\A\bxob(\A\bxob\A))} \\
	\\
	{(\A\bxob\A)\bxob\A} \&\& {(\A\bxob\A)\bxob\A} \&\& {\A\bxob(\A\bxob\A)} \&\& {\A\bxob(\A\bxob\A)} \\
	\\
	\& {\A\bxob\A} \&\&\&\& {\A\bxob\A} \\
	\\
	\&\&\& \A
	\arrow["\Delta"', from=1-4, to=3-2]
	\arrow["\Delta", from=1-4, to=3-6]
	\arrow["{\Delta\bxone\mathbbm{1}_\Cc}"', from=3-2, to=5-1]
	\arrow[""{name=0, anchor=center, inner sep=0}, "{\Delta\bxone\mathbbm{1}_\Cc}", from=3-2, to=5-3]
	\arrow["{\mathbbm{1}_\Cc\bxone \Delta}"', from=3-6, to=5-5]
	\arrow["{\mathbbm{1}_\Cc\bxone \Delta}", from=3-6, to=5-7]
	\arrow[""{name=1, anchor=center, inner sep=0}, "{(\Delta\bxone\mathbbm{1}_\Cc)\bxone\mathbbm{1}_\Cc}"', from=5-1, to=7-1]
	\arrow["{\asstwo_{\Cc,\Cc,\Cc}}", from=5-3, to=5-5]
	\arrow["{(\mathbbm{1}_\Cc\bxone \Delta)\bxone\mathbbm{1}_\Cc}"', from=5-3, to=7-3]
	\arrow[""{name=2, anchor=center, inner sep=0}, "{\mathbbm{1}_\Cc\bxone (\Delta\bxone\mathbbm{1}_\Cc)}", from=5-5, to=7-5]
	\arrow["{\mathbbm{1}_\Cc\bxone(\mathbbm{1}_\Cc\bxone \Delta)}", from=5-7, to=7-7]
	\arrow["{\asstwo_{\Cc,\Cc,\Cc}\bxone\mathbbm{1}_\Cc}", from=7-1, to=7-3]
	\arrow["{((\bfK\bxone\bfL)\bxone\bfM)\bxone \bfN}"', from=7-1, to=9-1]
	\arrow[between={0.4}{0.6}, equals, from=7-3, to=5-5]
	\arrow["{\asstwo_{\Cc,\Cc\bxob\Cc,\Cc}}", from=7-3, to=7-5]
	\arrow["{(\bfK\bxone(\bfL\bxone\bfM))\bxone\bfN}"', from=7-3, to=9-3]
	\arrow["{\mathbbm{1}_{\Cc}\bxone\asstwo_{\Cc,\Cc,\Cc}}", from=7-5, to=7-7]
	\arrow["{\bfK\bxone((\bfL\bxone\bfM)\bxone\bfN)}", from=7-5, to=9-5]
	\arrow["{\bfK\bxone(\bfL\bxone(\bfN\bxone \bfM))}", from=7-7, to=9-7]
	\arrow[between={0.4}{0.6}, equals, from=9-1, to=7-3]
	\arrow["{\asstwo_{\A,\A,\A}\bxone\mathbbm{1}_{\A}}", from=9-1, to=9-3]
	\arrow["{(\what\bxone\mathbbm{1}_{\A})\bxone\mathbbm{1}_{\A}}"', from=9-1, to=11-1]
	\arrow[between={0.4}{0.6}, equals, from=9-3, to=7-5]
	\arrow["{\asstwo_{\A,\A\bxob\A,\A}}", from=9-3, to=9-5]
	\arrow[""{name=3, anchor=center, inner sep=0}, "{(\mathbbm{1}_{\A}\bxone \what)\bxone \mathbbm{1}_{\A}}"', from=9-3, to=11-3]
	\arrow[between={0.4}{0.6}, equals, from=9-5, to=7-7]
	\arrow["{\mathbbm{1}\bxone \asstwo_{\A,\A,\A}}", from=9-5, to=9-7]
	\arrow["{\mathbbm{1}_{\A}\bxone (\what\bxone\mathbbm{1}_{\A})}", from=9-5, to=11-5]
	\arrow[""{name=4, anchor=center, inner sep=0}, "{\mathbbm{1}_{\A}\bxone(\mathbbm{1}_{\A}\bxone \what)}", from=9-7, to=11-7]
	\arrow["{\what\bxone\mathbbm{1}_{\A}}"', from=11-1, to=13-2]
	\arrow[between={0.4}{0.6}, equals, from=11-3, to=9-5]
	\arrow["{\asstwo_{\A,\A,\A}}", from=11-3, to=11-5]
	\arrow["{\what\bxone\mathbbm{1}_{\A}}", from=11-3, to=13-2]
	\arrow[""{name=5, anchor=center, inner sep=0}, "{\mathbbm{1}_{\A}\bxone \what}"', from=11-5, to=13-6]
	\arrow["{\mathbbm{1}_{\A}\bxone \what}", from=11-7, to=13-6]
	\arrow["\what"', from=13-2, to=15-4]
	\arrow["\what", from=13-6, to=15-4]
	\arrow["\delta", between={0.4}{0.6}, Rightarrow, from=0, to=3-6]
	\arrow["{\delta\bxtwo\mathbbm{1}_{\mathbbm{1}_\Cc}}", between={0.4}{0.6}, Rightarrow, from=1, to=5-3]
	\arrow["{\mathbbm{1}_{\mathbbm{1}_\Cc}\bxtwo \delta}", between={0.4}{0.6}, Rightarrow, from=2, to=5-7]
	\arrow["{\alpha\bxtwo\mathbbm{1}_{\mathbbm{1}_{\A}}}", between={0.4}{0.6}, Rightarrow, from=11-1, to=3]
	\arrow["{\mathbbm{1}_{\mathbbm{1}_{\A}}\bxtwo\alpha}", between={0.4}{0.6}, Rightarrow, from=11-5, to=4]
	\arrow["\alpha", between={0.4}{0.6}, Rightarrow, from=13-2, to=5]
\end{tikzcd}}
\end{align}

The codomain $\assthree_{\bfK, \bfL,\bfM\tenrepC\bfN } \,\assthree_{\bfK\tenrepC \bfL,\bfM,\bfN}$, on the other hand, may be represented by:

\begin{align}\label{pentacodom}
\scalebox{0.5}{
\begin{tikzcd}[ampersand replacement=\&]
	\&\&\& \Cc \&\&\& \\
	\\
	\& {\Cc\bxob\Cc} \&\& {\Cc\bxob\Cc} \&\& {\Cc\bxob\Cc} \\
	\\
	{(\Cc\bxob\Cc)\bxob \Cc} \&\& {\Cc\bxob(\Cc\bxob\Cc)} \&\& {(\Cc\bxob\Cc)\bxob \Cc} \&\& {\Cc\bxob(\Cc\bxob \Cc)} \\
	\\
	{((\Cc\bxob\Cc)\bxob\Cc)\bxob\Cc} \&\&\& {(\Cc\bxob \Cc)\bxob(\Cc\bxob\Cc)} \&\&\& {\Cc \bxob(\Cc\bxob(\Cc\bxob\Cc))} \\
	\\
	{((\A\bxob\A)\bxob\A)\bxob\A} \&\&\& {(\A\bxob\A)\bxob(\A\bxob\A)} \&\&\& {\A\bxob(\A\bxob(\A\bxob\A))} \\
	\\
	{(\A\bxob\A)\bxob\A} \&\& {\A\bxob(\A\bxob\A)} \&\& {(\A\bxob\A)\bxob\A} \&\& {\A\bxob(\A\bxob\A)} \\
	\\
	\& {\A\bxob\A} \&\& {\A\bxob\A} \&\& {\A\bxob\A} \\
	\\
	\&\&\& \A
	\arrow["\Delta"', from=1-4, to=3-2]
	\arrow["\Delta", from=1-4, to=3-4]
	\arrow["\Delta"', from=1-4, to=3-6]
	\arrow[""{name=0, anchor=center, inner sep=0}, "{\Delta\bxone\mathbbm{1}_\Cc}", from=3-2, to=5-1]
	\arrow["{\mathbbm{1}_\Cc\bxone \Delta}"', from=3-4, to=5-3]
	\arrow[""{name=1, anchor=center, inner sep=0}, "{\Delta\bxone\mathbbm{1}_\Cc}", from=3-4, to=5-5]
	\arrow["{\mathbbm{1}_\Cc\bxone \Delta}"', from=3-6, to=5-7]
	\arrow["{\asstwo_{\Cc,\Cc,\Cc}}", from=5-1, to=5-3]
	\arrow["{(\Delta\bxone\mathbbm{1}_\Cc)\bxone\mathbbm{1}_\Cc}"', from=5-1, to=7-1]
	\arrow[between={0.4}{0.6}, equals, from=5-3, to=5-5]
	\arrow["{\Delta\bxone (\mathbbm{1}_\Cc\bxone\mathbbm{1}_\Cc)}", from=5-3, to=7-4]
	\arrow["{\asstwo_{\Cc,\Cc,\Cc}}", from=5-5, to=5-7]
	\arrow["{(\mathbbm{1}_\Cc\bxone\mathbbm{1}_\Cc)\bxone \Delta}", from=5-5, to=7-4]
	\arrow["{\mathbbm{1}_\Cc\bxone(\mathbbm{1}_\Cc\bxone \Delta)}", from=5-7, to=7-7]
	\arrow[between={0.4}{0.6}, equals, from=7-1, to=5-3]
	\arrow["{\asstwo_{\Cc\bxob\Cc,\Cc,\Cc}}", from=7-1, to=7-4]
	\arrow["{((\bfK\bxone\bfL)\bxone\bfM)\bxone\bfN}"', from=7-1, to=9-1]
	\arrow[between={0.4}{0.6}, equals, from=7-4, to=5-7]
	\arrow["{\asstwo_{\Cc,\Cc,\Cc\bxob\Cc}}", from=7-4, to=7-7]
	\arrow["{(\bfK\bxone\bfL)\bxone(\bfM\bxone\bfN)}", from=7-4, to=9-4]
	\arrow["{\bfK\bxone(\bfL\bxone(\bfM\bxone\bfN))}", from=7-7, to=9-7]
	\arrow[between={0.4}{0.6}, equals, from=9-1, to=7-4]
	\arrow["{\asstwo_{\A\bxob\A,\A,\A}}", from=9-1, to=9-4]
	\arrow["{(\what\bxone\mathbbm{1}_{\A})\bxone \mathbbm{1}_{\A}}"', from=9-1, to=11-1]
	\arrow[between={0.4}{0.6}, equals, from=9-4, to=7-7]
	\arrow["{\asstwo_{\A,\A,\A\bxob\A}}", from=9-4, to=9-7]
	\arrow["{\what\bxone (\mathbbm{1}_{\A}\bxone\mathbbm{1}_{\A})}", from=9-4, to=11-3]
	\arrow["{(\mathbbm{1}_{\A}\bxone\mathbbm{1}_{\A})\bxone \what}", from=9-4, to=11-5]
	\arrow["{\mathbbm{1}_{\A}\bxone(\mathbbm{1}_{\A}\bxone \what)}", from=9-7, to=11-7]
	\arrow[between={0.4}{0.6}, equals, from=11-1, to=9-4]
	\arrow["{\asstwo_{\A,\A,\A}}", from=11-1, to=11-3]
	\arrow["{\what\bxone\mathbbm{1}_{\A}}", from=11-1, to=13-2]
	\arrow[between={0.4}{0.6}, equals, from=11-3, to=11-5]
	\arrow[""{name=2, anchor=center, inner sep=0}, "{\mathbbm{1}_{\A}\bxone \what}"', from=11-3, to=13-4]
	\arrow[between={0.4}{0.6}, equals, from=11-5, to=9-7]
	\arrow["{\asstwo_{\A,\A,\A}}", from=11-5, to=11-7]
	\arrow["{\what \bxone \mathbbm{1}_{\A}}", from=11-5, to=13-4]
	\arrow[""{name=3, anchor=center, inner sep=0}, "{\mathbbm{1}_{\A}\bxone \what}"', from=11-7, to=13-6]
	\arrow["\what"', from=13-2, to=15-4]
	\arrow["\what", from=13-4, to=15-4]
	\arrow["\what"', from=13-6, to=15-4]
	\arrow["\delta", between={0.4}{0.6}, Rightarrow, from=0, to=3-4]
	\arrow["\delta", between={0.4}{0.6}, Rightarrow, from=1, to=3-6]
	\arrow["\alpha", between={0.4}{0.6}, Rightarrow, from=13-2, to=2]
	\arrow["\alpha"', between={0.4}{0.6}, Rightarrow, from=13-4, to=3]
\end{tikzcd}}
\end{align}

We begin by carefully unpacking the domain and codomain of
$\pithree_{\bfK, \bfL,\bfM,\bfN}$.

We have that
\begin{align*}
(\assthree_{\bfK, \bfL, \bfM}\tenrepC_1\mathbbm{1}_{\bfN}) &= \what \star (\assthree_{\bfK, \bfL, \bfM}\bxtwo  \mathbbm{1}_{\bfN})\star \Delta \\
&\!
\begin{multlined}[0.6\columnwidth]
	=\what \star 
	\big[
	\big(
		((\bfK\tenrep_0(\bfL\tenrep_0\bfM))\bxone \bfN) \star
		(\delta \bxtwo \mathbbm{1}_{\mathbbm{1}_{\Cc}})
\big)\\ 
\big(
	(\alpha \bxtwo  \mathbbm{1}_{\mathbbm{1}_{\A}}) \star ((((\bfK\bxone\bfL)\bxone\bfM)\bxone \bfN)((\Delta \bxone \mathbbm{1}_\Cc)\bxone \mathbbm{1}_\Cc)(\Delta \bxone \mathbbm{1}_\Cc))
\big) 
\big]
\star \Delta
\end{multlined}\\
&\!
\begin{multlined}[0.6\columnwidth]
	=
	\big(
		((\bfK\tenrep_0(\bfL\tenrep_0\bfM))\tenrep_0 \bfN) \star
		(\delta \bxtwo \mathbbm{1}_{\mathbbm{1}_{\Cc}}) \star \Delta
\big)\\ 
\big(
	\what \star (\alpha \bxtwo  \mathbbm{1}_{\mathbbm{1}_{\A}}) \star ((((\bfK\bxone\bfL)\bxone\bfM)\bxone \bfN)((\Delta \bxone \mathbbm{1}_\Cc)\bxone \mathbbm{1}_\Cc)(\Delta \bxone \mathbbm{1}_\Cc)\Delta)
\big) .
\end{multlined}
\end{align*}

Let 
\[
\overline{\delta^{\ell}} = ((\bfK\tenrep_0(\bfL\tenrep_0\bfM))\tenrep_0 \bfN) \star
		(\delta \bxtwo \mathbbm{1}_{\mathbbm{1}_{\Cc}}) \star \Delta
\]
and 
\[
\overline{\alpha^{\ell}} = \what \star (\alpha \bxtwo  \mathbbm{1}_{\mathbbm{1}_{\A}}) \star ((((\bfK\bxone\bfL)\bxone\bfM)\bxone \bfN)((\Delta \bxone \mathbbm{1}_\Cc)\bxone \mathbbm{1}_\Cc)(\Delta \bxone \mathbbm{1}_\Cc)\Delta).
\]

Similarly 
\begin{align*}
(\mathbbm{1}_{\bfK}\tenrepC_1 \assthree_{\bfL, \bfM, \bfN}) &= \what \star ( \mathbbm{1}_{\bfK}\bxtwo \assthree_{\bfL, \bfM, \bfN})\star \Delta \\
&\!
\begin{multlined}[0.6\columnwidth]
	=\what \star 
	\big[
	\big(
		(\bfK\bxone (\bfL\tenrep_0(\bfM\tenrep_0\bfN))) \star
		(\mathbbm{1}_{\mathbbm{1}_{\Cc}}\bxtwo \delta )
\big)\\ 
\big(
	(\mathbbm{1}_{\mathbbm{1}_{\A}} \bxtwo \alpha ) \star ((\bfK\bxone((\bfL\bxone\bfM)\bxone \bfN)(\mathbbm{1}_\Cc\bxone(\Delta \bxone \mathbbm{1}_\Cc))( \mathbbm{1}_\Cc\bxone \Delta))
\big) 
\big]
\star \Delta
\end{multlined}\\
&\!
\begin{multlined}[0.6\columnwidth]
	=
	\big(
		(\bfK\tenrep_0 (\bfL\tenrep_0(\bfM\tenrep_0\bfN))) \star
		(\mathbbm{1}_{\mathbbm{1}_{\Cc}}\bxtwo \delta ) \star \Delta
\big)\\ 
\big(
	\what \star (\mathbbm{1}_{\mathbbm{1}_{\A}} \bxtwo \alpha ) \star ((\bfK\bxone((\bfL\bxone\bfM)\bxone \bfN)(\mathbbm{1}_\Cc\bxone(\Delta \bxone \mathbbm{1}_\Cc))( \mathbbm{1}_\Cc\bxone \Delta)\Delta)
\big) .
\end{multlined}
\end{align*}
Let 
\[
\overline{\delta^r} = (\bfK\tenrep_0 (\bfL\tenrep_0(\bfM\tenrep_0\bfN))) \star
		(\mathbbm{1}_{\mathbbm{1}_{\Cc}}\bxtwo \delta ) \star \Delta
\]
and
\[
\overline{\alpha^r} = \what \star (\mathbbm{1}_{\mathbbm{1}_{\A}} \bxtwo \alpha ) \star ((\bfK\bxone((\bfL\bxone\bfM)\bxone \bfN)(\mathbbm{1}_\Cc\bxone(\Delta \bxone \mathbbm{1}_\Cc))( \mathbbm{1}_\Cc\bxone \Delta)\Delta).
\]

Now

\begin{align*}
\assthree_{\bfK\tenrepC \bfL,\bfM,\bfN} 
&\!
\begin{multlined}[t][0.6\columnwidth]
	=
	\big(
		((\bfK\tenrepC_0 \bfL)\tenrep_0(\bfM\tenrep_0\bfN)) \star\delta 
\big)\\ 
\big(
	\alpha \star ((((\bfK\tenrepC_0\bfL)\bxone\bfM)\bxone \bfN)(\Delta \bxone \mathbbm{1}_\Cc)\Delta)
\big) 
\end{multlined}
\end{align*}

\begin{align*}
\assthree_{\bfK, \bfL \tenrepC_0\bfM, \bfN} 
&\!
\begin{multlined}[t][0.6\columnwidth]
	=
	\big(
		(\bfK\tenrep_0 (\bfL\tenrepC_0\bfM)\tenrep_0\bfN) \star\delta 
\big)\\ 
\big(
	\alpha \star ((\bfK\bxone(\bfL\tenrepC_0\bfM)\bxone \bfN)(\Delta \bxone \mathbbm{1}_\Cc)\Delta)
\big) 
\end{multlined}
\end{align*}
and

\begin{align*}
\assthree_{\bfK, \bfL,\bfM\tenrepC\bfN }
&\!
\begin{multlined}[t][0.6\columnwidth]
	=
	\big(
		(\bfK\tenrep_0 (\bfL\tenrep_0(\bfM\tenrepC_0\bfN))) \star\delta 
\big)\\ 
\big(
	\alpha \star ((\bfK\bxone(\bfL\bxone(\bfM\tenrepC_0\bfN)))(\Delta \bxone \mathbbm{1}_\Cc)\Delta)
\big) .
\end{multlined}
\end{align*}
Let 
\begin{align*}\overline{\delta^{m\ell}}
&=((\bfK\tenrepC_0 \bfL)\tenrep_0(\bfM\tenrep_0\bfN)) \star\delta, \\
\overline{\delta^{mm}}
&= (\bfK\tenrep_0 (\bfL\tenrepC_0\bfM)\tenrep_0\bfN) \star\delta, \text{ and} \\
\overline{\delta^{mr}}
&=(\bfK\tenrep_0 (\bfL\tenrep_0(\bfM\tenrepC_0\bfN))) \star\delta .
\end{align*}
Similarly let
\begin{align*}
\overline{\alpha^{m\ell}}
&= \alpha \star ((((\bfK\tenrepC_0\bfL)\bxone\bfM)\bxone \bfN)(\Delta \bxone \mathbbm{1}_\Cc)\Delta)\\
\overline{\alpha^{mm}}
&= \alpha \star ((\bfK\bxone(\bfL\tenrepC_0\bfM)\bxone \bfN)(\Delta \bxone \mathbbm{1}_\Cc)\Delta), \text{ and}\\
\overline{\alpha^{mr}}
&= \alpha \star (((\bfK\bxone\bfL)\bxone(\bfM\tenrepC_0\bfN)))(\Delta \bxone \mathbbm{1}_\Cc)\Delta).
\end{align*}

Then the domain $(\mathbbm{1}_{\bfK}\tenrepC_1 \assthree_{ \bfL, \bfM,\bfN})\assthree_{\bfK, \bfL\tenrepC_0 \bfM,\bfN}(\assthree_{\bfK, \bfL, \bfM}\tenrepC_1\mathbbm{1}_{\bfN})$ of $\pithree_{\bfK, \bfL,\bfM,\bfN}$ is given by 
\[\overline{\delta^{r}}\;\overline{\alpha^{r}}\;\overline{\delta^{mm}}\;\overline{\alpha^{mm}}\;\overline{\delta^{\ell}}\;\overline{\alpha^{\ell}}\] 
and the codomain $\assthree_{\bfK, \bfL,\bfM\tenrepC\bfN } \,\assthree_{\bfK\tenrepC \bfL,\bfM,\bfN}$ of $\pithree_{\bfK, \bfL,\bfM,\bfN}$ is given by \[\overline{\delta^{mr}}\;\overline{\alpha^{mr}}\;\overline{\delta^{m\ell}}\;\overline{\alpha^{m\ell}}\;.\]

Now we define two new strong transformations given by composing the $2$-cells in the diagrams \eqref{pentadom} and  \eqref{pentacodom} in a different order. To this end, set
\begin{align*}
\overline{\overline{\alpha^{r}}}&=
\begin{multlined}[t][0.8\columnwidth]\what \star (\mathbbm{1}_{\mathbbm{1}_{\A}} \bxtwo \alpha ) \star \big(
(\bfK\bxone((\bfL\bxone\bfM)\bxone \bfN)) \\\asstwo_{\Cc,\Cc\bxob \Cc,\Cc}
(\asstwo_{\Cc,\Cc,\Cc}\bxone\mathbbm{1}_\Cc) (( \Delta \bxone\mathbbm{1}_\Cc)\bxone\mathbbm{1}_\Cc)( \Delta \bxone\mathbbm{1}_\Cc)\Delta)\big)
 ,
\end{multlined}\\
\overline{\overline{\alpha^{mm}}}&=  \alpha \star \big( ((\bfK\bxone(\bfL\tenrep_0\bfM))\bxone \bfN)
(\asstwo_{\Cc,\Cc,\Cc}\bxone\mathbbm{1}_\Cc) (( \Delta \bxone\mathbbm{1}_\Cc)\bxone\mathbbm{1}_\Cc)( \Delta \bxone\mathbbm{1}_\Cc)\Delta \big)
,\\
\overline{\overline{\alpha^{mr}}}&=\alpha \star \big(((\bfK\bxone\bfL)\bxone(\bfM\tenrep_0\bfN))  \asstwo_{\Cc\bxob\Cc,\Cc,\Cc}  ((\Delta \bxone \mathbbm{1}_\Cc) \bxone \mathbbm{1}_\Cc)   (\Delta \bxone \mathbbm{1}_\Cc)\Delta\big)
,\\
\overline{\overline{\delta^{\ell}}}&=\big((\bfK\tenrep_0(\bfL\tenrep_0(\bfM\tenrep_0 \bfN)))(\mathbbm{1}_\Cc\bxone\asstwo_{\Cc,\Cc,\Cc})  \asstwo_{\Cc,\Cc\bxob \Cc,\Cc}\big)\star (\delta \bxtwo \mathbbm{1}_{\mathbbm{1}_{\Cc}}) \star \Delta
,\\
\overline{\overline{\delta^{m\ell}}}&= \big( (\bfK\tenrep_0(\bfL\tenrep_0(\bfM\tenrep_0 \bfN)))  \asstwo_{\Cc,\Cc,\Cc\bxob \Cc} ( \Delta \bxone(\mathbbm{1}_\Cc\bxone\mathbbm{1}_\Cc))\big)\star\delta,
 \text{ and}\\
\overline{\overline{\delta^{mm}}}&=   \big( (\bfK\tenrep_0(\bfL\tenrep_0(\bfM\tenrep_0 \bfN)))(\mathbbm{1}_\Cc\bxone\asstwo_{\Cc,\Cc,\Cc})    (\mathbbm{1}_\Cc \bxone( \Delta\bxone\mathbbm{1}_\Cc))\big)\star \delta
\end{align*}
and define the two strong transformations
\begin{equation}\label{domordered}
\overline{\delta^{r}}\;\overline{\overline{\delta^{mm}}}\;\overline{\overline{\delta^{\ell}}}\;\overline{\overline{\alpha^{r}}}\;\overline{\overline{\alpha^{mm}}}\;\overline{\alpha^{\ell}}
\end{equation}
and 
\begin{equation}\label{codomordered}
\overline{\delta^{mr}}\;\overline{\overline{\delta^{m\ell}}}\; \overline{\overline{\alpha^{mr}}} \;\overline{\alpha^{m\ell}},
\end{equation}
where, to make sense of the compositions, we freely use the strict commutativity of the squares in the diagram labelled with equal signs.

Post-whiskering the equality of the two strong transformations in condition \eqref{pentacond} 
in the definition of a comonoidal structure with $ (\bfK\tenrep_0(\bfL\tenrep_0(\bfM\tenrep_0 \bfN)))$, we see that  
\[\overline{\delta^{r}}\;\overline{\overline{\delta^{mm}}}\;\overline{\overline{\delta^{\ell}}} =\overline{\delta^{mr}}\;\overline{\overline{\delta^{m\ell}}}.\]
Thus, to show that the codomain and the domain of the pentagonator are equal, it suffices to show that
\begin{equation}\label{alphaeq}
  \overline{\overline{\alpha^{r}}}\;\overline{\overline{\alpha^{mm}}}\;\overline{\alpha^{\ell}}  = \overline{\overline{\alpha^{mr}}} \;\overline{\alpha^{m\ell}}
  \end{equation}
and that the domain of $\pithree_{\bfK, \bfL,\bfM,\bfN}$ is equal to \eqref{domordered}, while the codomain is equal to \eqref{codomordered}.

To verify \eqref{alphaeq}, we ignore the pre-whiskering with $ (( \Delta \bxone\mathbbm{1}_\Cc)\bxone\mathbbm{1}_\Cc)( \Delta \bxone\mathbbm{1}_\Cc)\Delta$ and verify that the strong transformations
\begin{multline}\label{onealphaside}
\big[\what \star (\mathbbm{1}_{\mathbbm{1}_{\A}} \bxtwo \alpha ) \star \big(
(\bfK\bxone((\bfL\bxone\bfM)\bxone \bfN)) \asstwo_{\Cc,\Cc\bxob \Cc,\Cc}
(\asstwo_{\Cc,\Cc,\Cc}\bxone\mathbbm{1}_\Cc) \big)\big] \\
\big[\alpha \star \big((\bfK\bxone(\bfL\tenrep_0\bfM)\bxone \bfN)(\asstwo_{\Cc,\Cc,\Cc}\bxone\mathbbm{1}_\Cc)\big) \big] 
\big[  \what \star (\alpha \bxtwo  \mathbbm{1}_{\mathbbm{1}_{\A}}) \star (((\bfK\bxone\bfL)\bxone\bfM)\bxone \bfN)\big] 
\end{multline}
and 
\begin{equation}\label{otheralphaside}
\big[\alpha \star\big(((\bfK\bxone\bfL)\bxone(\bfM\tenrep_0\bfN))  \asstwo_{\Cc\bxob\Cc,\Cc,\Cc} \big) \big] 
\big[   \alpha \star (((\bfK\tenrep_0\bfL)\bxone\bfM)\bxone \bfN)\big] 
\end{equation}
are equal. Note that these have domain and codomain the pseudofunctors from $((\Cc\bxob\Cc)\bxob\Cc)\bxob\Cc$ to $\A$ given by $((\bfK\tenrep_0\bfL)\tenrep_0\bfM)\tenrep_0 \bfN$ and \[(\bfK\tenrep_0(\bfL\tenrep_0(\bfM\tenrep_0 \bfN)))(\mathbbm{1}_\Cc\bxone\asstwo_{\Cc,\Cc,\Cc})\asstwo_{\Cc,\Cc\bxone\Cc,\Cc}(\asstwo_{\Cc,\Cc,\Cc}\bxone\mathbbm{1}_\Cc)=(\bfK\tenrep_0(\bfL\tenrep_0(\bfM\tenrep_0 \bfN)))\asstwo_{\Cc,\Cc,\Cc\bxone\Cc}\asstwo_{\Cc\bxone\Cc,\Cc,\Cc}, \]
where the equality follows from Proposition \ref{pentagtwo}.

The $1$-cell component of each, at an object $(((\ob{i}, \ob{j}),\ob{k}),\ob{l})$, is given by a functor from $((\bfK(\ob{i})\monob \bfL(\ob{j}))\monob\bfM(\ob{k}))\monob \bfN(\ob{l})$ to $\bfK(\ob{i})\monob(\bfL(\ob{j})\monob(\bfM(\ob{k})\monob \bfN(\ob{l})))$. The $1$-cell component of \eqref{onealphaside} acts on an 
object $\bigoplus_t(\bigoplus_s(\bigoplus_r(W_{rst},X_{rst}), Y_{st}), Z_t)$ by
\begin{equation*}\begin{split}
\bigoplus_t(\bigoplus_s(\bigoplus_r(W_{rst},X_{rst}), Y_{st}), Z_t) &\xrightarrow{\big[  \what \star (\alpha \bxtwo  \mathbbm{1}_{\mathbbm{1}_{\A}}) \star (((\bfK\bxone\bfL)\bxone\bfM)\bxone \bfN)\big]_{(((\ob{i}, \ob{j}),\ob{k}),\ob{l})}}  \\
\bigoplus_t(\bigoplus_s\bigoplus_r(W_{rst},(X_{rst}, Y_{st})), Z_t)&\xrightarrow{\big[\alpha \star \big((\bfK\bxone(\bfL\tenrep_0\bfM)\bxone \bfN)(\asstwo_{\Cc,\Cc,\Cc}\bxone\mathbbm{1}_\Cc)\big) \big]_{(((\ob{i}, \ob{j}),\ob{k}),\ob{l})}  }\\
\bigoplus_t\bigoplus_s\bigoplus_r(W_{rst},((X_{rst}, Y_{st}), Z_t))&\xrightarrow{\big[\what \star (\mathbbm{1}_{\mathbbm{1}_{\A}} \bxtwo \alpha ) \star \big(
(\bfK\bxone((\bfL\bxone\bfM)\bxone \bfN)) \asstwo_{\Cc,\Cc\bxob \Cc,\Cc}
(\asstwo_{\Cc,\Cc,\Cc}\bxone\mathbbm{1}_\Cc) \big)\big]_{(((\ob{i}, \ob{j}),\ob{k}),\ob{l})}   }\\
\bigoplus_t\bigoplus_s\bigoplus_r(W_{rst},(X_{rst},( Y_{st}, Z_t)))
\end{split}\end{equation*}
while the $1$-cell component of \eqref{otheralphaside} acts on the same object as
\begin{equation*}\begin{split}
\bigoplus_t(\bigoplus_s(\bigoplus_r(W_{rst},X_{rst}), Y_{st}), Z_t) &\xrightarrow{ \big[   \alpha \star (((\bfK\tenrep_0\bfL)\bxone\bfM)\bxone \bfN)\big]_{(((\ob{i}, \ob{j}),\ob{k}),\ob{l})} }\\
\bigoplus_t\bigoplus_s(\bigoplus_r(W_{rst},X_{rst}), (Y_{st}, Z_t)) &\xrightarrow{ \big[\alpha \star\big(((\bfK\bxone\bfL)\bxone(\bfM\tenrep_0\bfN))  \asstwo_{\Cc\bxob\Cc,\Cc,\Cc} \big) \big]_{(((\ob{i}, \ob{j}),\ob{k}),\ob{l})}   }\\
\bigoplus_t\bigoplus_s\bigoplus_r(W_{rst},(X_{rst}, (Y_{st}, Z_t)))
\end{split}\end{equation*}
so the two functors are the same resummations on objects and hence also on morphisms. As before, since the $2$-cell components of both strong transformations have the same source and target, and all involved components are just resummations, they also coincide and \eqref{alphaeq} holds.

To prove that the domain of $\pithree_{\bfK, \bfL,\bfM,\bfN}$ is equal to \eqref{domordered}, while the codomain is equal to \eqref{codomordered}, we first need the following lemma.

\begin{lemma}
For any $\ob{i}\in\Cc$ and $\j\in\{1,2,3\}$, we can write $\Delta(\delta_{\ob{i},j}) = (\delta_{\ob{i},j}^{(1)}, \delta_{\ob{i},j}^{(2)})$ with no outside direct sums.
\end{lemma}

\proof
We verify this by comparing the $1$-cell components of both sides of condition \eqref{pentacond} in the definition of a comonoidal structure. 
Suppose $\Delta(\delta_{\ob{i},j}) = \bigoplus_{s_j}(\delta_{\ob{i},j, s_j}^{(1)}, \delta_{\ob{i},j, s_j}^{(2)})$ for $j\in\{1,2,3\}$.
For the the left hand side, we compute
\begin{equation*}\begin{split}
\big(&(\mathbbm{1}_{\mathbbm{1}_{\Cc}}\bxtwo \delta)\star \Delta\big)_{\ob{i}} = (\mathbbm{1}_{\ob{i}_\ell}, ( \delta_{\ob{i}_r,1}  , ( \delta_{\ob{i}_r,2}   ,  \delta_{\ob{i}_r,3}  )))\\
\big( & ( (\mathbbm{1}_{\Cc}\bxone \asstwo_{\Cc,\Cc,\Cc}) (\mathbbm{1}_{\Cc}\bxone(\Delta \bxone \mathbbm{1}_{\Cc}))) \star \delta \big)_{\ob{i}} \\
&=( (\mathbbm{1}_{\Cc}\bxone \asstwo_{\Cc,\Cc,\Cc}) (\mathbbm{1}_{\Cc}\bxone(\Delta \bxone \mathbbm{1}_{\Cc}))) ( \delta_{\ob{i},1}  , ( \delta_{\ob{i},2}   ,  \delta_{\ob{i},3}  ))\\
&=( (\mathbbm{1}_{\Cc}\bxone \asstwo_{\Cc,\Cc,\Cc}) ( ( \delta_{\ob{i},1}  , ( \bigoplus_{s_2}(\delta_{\ob{i},j, s_2}^{(1)}, \delta_{\ob{i},j, s_2}^{(2)})  ,  \delta_{\ob{i},3}  ))\\
&= ( \delta_{\ob{i},1}  , \bigoplus_{s_2}(\delta_{\ob{i},j, s_2}^{(1)}, (\delta_{\ob{i},j, s_2}^{(2)}  ,  \delta_{\ob{i},3} ) ))\\
\big( &( (\mathbbm{1}_{\Cc}\bxone \asstwo_{\Cc,\Cc,\Cc})\asstwo_{\Cc,\Cc\bxob\Cc,\Cc}) \star (\delta \bxone\mathbbm{1}_{\mathbbm{1}_{\Cc}} )\star \Delta \big)_{\ob{i}}
 = (( \delta_{\ob{i}_\ell,1}  , ( \delta_{\ob{i}_\ell,2}   ,  (\delta_{\ob{i}_\ell,3}, \mathbbm{1}_{\ob{i}_r}  )))
\end{split}\end{equation*}
so the $1$-cell component of the left hand side is given by
\begin{equation*}\begin{split}
( \delta_{\ob{i},1} \delta_{\ob{i}_\ell,1}   , \bigoplus_{s_2}(\delta_{\ob{i}_r,1} \delta_{\ob{i},j, s_2}^{(1)}\delta_{\ob{i}_\ell,2} , (\delta_{\ob{i}_r,2} \delta_{\ob{i},j, s_2}^{(2)}\delta_{\ob{i}_\ell,3}  ,  \delta_{\ob{i}_r,3}  \delta_{\ob{i},3} ) )).
\end{split}\end{equation*}

For the right hand side, we compute
\begin{equation*}\begin{split}
\big( &(\mathbbm{1}_\Cc \bxone (\mathbbm{1}_\Cc \bxone \Delta))\star \delta\big)_{\ob{i}} = ( \delta_{\ob{i},1}  , ( \delta_{\ob{i},2}   ,  \bigoplus_{s_3}(\delta_{\ob{i},3, s_3}^{(1)}, \delta_{\ob{i},3, s_3}^{(2)}) ))\\
\big( & (\asstwo_{\Cc,\Cc,\Cc\bxob\Cc}( \Delta\bxone( \mathbbm{1}_\Cc \bxone\mathbbm{1}_\Cc ) )) \star \delta)_{\ob{i}}  =  
 \asstwo_{\Cc,\Cc,\Cc\bxob\Cc} ( \bigoplus_{s_1}(\delta_{\ob{i},j, s_1}^{(1)}, \delta_{\ob{i},j, s_1}^{(2)}) , ( \delta_{\ob{i},2}   ,  \delta_{\ob{i},3}  ))\\
 &= \bigoplus_{s_1}(\delta_{\ob{i},1, s_1}^{(1)},( \delta_{\ob{i},1, s_1}^{(2)} , ( \delta_{\ob{i},2}   ,  \delta_{\ob{i},3}  )))
\end{split}\end{equation*}
so the $1$-cell component of the right hand side is given by
\[ \bigoplus_{s_1}(\delta_{\ob{i},1} \delta_{\ob{i},1, s_1}^{(1)},(  \delta_{\ob{i},2} \delta_{\ob{i},1, s_1}^{(2)} , \bigoplus_{s_3}( \delta_{\ob{i},3, s_3}\delta_{\ob{i},2}   ,  \delta_{\ob{i},3, s_3}^{(2)}\delta_{\ob{i},3}  )))  \]

For both $1$-cells to be equal, the sums over $s_1,s_2,s_3$ must all be trivial and the claim follows.
\endproof

To compare the codomain $\overline{\delta^{mr}}\;\overline{\alpha^{mr}}\;\overline{\delta^{m\ell}}\;\overline{\alpha^{m\ell}}$ of $\pithree_{\bfK, \bfL,\bfM,\bfN}$ with $\overline{\delta^{mr}}\;\overline{\overline{\delta^{m\ell}}}\; \overline{\overline{\alpha^{mr}}} \;\overline{\alpha^{m\ell}}$, we need to swap the middle two $2$-cells given by $\delta$ and $\alpha$. By Lemma \ref{swaptwocells}, the component $2$-cell of the modification from $\overline{\alpha^{mr}}\;\overline{\delta^{m\ell}}$ to $\overline{\overline{\delta^{m\ell}}}\; \overline{\overline{\alpha^{mr}}}$ is given by 
\begin{equation*}\begin{split}
\alpha_{((\bfK\bxone\bfL)\bxone(\bfM\tenrep_0\bfN))(\Delta\bxone(\mathbbm{1}_\Cc \bxone\mathbbm{1}_\Cc))(\delta_{\ob{i}})}=\alpha_{ ((\bfK(\delta_{\ob{i},1}^{(1)}), \bfL(\delta_{\ob{i},2}^{(2)}))  , \bfM (\delta_{\ob{i},2} ) \monmor \bfN(\delta_{\ob{i},3} ) )}
\end{split}\end{equation*}
and considering the construction in Proposition \ref{alpha}, the lack of summations implies that this $2$-cell component is indeed the identity.

To compare the domain 
$\overline{\delta^{r}}\;\overline{\alpha^{r}}\;\overline{\delta^{mm}}\;\overline{\alpha^{mm}}\;\overline{\delta^{\ell}}\;\overline{\alpha^{\ell}}$ of $\pithree_{\bfK, \bfL,\bfM,\bfN}$ with $\overline{\delta^{r}}\;\overline{\overline{\delta^{mm}}}\;\overline{\overline{\delta^{\ell}}}\;\overline{\overline{\alpha^{r}}}\;\overline{\overline{\alpha^{mm}}}\;\overline{\alpha^{\ell}}$,
we need an intermediate step. Define
\[\underline{\alpha^r} =  \what \star (\mathbbm{1}_{\mathbbm{1}_{\A}} \bxtwo \alpha ) \star ((\bfK\bxone((\bfL\bxone\bfM)\bxone \bfN))\asstwo_{\Cc,\Cc\bxob\Cc,\Cc}((\mathbbm{1}_\Cc\bxone\Delta )\bxone \mathbbm{1}_\Cc)( \Delta\bxone \mathbbm{1}_\Cc)\Delta)\]
and
\[
\underline{\delta^\ell} = (\bfK\tenrep_0((\bfL\tenrep_0\bfM)\tenrep_0 \bfN))\asstwo_{\Cc,\Cc\bxob\Cc,\Cc} \star (\delta \bxtwo \mathbbm{1}_{\mathbbm{1}_{\Cc}}) \star \Delta
\]
we have modifications from $\overline{\alpha^{r}}\;\overline{\delta^{mm}}$ to $\overline{\overline{\delta^{mm}}}\;\underline{\alpha^r} $ and from $\overline{\alpha^{mm}}\;\overline{\delta^{\ell}}$ to $\underline{\delta^\ell}\;\overline{\overline{\alpha^{mm}}}$ whose $2$-cell components are given by
\[
\what((\mathbbm{1}_{\mathbbm{1}_{\A}} \bxtwo \alpha)_{(\bfK\bxone((\bfL\bxone\bfM)\bxone \bfN))(\mathbbm{1}_\Cc\bxone(\Delta \bxone \mathbbm{1}_\Cc))(\delta_{\ob{i}})}) = \mathbbm{1}_{\bfK(\ob{i}_\ell)}\otimes_2\alpha_{((\bfL(\delta_{\ob{i},2}^{(1)}),\bfM(\delta_{\ob{i},2}^{(2)})  ), \bfN(\delta_{\ob{i},3}) )}
\]
and
\[   
\alpha_{ ((\bfK\bxone(\bfL\tenrep_0\bfM))\bxone \bfN)(\delta_{\ob{i}_\ell}, \mathbbm{1}_{\ob{i}_r} )  } = \alpha_{((\bfK(\delta_{\ob{i}_\ell,1}), \bfL(\delta_{\ob{i}_\ell,2})\monmor\bfM(\delta_{\ob{i}_\ell,3})),\bfN(\mathbbm{1}_{\ob{i}_r}) )}
,\]
respectively. Again, given the lack of summations, the construction in Proposition \ref{alpha} shows that these are the identities, so 
\[
\overline{\delta^{r}}\;\overline{\alpha^{r}}\;\overline{\delta^{mm}}\;\overline{\alpha^{mm}}\;\overline{\delta^{\ell}}\;\overline{\alpha^{\ell}}
 = 
\overline{\delta^{r}}\;  \overline{\overline{\delta^{mm}}}\;\underline{\alpha^r} \;\underline{\delta^\ell}\;\overline{\overline{\alpha^{mm}}}\;\overline{\alpha^{\ell}}.
\]
Finally we have a modification from $\underline{\alpha^r} \;\underline{\delta^\ell}$ to $\overline{\overline{\delta^{\ell}}}\;\overline{\overline{\alpha^{r}}}$ whose $2$-cell components are given by
\[
\what((\mathbbm{1}_{\mathbbm{1}_{\A}} \bxtwo \alpha)_{(\bfK\bxone((\bfL\bxone\bfM)\bxone \bfN)) \asstwo_{\Cc,\Cc\bxob\Cc,\Cc} (\delta_{\ob{i}_\ell}, \mathbbm{1}_{\ob{i}_r} )    } =  \mathbbm{1}_{\bfK(\ob{i}_{\ell\ell})} \otimes_2 \alpha_{ ((\bfL(\delta_{\ob{i_\ell}_2}), \bfM(\delta_{\ob{i_\ell}_3})), \mathbbm{1}_{\bfN(\ob{i}_r)} )}
\]
which, similarly, is the identity, so we have indeed shown that 
\[
\overline{\delta^{r}}\;\overline{\alpha^{r}}\;\overline{\delta^{mm}}\;\overline{\alpha^{mm}}\;\overline{\delta^{\ell}}\;\overline{\alpha^{\ell}}
 = 
\overline{\delta^{r}}\;\overline{\overline{\delta^{mm}}}\;\overline{\overline{\delta^{\ell}}}\;\overline{\overline{\alpha^{r}}}\;\overline{\overline{\alpha^{mm}}}\;\overline{\alpha^{\ell}}
\]
as desired. Thus the pentagonator can be chosen as the identity.

\subsection{The main theorem}

\begin{theorem}
Let $\Cc\in \twoAk$ be equipped with a comonoidal structure. Then $\Cc\amod$ is a monoidal $2$-category with monoidal product $\tenrepC$ (Proposition \ref{tensorthree}), unit $\One_{\Cc\amod}$  (Proposition \ref{unitthree}), associator $\assthree$ (Proposition \ref{assthree}), left and right unitors
$\lunthree$ and $\runthree$ (Proposition \ref{unitorthree}), and identities for the pentagonator $\pithree$ and $2$-unitors $\muthree, \lambdathree$ and $\rhothree$.
\end{theorem}

\proof
It is easy to see that the diagrams (11.2.14), (11.2.16) and (11.2.17) in \cite[Def.11.2.4]{JY} commute trivially, since all higher coherences are given by identities.\endproof

\section{Examples}\label{sec:ex}

\subsection{Trivial on the level of $1$-morphisms and  $2$-morphisms}

Let $\mathsf{C}$ be a comonoid with comultiplication $\tilde\Delta$. Consider the $2$-category $\Cc$ whose objects are given by the elements of $\mathsf{C}$, and whose morphism categories $\Cc(\ob{i},\ob{j})$ for $\ob{i},\ob{j}\in \mathsf{C}$ are defined as follows. For $\ob{i}=\ob{j}$, we have $\Cc(\ob{i},\ob{i}) =  \Add \{\mathbbm{1}_{\ob{i}}\}$ with $\Cc(\ob{i},\ob{i})(\mathbbm{1}_{\ob{i}},\mathbbm{1}_{\ob{i}}) = \Bbbk\cdot \mathbbm{1}_{\mathbbm{1}_{\ob{i}}}$. For
 $\ob{i}\neq\ob{j}$, we have $\Cc(\ob{i},\ob{j}) =  \Add \{\F_{\ob{j},\ob{i}}\}$ with $\Cc(\ob{i},\ob{i})(\F_{\ob{j},\ob{i}},\F_{\ob{j},\ob{i}}) = \Bbbk\cdot \mathbbm{1}_{\mathbbm{1}_{\ob{i}}}$. Horizontal composition is defined by $\F_{\ob{k},\ob{j}}\F_{\ob{j},\ob{i}} =\F_{\ob{k},\ob{i}}$, in particular $\F_{\ob{j},\ob{i}}$ and $\F_{\ob{i},\ob{j}}$ are mutual inverses.
 
There is a $2$-functor $\Delta\colon \Cc\to\Cc\bxob\Cc$ uniquely defined by $\Delta(\ob{i}) = (\ob{i}_1,\ob{i}_2)$ if $\tilde\Delta(\ob{i}) = (\ob{i}_1,\ob{i}_2)$ and $\Delta(\F_{\ob{j},\ob{i}}) =( \F_{\ob{j}_1,\ob{i}_1},\F_{\ob{j}_2,\ob{i}_2})$. Moreover we can define a $2$-functor $\Theta\colon \Cc\to \Vv$ uniquely by sending each $\ob{i}$ to $\bullet$. Then $\delta$ is the identity and $\theta^\ell,\theta^{\ell *}, \theta^r, \theta^{r*}$ are given by their $1$-cell components $\theta_{\ob{i}}^\ell = \F_{\ob{i},\ob{i}_2}$, $\theta_{\ob{i}}^{\ell *} = \F_{\ob{i}_2,\ob{i}}$,  $\theta_{\ob{i}}^r = \F_{\ob{i},\ob{i}_1}$ and $\theta_{\ob{i}}^{r*} = \F_{\ob{i}_1,\ob{i}}$ with their $2$-cell components simply being identities.

Any $\bfM\in \Cc\amod$ is given by isomorphic categories $\bfM(\ob{i})$ for $\ob{i}\in \mathsf{C}$ together with isomorphisms $\bfM(\F_{\ob{j},\ob{i}} )\colon \bfM(\ob{i})\to \bfM(\ob{j})$ such that $\bfM(\F_{\ob{i},\ob{j}})\bfM(\F_{\ob{j},\ob{i}}) = \mathbbm{1}_{\bfM(\ob{i})}$ and $\bfM(\F_{\ob{k},\ob{i}}) = \bfM(\F_{\ob{k},\ob{j}})\bfM(\F_{\ob{j},\ob{i}})$. For two such $2$-representations $\bfM, \bfN$, we then simply have
$(\bfM\tenrepC_0\bfN)(\ob{i}) = \bfM(\ob{i}_1)\monob \bfN(\ob{i}_2)$ and 
$(\bfM\tenrepC_0\bfN)(\F_{\ob{j},\ob{i}})= \bfM(\F_{\ob{j}_1,\ob{i}_1})\monmor \bfN(\F_{\ob{j}_2,\ob{i}_2})$.

\subsection{Trivial on the level of objects and $2$-morphisms}\label{1level}

Let $A$ be a bialgebra with a fixed basis $\{a_i\,\vert\, i\in I\}$ containing $a_1 = 1_A$ with respect to which the structure constants are positive integral for the multiplication, comultiplication and the counit. Define a 2-category $\Cc$ with one object $\ast$ and morphism category $\Cc(\ast,\ast) = \Add\{ \F_i\,\vert\, i\in I\}$ with $\Cc(\ast,\ast)(\F_i,\F_i) = \Bbbk\cdot \mathbbm{1}_{\F_i}$ and $\Cc(\ast,\ast)(\F_i,\F_j) = 0$ for $i\neq j$. Horizontal composition is defined by $\F_i\F_j = \bigoplus_{k\in I} \F_k^{\oplus \lambda_{i,j}^k}$ if $a_ia_j = \sum_{k\in I} \lambda_{i,j}^k a_k $. 

Define $\Delta\colon \Cc\to \Cc\bxob \Cc$ uniquely by $\Delta(\F_i) = \bigoplus_{j,k\in I} (\F_j,\F_k)^{\oplus \mu_i^{j,k}}$ for the structure constants of the comultiplication. 
Moreover, we can define a $2$-functor $\Theta\colon \Cc\to \Vv$ by setting $\Theta(\F_i) = \mathbbm{1}_\bullet^{\oplus t_i}$ if the counit on $A$ sends $a_i$ to $t_i$.

Then all of $\theta^\ell,\theta^{\ell *}, \theta^r$ and $\theta^{r*}$ are the identity, and $\delta$ has the identity as a $1$-cell component and a canonical resummation as a $2$-cell component. Commutativity of the necessary diagrams follows from the usual argument about all morphisms being canonical resummations with the same source and target and hence being equal.

An explicit example is given by $A = \Bbbk[x]$ with $x$ primitive and annihilated by the counit. In this case, writing $\F_i$ for $\F_{x^i}$, we have $\Cc(\ast,\ast) = \Add\{ \F_i\,\vert\, i\in \mathbb{Z}_\geq 0\}$ with $\F_0 = \mathbbm{1}_{\ast}$ and $\F_i\F_j = \F_{i+j}$. The $\Delta(\F_i) = \bigoplus_{j=0}^i (\F_j, \F_{i-j})$ and  $\Theta(\F_i)=0$ unless $i=0$.  

\subsection{$G$-graded vector spaces}

A special case of a $2$-category as in Section \ref{1level} is given by the monoidal category of $G$-graded vector spaces for a finite group $G$. In this case $A = \Bbbk G$ is the group algebra of a finite group $G$ with comultiplication $g\mapsto g\otimes g$ and counit $g\mapsto 1$. Then $\Cc$ is defined as having one object $\ast$ and $\Cc(\ast,\ast) = \Add\{ \F_g\,\vert\, g\in G\}$ with $\Cc(\ast,\ast)(\F_g,\F_g) = \Bbbk\cdot \mathbbm{1}_{\F_g}$ and $\Cc(\ast,\ast)(\F_g,\F_h) = 0$ for $g\neq h$. Horizontal composition is given by $\F_g\F_h=\F_{gh}$, in particular, $\F_1 = \mathbbm{1}_\ast$. 
Moreover, $\Delta(\F_g) = (\F_g,\F_g)$ and $\Theta(\F_g) = \mathbbm{1}_\bullet$.

Simple $2$-representations $\bfM_{H,\omega}$ are classified by pairs $(H,\omega)$ of a subgroup $H$ and a $2$-cocycle on $H$ up to $G$-conjugacy, see e.g. \cite[Example 7.4.10]{EGNO}. More precisely, $\bfM_{H,\omega}$ is the category of (necessarily projective) right modules in $\Cc$ for the algebra $1$-morphism $\mathrm{A}_{H, \omega}= \bigoplus_{h\in H} \F_h$ with multiplication induced by $\F_h\F_{h'}\xrightarrow{\omega_{h,h'}}\F_{hh'}$. The isomorphism classes of indecomposable $1$-morphisms are parameterised by the cosets of $G/H$.

Consider $\bfM_{H,\omega}\tenrepC_0 \bfM_{K,\tau}$. It is well-known that $G/H\times G/K$ decomposes, under the diagonal action into the disjoint union of orbits $G\cdot (eH, sK)$ where s runs through a set $\Gamma$ of representatives of $H\backslash G/K$. The stabiliser of $(eH, sK)$ consists of $H\cap sKs^{-1}$. Thus we obtain an equivalence of $2$-representations
\[\bfM_{H,\omega}\tenrepC_0 \bfM_{K,\tau}\cong \bigoplus_{g\in \Gamma} \bfM_{H\cap gKg^{-1},\sigma_g}\]
for some $2$-cocycles $\sigma_g$. By explicitly computing the algebra structure on the internal hom of a generator $(\mathrm{A}_{H, \omega}, \F_g\mathrm{A}_{K, \tau})$ of the component $\bfM_{H\cap gKg^{-1},\sigma_g}$ via the evaluation map, one can see that $\sigma_g = \omega \tau^g$ where $\tau^g$ denotes the corresponding conjugate of $\tau$, i.e. for $s,s'\in H\cap gKg^{-1}$, we have $\omega \tau^g(s,s')= \omega(s,s')\tau(g^{-1}sg)$.

\subsection{Trivial on the level of objects and $1$-morphisms}

Let $\Cc$ be a $2$-category with one object $\ast$, morphism category $\Cc(\ast,\ast) = \Add\{\mathbbm{1}_{\ast}\}$ and $\Cc(\ast,\ast)(\mathbbm{1}_{\ast},\mathbbm{1}_{\ast}) = A$ where $A$ is a bialgebra.

Define $\Delta\colon \Cc\to\Cc\bxob\Cc$ and $\Theta\colon \Cc\to \Vv$, which necessarily send $\ast$ to $(\ast,\ast)$ resp.\ $\bullet$, as well as $\mathbbm{1}_{\ast}$ to $(\mathbbm{1}_{\ast},\mathbbm{1}_{\ast})$ resp.\ $\mathbbm{1}_{\bullet}$, and which acts by the comultiplication, resp.\ the counit on the level of $2$-morphisms. 

Again, all of $\delta, \theta^\ell,\theta^{\ell *}, \theta^r$ and $\theta^{r*}$ are the identity.

Note that a $2$-representation of $\Cc$ is simply given by a category $\bfM(\ast)\in\A$ together with an algebra morphism from $A$ to the endomorphism algebra of $\mathbbm{1}_{\bfM(\ast)}$. If  $\bfM_1(\ast)= B_1\text{-}\mathrm{proj}$ and $\bfM_2(\ast)= B_1\text{-}\mathrm{proj}$, the categories of finitely generated projective modules over associative $\Bbbk$-algebras $B_1, B_2$, this corresponds to morphisms $\phi_i\colon A\to Z(B_i)$ for $i=1,2$, where $Z(B_i)$ denotes the centre of $B_i$. In this case, $(\bfM_1\tenrepC_0\bfM_2)(\ast)\cong (B_1\otimes_\Bbbk B_2)\text{-}\mathrm{proj}$ and the algebra morphism is given by $$A\to Z(B_1\otimes_\Bbbk B_2)=Z(B_1)\otimes_\Bbbk Z(B_2), a\mapsto \phi_1(a)\otimes \phi_2(a).$$

The example in this subsection can be combined with that in the previous one to define a comonoidal structure on a $Gr$-category $(G,A,\omega)$ provided $A$ is a bialgebra (see \cite[Section 2.11]{EGNO} for notation).

\end{document}